\documentclass[twocolumn]{svjour3}
\usepackage{amsmath,amssymb,bm}
\usepackage{graphicx,color,url}
\usepackage{mathptmx}

\usepackage{microtype}
\usepackage{xspace}

\newcommand{\rt}[1]{{#1}}

\newcounter{myremark}
\newenvironment{myremark}[1][]{\refstepcounter{myremark}\par\medskip
   \noindent \textbf{Remark~\themyremark. #1} \rmfamily}{}

\graphicspath{{pics/},{figures/}} 
\newcommand{\interface}{\Gamma}
\newcommand{\targetCurve}{\tilde{\interface}}
\newcommand{\conformalLayer}{\Omega^c}
\newcommand{\backDomainIn}{\Omega^{\textnormal{in}}} 

\begin{document}

\title{Immersed Boundary-Conformal Isogeometric Method for Linear Elliptic Problems}

\author{Xiaodong Wei \and Benjamin Marussig \and Pablo Antolin \and Annalisa Buffa}

\institute{X. Wei \and P. Antolin \at Institute of Mathematics, \'Ecole Polytechnique F\'ed\'erale de Lausanne, 1015 Lausanne, Switzerland \\
\and
A. Buffa \at Institute of Mathematics, \'Ecole Polytechnique F\'ed\'erale de Lausanne, 1015 Lausanne, Switzerland 
\at Istituto di Matematica Applicata e Tecnologie Informatiche ``Enrico Magenes'' del CNR, Pavia, Italy\\
\and
B. Marussig \at Institute of Applied Mechanics, Graz University of Technology, 8010 Graz, Austria
}

\maketitle

\begin{abstract}

We present a novel isogeometric method, namely the \emph{Immersed Boundary-Conformal Method} (IBCM), that features a layer of discretization conformal to the boundary while employing a simple background mesh for the remaining domain. In this manner, we leverage the geometric flexibility of the immersed boundary method with the advantages of a conformal discretization, such as intuitive control of mesh resolution around the boundary, higher accuracy per degree of freedom, automatic satisfaction of interface kinematic conditions, and the ability to strongly impose Dirichlet boundary conditions. In the proposed method, starting with a boundary representation of a geometric model, we extrude it to obtain a corresponding conformal layer. Next, a given background B-spline mesh is cut with the conformal layer, leading to two disconnected regions: an exterior region and an interior region. Depending on the problem of interest, one of the two regions is selected to be coupled with the conformal layer through Nitsche's method. Such a construction involves Boolean operations such as difference and union, which therefore require proper stabilization to deal with arbitrarily cut elements. In this regard, we follow our precedent work called the minimal stabilization method~\cite{ref:wei19u}. In the end, we solve several 2D benchmark problems to demonstrate improved accuracy and expected convergence with IBCM. Two applications that involve complex geometries are also studied to show the potential of IBCM, including a spanner model and a fiber-reinforced composite model. \rt{Moreover, we demonstrate the effectiveness of IBCM in an application that exhibits boundary-layer phenomena.}

\end{abstract}

\keywords{Immersed Method \and Conformal Boundary/Interface \and Boolean operations \and Stabilized Method \and Isogeometric Analysis \and Boundary Layer}


\section{Introduction}

%

Isogeometric analysis (IGA) was proposed to tightly integrate computer-aided design (CAD) with downstream engineering simulation \cite{ref:hughes05,ref:cottrell09}, which is enabled by employing the same basis of CAD in simulation. While significant advances have been made in almost every aspect of IGA, dealing with boundary representation (B-rep) remains one of the biggest challenges in the field. \rt{In many current commercial CAD systems}, solid modeling is built upon B-reps, where a solid model is represented by its boundary only. However, B-reps are not immediately ready to be used in IGA. A key issue is the lack of a true volumetric description. They need to be processed through either volume parameterization (from the geometry point of view) or the immersed method (from the analysis point of view). Volume parameterization conformal to a given B-rep is appealing but remains an open problem at large \cite{ref:martin09,ref:zhang13}. Alternatively, the immersed boundary method, which embeds a given B-rep into a simple background mesh (e.g., a Cartesian grid), eliminates the need for conformal and quality volume parameterization and \rt{has gained} wide popularity in IGA~ \cite{ref:rank12,ref:ruberg12,ref:schillinger12,ref:kamensky15,ref:hoang19,ref:casquero20im}.

The term ``immersed boundary method" indicates a large family of methods and is sometimes used interchangeably with other terms such as ``non-boundary-fitted method", ``finite cell method", ``fictitious domain method" and ``embedded domain method". Despite their various origins, all of such methods share the same idea: embedding a B-rep into a simple background mesh. We use the term only to reflect this key idea. Interested readers may refer to \cite{ref:zli06} and references therein for a comprehensive review.

In the context of IGA, the immersed idea is also closely related to Boolean operations (i.e., intersection, difference, and union) and the related trimming operations, which allow the specification of arbitrary interfaces within the elements of a tensor product surface. Indeed, an immersed method needs to address the same challenging problems encountered in dealing with Boolean operations, such as numerical integration of cut elements \cite{ref:kim09,ref:nagy15,ref:kudela15,ref:antolin19a,ref:massarwi19}, stabilization for small cut elements \cite{ref:marussig17,ref:marussig18,ref:guo18,ref:puppi19}, and imposition of boundary/interface conditions \cite{ref:bazilevs07,ref:embar10}.

Motivated to accurately capture boundary/interface features while retaining geometric flexibility of immersed methods, we propose a novel method that combines the immersed idea with a conformal discretization. We call it the \emph{immersed boundary-conformal method} (IBCM). The geometric construction of IBCM is conceptually simple. \rt{It starts with a B-rep of the domain of interest. The B-rep is first extruded to yield a boundary layer that is clearly conformal to the input B-rep}, so we call it the \emph{conformal layer}. Next, we embed the conformal layer into a sufficiently large background mesh such as a Cartesian grid. As a result, the background mesh is cut into two disconnected regions: an exterior region and an interior region. Depending on the problem of interest, one of the two regions is selected and coupled with the conformal layer to constitute (a part of) underlying computation domain. \rt{It is noted that we focus on 2D planar domains in this paper.}

An IBCM construction involves Boolean operations such as difference and union. Properly dealing with them in analysis is the key to the success of IBCM, which includes numerical integration on cut elements and interfaces, and stabilized formulation regardless of how elements are cut. To address the integration issue, we present an improved decomposition method to minimize the number of resulting quadrature cells. On the other hand, we adopt our previous work, the minimal stabilization method~\cite{ref:wei19u,ref:puppi19}, to guarantee stability for its simplicity and effectiveness.

Two kinds of applications are studied with IBCM: (1) accuracy enhancement via adding conformal layers to geometric features (e.g., holes, reentrant corners, etc), and (2) stress analysis on domains composed of multiple materials (e.g., inclusions). It is straightforward for IBCM to control the mesh resolution near boundaries/interfaces due to the conformal layers. As singularities and/or large gradients often happen close to the boundary, it is often vital to have a sufficiently fine mesh near the boundary to obtain accurate solutions. The conformal layer also provides the possibility to strongly impose Dirichlet boundary conditions, which, in contrast, is not possible in immersed methods. In the case of modeling composite materials, material interfaces can be represented with conformal meshes in IBCM, and thus kinematic constraints are automatically satisfied on the interfaces by construction. \rt{While strong imposition of boundary/interface conditions is enabled by the conformal layer, coupling it with the background mesh is done weakly with the control over where to put the coupling interface.} In addition to all these analysis benefits, we will further show that IBCM indeed can easily model complex geometries via the immersed idea. 

It is worth mentioning that there exist several methods \rt{that IBCM shares the essential idea with. First, a popular group of methods in computational fluid dynamics, called overset grid methods \cite{ref:benek83,ref:liou94,ref:henshaw94,ref:tang03,ref:kannan07,ref:duan20}, also employ boundary-fitted meshes near geometric features while keeping discretization simple elsewhere. However, the two methods differ in the key building blocks such as geometric representations and treatment of coupling interfaces. Overset grid methods are typically developed on linear elements, and different meshes are coupled via, for instances, hybrid meshes \cite{ref:liou94}, and interpolation schemes \cite{ref:benek83,ref:tang03} that solve the problem iteratively to reach a convergent solution.} Second, in \cite{ref:bouclier16}, disk-shaped patches (represented by degenerated NURBS) are added on top of a large domain (represented by a Cartesian grid) to model holes or inclusions. Different domains are coupled in a non-intrusive way that requires tens or even hundreds of iterations to obtain a converged solution for the linear elasticity problem. Moreover, the method primarily aims to capture local features inside a rectangle-shaped domain and thus does not support complex exterior geometries. Third, a multi-mesh finite element method was proposed in \cite{ref:dokken19} to optimize the shapes of disk- or ring-like objects that are embedded in a large background domain. \rt{Compared to the multi-mesh approach~\cite{ref:becker11,ref:dokken19,ref:johansson19}, the main differences of our proposed method lie in four aspects: (1) smooth discretization with IGA, (2) emphasis on parameterizations of conformal layers, (3) a minimal stabilization approach, and (4) use of conformal layers to improve accuracy and provide flexible control near the boundary/interface.}

The paper is organized as follows. We first introduce several basic concepts and notations in Section \ref{sec:notation}. Section \ref{sec:overview} introduces how to deal with Boolean operations in IGA, which lays the foundation to our proposed method. Section \ref{sec:IBCM} presents the proposed immersed boundary-conformal method. In Section \ref{sec:result}, we present several numerical examples to demonstrate the proposed method in terms of improved accuracy, ease of imposing Dirichlet boundary conditions, and flexibility in modeling complex geometries. Section \ref{sec:con} concludes the paper and suggests future work.

\section{Analysis-Aware Treatment of Boolean Operations}
\label{sec:overview}

In this section, we introduce how to handle difference and union operations in IGA. Their analysis-aware treatment lays the foundation to our proposed method. As the difference operation is built upon trimming, we discuss trimming in the following without loss of generality. Our description is based on 2D planar domains unless stated otherwise.

\subsection{B-spline/NURBS surfaces}
\label{sec:notation}

We first introduce several basic notations of B-spline/NURBS surfaces to facilitate our discussion.  Given a control mesh $\{\bm{P}_{ij}\}_{i=1,j=1}^{n,m}$ ($\bm{P}_{ij}\in\mathbb{R}^2$), degrees $p$ and $q$, and knot vectors $\{\xi_i\}_{i=1}^{n+p+1}$ and $\{\eta_j\}_{j=1}^{m+q+1}$, a bivariate B-spline basis function is defined as
\begin{equation}
B_{\bm{i},\bm{p}}(\bm{\xi}) = N_{i,p} (\xi) M_{j,q}(\eta),
\end{equation}
where $\bm{i}=(i,j)$, $\bm{p}=(p,q)$, $\bm{\xi} = (\xi,\eta)$, and $N_{i,p}(\xi)$ and $M_{j,q}(\eta)$ are univariate B-spline basis functions defined on $\{\xi_i\}_{i=1}^{n+p+1}$ and $\{\eta_j\}_{j=1}^{m+q+1}$, respectively. The geometric mapping of a B-spline surface, $\bm{F}:\,\hat{\Omega}\to\Omega$, is defined as
\begin{equation}
\bm{F}(\bm{\xi}) = \sum_{\bm{i} \in \mathcal{I}} \bm{P}_{\bm{i}} B_{\bm{i},\bm{p}} (\bm{\xi}), \quad \mathcal{I}=\{(i,j)\}_{i=1,j=1}^{n,m},
\label{eq:bsp_surf}
\end{equation}
where $\hat{\Omega}:=(\xi_1,\xi_{n+p+1})\times (\eta_1,\eta_{m+q+1})$ and $\Omega$ are called the \emph{parametric domain} and the \emph{physical domain}, respectively. \rt{We also call a surface a \emph{patch}, especially when multiple surfaces are involved.} $\hat{\Omega}$ is partitioned by the knot vectors, leading to a parametric mesh $\hat{\mathcal{M}}$,
\begin{equation}
\begin{aligned}
\hat{\mathcal{M}} = \{ & (\xi_i,\xi_{i+1}) \times (\eta_j,\eta_{j+1}):\quad \xi_i < \xi_{i+1} \text{ and } \eta_j < \eta_{j+1}, \\
& i=1,\ldots,n+p, \quad j=1,\ldots,m+q \}.
\end{aligned}
\end{equation}
The corresponding physical mesh $\mathcal{M}$ is obtained as the images of elements in $\hat{\mathcal{M}}$ through the geometric mapping,
\begin{equation}
\mathcal{M} = \{\bm{F}(K):\, K\in \hat{\mathcal{M}} \}.
\end{equation}
The element boundaries in $\hat{\mathcal{M}}$ are represented by knot lines, i.e., $\{\xi_i\}\times [\eta_1,\eta_{m+q+1}]$ and $\{\eta_j\}\times [\xi_1,\xi_{n+p+1}]$. Their images are called knot-line curves and they are the element boundaries in $\mathcal{M}$.

NURBS are a generalization of B-splines and they are able to exactly represent conic sections such as circles and ellipses. A weight $w_{\bm{i}} \in \mathbb{R}$ is associated to each control point in addition to its coordinates. The NURBS basis function corresponding to the B-spline $B_{\bm{i},\bm{p}}(\bm{\xi})$ is defined as
\begin{equation}
R_{\bm{i},\bm{p}}(\bm{\xi}) = \frac{ B_{\bm{i},\bm{p}}(\bm{\xi}) w_{\bm{i}}}{\sum_{\bm{j} \in \mathcal{I}} B_{\bm{j},\bm{p}}(\bm{\xi}) w_{\bm{j}} }.
\end{equation}
Note that NURBS and B-spline basis functions are equivalent when $w_{\bm{i}} =1$ $\forall \bm{i}\in\mathcal{I}$, i.e., $R_{\bm{i},\bm{p}} = B_{\bm{i},\bm{p}}(\bm{\xi})$. A NURBS surface representation is obtained by replacing $B_{\bm{i},\bm{p}}(\bm{\xi})$ with $R_{\bm{i},\bm{p}}(\bm{\xi})$ in Eq. \ref{eq:bsp_surf}. To simplify notations, we use $B_{\bm{i},\bm{p}}(\bm{\xi})$ to generically denote a B-spline or NURBS basis function in the rest of the paper.

The spline space on the physical domain $\Omega$ is defined as
\begin{equation}
\mathcal{B} = \mathrm{span} \{ B_{\bm{i},\bm{p}} \circ \bm{F}^{-1}(\bm{x}) : \, \bm{i} \in \mathcal{I} \}.
\end{equation}
Let $\Omega_a \subset \Omega$ be a subdomain. We further define the restriction of $\mathcal{B}$ to $\Omega_a$ as
\begin{equation}
\begin{aligned}
\mathcal{B}|_{\Omega_a} = \mathrm{span} \{ & B_{\bm{i},\bm{p}} \circ \bm{F}^{-1}(\bm{x}) : \, \bm{i} \in \mathcal{I},\\
& \mathrm{supp}B_{\bm{i},\bm{p}} \cap \bm{F}^{-1}(\Omega_a) \neq \emptyset \},
\end{aligned}
\end{equation}
where $\mathrm{supp}B_{\bm{i},\bm{p}}:=(\xi_i,\xi_{i+p+1})\times (\eta_j,\eta_{j+q+1})$ is the support of $B_{\bm{i},\bm{p}}(\bm{\xi})$. $\mathcal{B}|_{\Omega_a}$ will be used in trimmed geometries.

\subsection{Trimming}
\label{sec:trim}

\begin{figure*}[htb]
\centering
\includegraphics[width=0.7\textwidth]{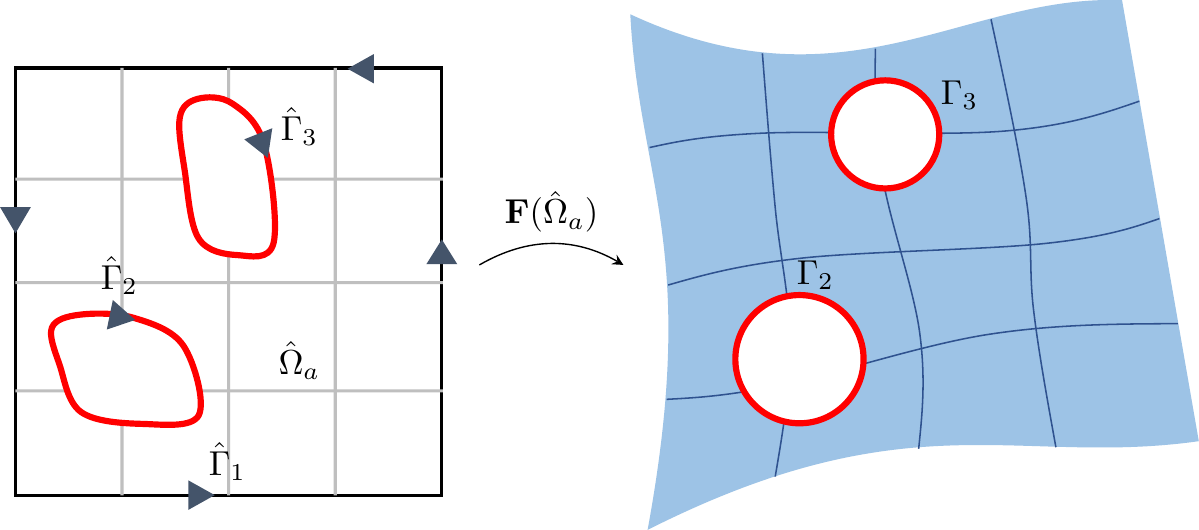}
\caption{Illustration of trimming, where $\Gamma_i$ ($i=2,3$) are trimming loops in the physical domain.}
\label{fig:trim}
\end{figure*}

From the geometric point of view, trimming provides a powerful and flexible tool to represent complex shapes that do not have a tensor-product structure. However, it is merely a visualization means to make part of a geometry invisible to users without changing anything else. More specifically, we consider a NURBS surface $\bm{F}:\ \hat{\Omega} \to \Omega$. Applying trimming is equivalent to restricting $\bm F$ to a visible or \emph{active} subdomain $\hat{\Omega}_a\subset \hat\Omega$. Correspondingly, the complementary portion $\bm{F}(\hat{\Omega}\backslash\hat{\Omega}_a)$ is invisible or trimmed away. We note that the only change in a trimming operation is the change in the definition domain: from $\hat{\Omega}$ to $\hat{\Omega}_a$. The trimming operation is illustrated in Fig.\ \ref{fig:trim}. We can observe its resemblance to the immersed idea, where neither of their background meshes aligns with the domain boundary.

The active subdomain $\hat{\Omega}_a$ is identified by a set of closed loops $\{\hat{\Gamma}_i\}$. Each trimming loop $\hat{\Gamma}_i$ is generally a collection of connected B-spline curves. Trimming loops are oriented such that all points in $\hat{\Omega}_a$ locate on the left side (in convention) of any $\hat{\Gamma}_i$. In practice, trimming curves are often given in the physical domain and may not even lie on the surface $\mathbf{F}(\hat\Omega)$. To obtain the parametric counterpart $\hat{\Gamma}_i$, points are sampled along $\Gamma_i$ and then pulled back to $\hat{\Omega}$ through the inversion algorithm \cite{ref:nurbsbook}. $\hat{\Gamma}_i$ is constructed based on these sampled points.

From the analysis point of view, numerical integration is performed on the active subdomain $\hat{\Omega}_a$ rather than the entire domain $\hat{\Omega}$. $\hat{\Omega}_a$ consists of non-trimmed elements and trimmed elements. Given an element $K\in\hat{\mathcal{M}}$, $K$ is cut if $0<|K\cap\hat{\Omega}_a|<|K|$, where $|\cdot|$ denotes area. Its active portion, $K\cap\hat{\Omega}_a$, is referred to as a trimmed element. Clearly, non-trimmed elements are those with $|K\cap\hat{\Omega}_a|=|K|$, whereas elements with no active area (i.e., $|K\cap\hat{\Omega}_a|=0$) are inactive in the sense that they are not used in geometric representation or in analysis. 

While it is straightforward to integrate non-trimmed elements, the challenge lies in the integration for a trimmed element because it is topologically a polygon with certain edges corresponding to part of a (high-order) trimming curve, where no quadrature rules are immediately available. One solution is to reparameterize trimmed elements such that the standard Gauss quadrature rule can be applied. Reparameterization of a trimmed element is a \emph{local} operation that has no influence on other elements. It consists of an approximation step and a decomposition step. The former approximates each involved trimming curve with a B\'{e}zier curve, whereas the latter decomposes trimmed element into quadrilateral/triangular cells. We explain related details in the following. Interested readers may refer to \cite{ref:marussig18r} for a comprehensive review on the treatment of trimming in IGA.

\textbf{Approximation step}. An element $K\in\hat{\mathcal{M}}$ may be cut by multiple trimming loops. For each involved trimming loop $\hat{\Gamma}_i$, we focus on the portion restricted to $K$, i.e., $\hat{\Gamma}_i \cap K$. In general, $\hat{\Gamma}_i \cap K$ may be composed of several B-spline curves that are piecewise polynomials, but Gauss quadrature can only be applied to polynomials. 

Therefore, we approximate $\hat{\Gamma}_i \cap K$ with a B\'{e}zier curve that has a polynomial representation. The degree of the target B\'{e}zier curve is chosen to be that of spline discretization even if the trimming curve has a higher degree. It has been proven that this choice guarantees optimal error estimates \cite{ref:antolin19a}. The remaining procedure of least-squares fitting the B\'{e}zier curve to the sampled points of $\hat{\Gamma}_i \cap K$ is straightforward.

Note that the above discussion is based on the assumption that \rt{the trimming curve $\hat{\Gamma}_i \cap K$ is sufficiently smooth. However, generally it involves points of reduced continuity or even $C^0$ points. In practice, we observe that to retain integration accuracy, it suffices to only respect those $C^0$ points, for which we split $\hat{\Gamma}_i \cap K$ at each of them and fit a B\'{e}zier curve separately to each resulting piece.}

\begin{figure}[htb]
\centering
\includegraphics[width=\columnwidth]{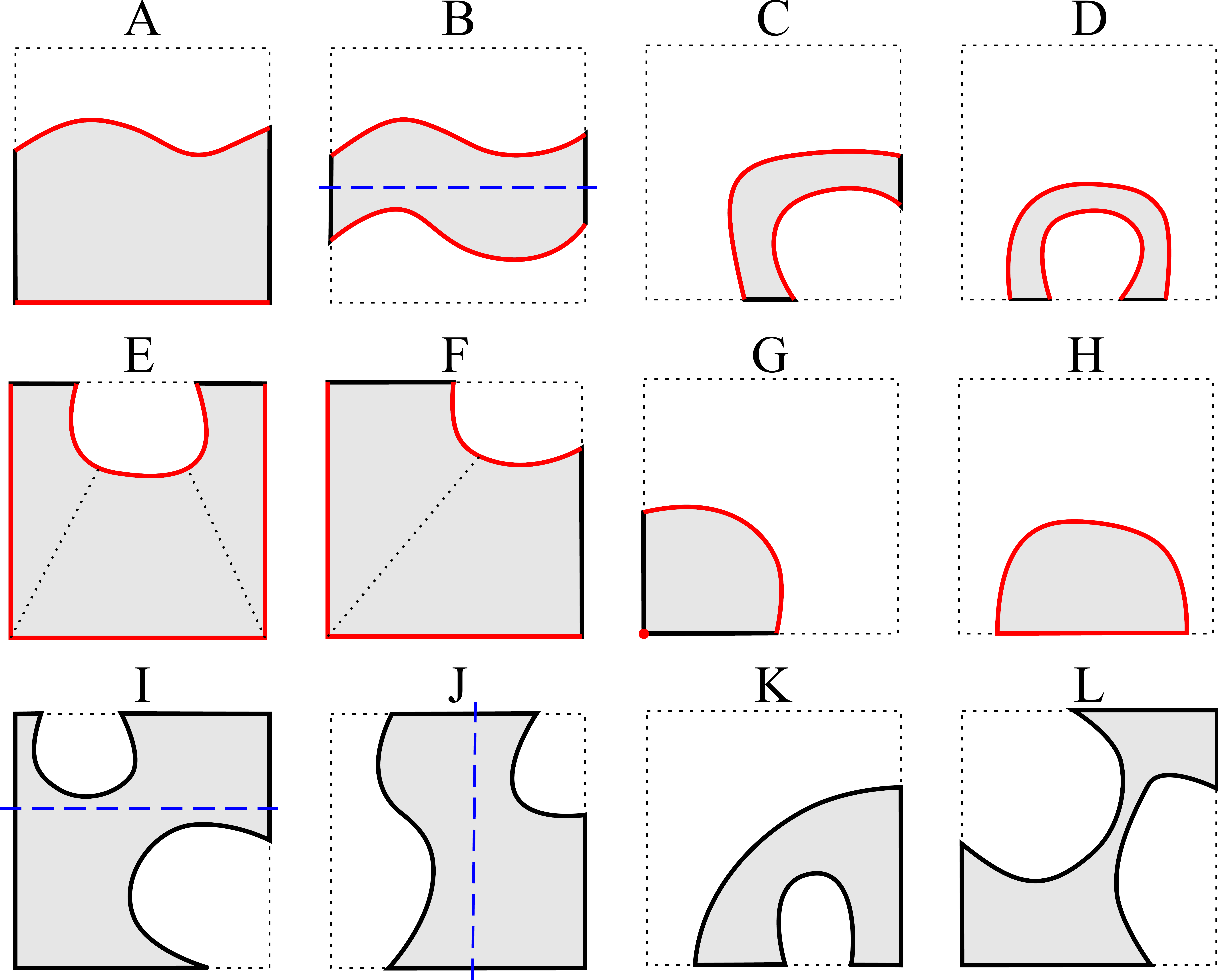}
\caption{Reparameterization of trimmed elements using \rt{quadrilateral} cells.}
\label{fig:trimming_cases}
\end{figure}

\textbf{Decomposition step}. We develop a divide-and-conquer approach to decompose trimmed elements into reparameterization cells.
\rt{The method applies to 2D only.
We first investigate how an element (in the parametric domain) is cut by one or two trimming curves.
The key is to find a pair of curves and create one or multiple ruled surfaces between them.
In each pair of curves, one has to be the trimming curve, whereas the other can be a straight segment (A, H), a series of concatenated segments (E, F), a point (K), or another trimming curve (B, C, D); see Fig. \ref{fig:trimming_cases}.
As a result, 8 base cases are identified.
Ruled surfaces are then created between a pair of such curves.}
Note that cases G and H present a single red curve and lead to degenerated quadrilateral reparameterization cells: For the case G, the ruled surface is created using the red curve and the \rt{opposite} red vertex; whereas in the case H, the cell is created by splitting the red curve at its two vertices, and generating a ruled surface between them. As already mentioned above, the curved segments of the trimming loop may present $C^0$ features. Using this strategy, these features will be explicitly represented in the generated reparameterization.
\rt{On the other hand, applying the same idea to finding base cases in 3D becomes intractable due to the proliferation of possible cases. One viable solution is to explore the unique base cases of marching cubes \cite{ref:mc} and design a corresponding decomposition for each of them. 
Alternatively, each cut element may be divided into a set of tetrahedra, which, however, leads to a proliferation of integration cells \cite{ref:fries17a,ref:fries17b}.
}

If the algorithm above fails (the generated ruled surface presents a change in its Jacobian's sign), or the trimmed element does not fall into one of the base categories A-H, a divide-and-conquer strategy is applied: The trimming loop is split with a line parallel to one of the two coordinate axes, and the algorithm is applied again for the resulting pieces. This is also the case in which the element presents nested trimming loops. \rt{Such a recursive method works up to a certain geometric tolerance, beyond which the splitting algorithm terminates and unresolved trimming features will be ignored. Adoption of geometric tolerances is a common practice in CAD when implementing trimming-related operations.}

For selecting the splitting line we consider the bounding boxes of all the curved segments on the trimming loop, cherry pick the one with the smallest volume, and consider a line that passes through that bounding box's mid point and is perpendicular to its longest coordinate direction.
Nevertheless, whenever possible, the splitting line is defined such that different curved segments are not intersected by the line and they are separated into different reparameterization cells (see, e.g., cases B, I, and J in Fig.\ \ref{fig:trimming_cases}).
This is not always possible; see for instance cases C, D, K, and L.

\subsection{Union}
\label{sec:union}

\begin{figure}[htb]
\centering
\includegraphics[width=\columnwidth]{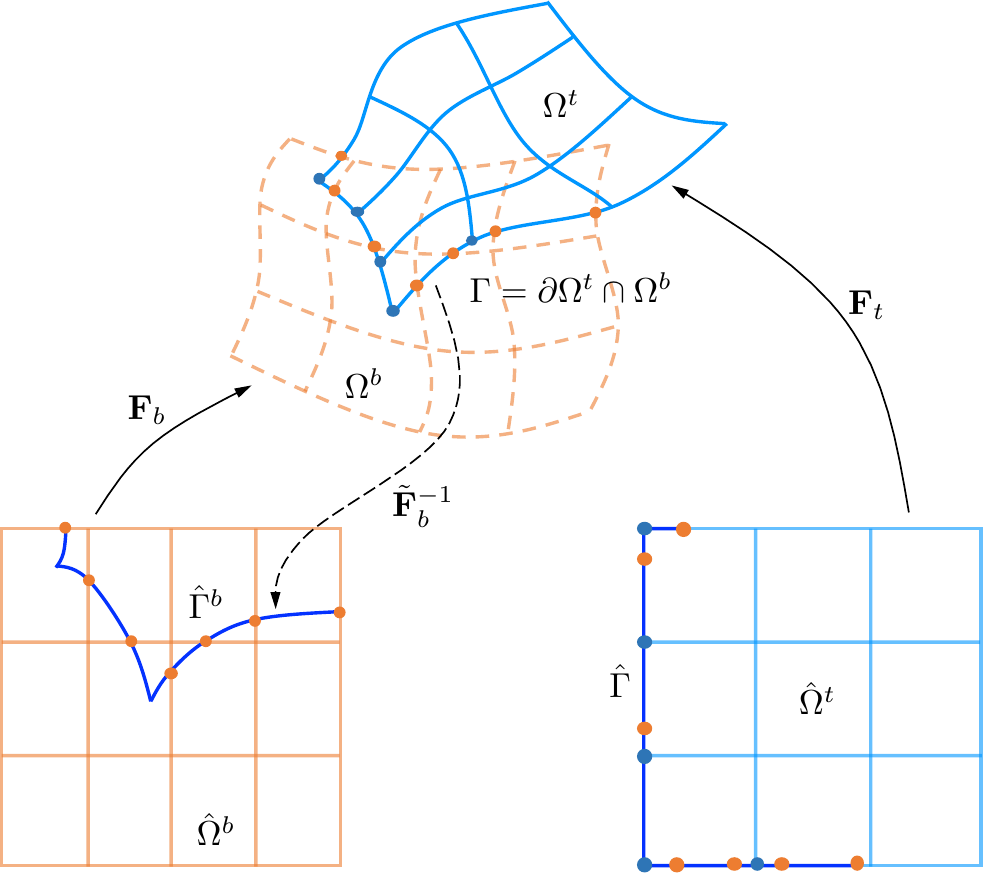}
\caption{Union of two overlapping patches and generation of their interface quadrature mesh.}
\label{fig:intf_quad}
\end{figure}

Following \cite{ref:wei19u,ref:johansson19}, we adopt a hierarchy of overlapping patches to perform the union operation. We consider a two-patch union for simplicity, which is actually sufficient to help explain our proposed method in Section \ref{sec:IBCM}. Interested readers may refer to \cite{ref:wei19u} for more general constructions. Given a pair of domains, we put one on top of another. The top and bottom domains have the geometric mappings $\bm{F}_t:\ \hat{\Omega}^t \to \Omega^t$ and $\bm{F}_b:\ \hat{\Omega}^b \to \Omega^b$, respectively. Their union is created by first trimming $\Omega^b$ with $\Omega^t$, leading to a trimmed (bottom) domain $\Omega_a^b = \Omega^b\backslash\Omega^t$. Next, $\Omega_a^b$ is \rt{combined} with $\Omega^t$ to constitute the computational domain $\Omega=\Omega_a^b \cup \Omega^t$; see Fig. \ref{fig:intf_quad}. In such an overlapping construction, the bottom patch is always trimmed whereas the top patch is intact. Note that we can choose either of the two patches to be on the top. It has been shown that different arrangements do not influence solution accuracy or matrix conditioning in linear elliptic problems~\cite{ref:wei19u}.

With trimming handled according to Section \ref{sec:trim}, analysis-aware treatment of the union operation centers on weakly coupling patches through their interfaces, where we choose Nitsche's method for its consistent and symmetric formulation \cite{ref:nitsche71}. We then need to address the following challenging issues: (1) generation of interface quadrature meshes to accurately compute involved interface integrals, and (2) stabilization of flux terms in the interface integral to guarantee the well-posedness of the problem regardless of how elements in the bottom patch are cut. 


The key to generating an interface quadrature mesh is to find a mesh intersection on the interface. This ensures that on each quadrature cell, all the involved basis functions from both top and bottom patches are polynomials rather than piecewise polynomials. In 2D, this is to find curve-curve intersections between the interface and the knot-line curves of both patches. Note that in the overlapping construction of union, the interface $\Gamma$ (in the physical domain) is part of the trimming loop in the bottom patch $\Omega^b$. Moreover, $\Gamma$ itself is part of the top patch boundary, i.e., $\Gamma = \partial\Omega^t \cap \Omega^b \subset \partial\Omega^t$. In other words, $\Gamma$ already has the mesh information of $\Omega^t$; see, for example, the blue dots in Fig.~\ref{fig:intf_quad}. It is left to find its intersections with the knot-line curves of $\Omega^b$. We proceed in the parametric domain of the bottom patch $\hat{\Omega}^b$.  An approximate preimage of $\Gamma$ \rt{is} first found in $\hat{\Omega}^b$ through the inversion algorithm and spline fitting, i.e., $\hat{\Gamma}^b=\tilde{\bm{F}}_b^{-1} \circ \Gamma$, where $\tilde{\bm{F}}_b^{-1}$ indicates an approximate inversion map.  Next, we can easily find the intersections of $\hat{\Gamma}^b$ with the axis-aligned knot lines in $\hat{\Omega}^b$; see the orange dots in $\hat{\Omega}^b$ in Fig.~\ref{fig:intf_quad}. For each resulting intersection point $\bm{\xi}$, we find its image $\bm{F}_b(\bm{\xi})$ in the physical domain and further bring it onto $\partial\hat{\Omega}^t$, leading to a desired intersection $\tilde{\bm{F}}_t^{-1}\circ \bm{F}_b(\bm{\xi})$. Once all the intersections are found, a 1D mesh can be readily constructed as the interface quadrature mesh, for example, based on the blue and orange dots in $\hat{\Omega}^t$ in Fig.~\ref{fig:intf_quad}.

\subsection{Stabilized formulations}
\label{sec:stabf}

We choose Nitsche's method \cite{ref:nitsche71} to couple independent domains for its consistent and symmetric formulation. However, the flux terms (those involving normal derivatives) in Nitsche's formulation may give rise to an instability issue when elements adjacent to an interface are badly cut, that is, only an extremely small portion of an element is left after trimming \cite{ref:prenter17,ref:puppi19}. The lack of stability indicates the potential violation of the coercivity condition in the bilinear form, leading to that the well-posedness of the problem is not guaranteed. Therefore, a proper stabilization method is needed and it often poses as one of the most challenging problems in employing Boolean operations in IGA.

In this regard, we adopt the \emph{minimal stabilization} method that was developed in our previous work on the overlapping construction of union \cite{ref:wei19u}. It was originally introduced to address the stability issue in Nitsche's formulation on trimmed boundaries \cite{ref:puppi19}. It has shown optimal convergence behaviors as well as trimming-independent conditioning in various tests. As stabilization methods are problem-dependent, let us take two linear elliptic problems, Poisson's problem and linear elasticity, as the model problems to explain the formulation details.

\textbf{Poisson's problem}. In the two-domain setting, the strong form of the Poisson's problem is stated as follows. Given $f:\,\Omega\to\mathbb{R}$,  $g_D:\,\Gamma_D\to\mathbb{R}$, and $g_N:\,\Gamma_N\to\mathbb{R}$, find $u:\,\Omega\to\mathbb{R}$ such that
\begin{equation}
\left\{
\begin{aligned}
-\Delta u = f \quad & \text{in} \quad \Omega=\Omega^t \cup \Omega_a^b, \\
u^t - u^b = 0 \quad & \text{on} \quad \Gamma = \partial \Omega^t \cap \Omega^b, \\
\nabla u^t \cdot \bm{n}^t + \nabla u^b \cdot \bm{n}^b = 0 \quad & \text{on} \quad \Gamma= \partial \Omega^t \cap \Omega^b, \\
u = g_D \quad & \text{on} \quad \Gamma_D, \\
\nabla u \cdot \bm{n} = g_N \quad & \text{on} \quad \Gamma_N, \\
\end{aligned}
\right.
\label{eq:poisson_strong}
\end{equation}
where $u^t=u|_{\Omega^t}$, $u^b=u|_{\Omega_a^b}$ are restrictions of $u$ to respective domains; $\bm{n}^t$, $\bm{n}^b$, and $\bm{n}$ are outwards unit normals of $\partial\Omega^t$, $\partial\Omega_a^b$ and $\partial\Omega$ respectively; and $\Gamma_D$ and $\Gamma_N$ ($\Gamma_D\cap\Gamma_N=\emptyset$, $\overline{\Gamma_D\cup\Gamma_N}=\partial\Omega$) are Dirichlet and Neumann boundaries, respectively. The second and third equations are transmission conditions across the interface $\Gamma$. We have also assumed that the Dirichlet boundary $\Gamma_D$ is not trimmed to simplify formulations. One may refer to \cite{ref:puppi19} for the case of trimmed $\Gamma_D$.

Before presenting the weak formulation of Problem \eqref{eq:poisson_strong}, we first introduce the following generic approximation space,
\begin{equation}
\begin{aligned}
\mathcal{V}_h^{\alpha} = &\{v_h^t \in \mathcal{B}^t:\, v_h^t|_{\Gamma_D \cap \partial\Omega^t} = \alpha \} \\
\oplus &\{v_h^b \in \mathcal{B}^b|_{\Omega_a^b}:\, v_h^b|_{\Gamma_D\cap \partial\Omega^b} = \alpha \} ,
\end{aligned}
\label{eq:approx_space}
\end{equation}
where \rt{$\alpha:\,\Gamma_D\to\mathbb{R}$} is a generic scalar function with suitable regularity, $\mathcal{B}^t$ and $\mathcal{B}^b$ are B-spline (or NURBS) spaces on $\Omega^t$ and $\Omega^b$, respectively, and $\mathcal{B}^b|_{\Omega_a^b}$ is the restriction of $\mathcal{B}^b$ to the active subdomain $\Omega_a^b$. We have the approximation spaces of trial functions $\mathcal{V}_h^{g_D}$ and test functions $\mathcal{V}_h^{0}$ when $\alpha=g_D$ and $\alpha=0$, respectively.

The discrete weak form of Problem \eqref{eq:poisson_strong} is stated as follows: Find $u_h \in \mathcal{V}_h^{g_D}$ such that 
\begin{equation}
a_h(u_h,v_h) = l (v_h), \quad \forall v_h \in \mathcal{V}_h^0,
\end{equation}
where
\begin{equation}
\begin{aligned}
a_h(u_h,v_h) &= \int_{\Omega} \nabla u_h \cdot \nabla v_h \\
&- \int_{\Gamma} \langle\nabla u_h \cdot \rt{\bm{n}^t} \rangle [v_h] - \int_{\Gamma} \langle\nabla v_h \cdot \rt{\bm{n}^t} \rangle [u_h] \\
&+ \beta (h_t^{-1} + h_b^{-1}) \int_{\Gamma} [u_h] [v_h], \\
\end{aligned}
\label{eq:bilinear}
\end{equation}
and
\begin{equation}
l(v_h) = \int_{\Omega} f \, v_h + \int_{\Gamma_N} g_N \, v_h.
\end{equation}
The transmission conditions are weakly enforced by Nitsche's method. In Eq. \eqref{eq:bilinear}, $[u_h]:=u_h^t|_{\Gamma}-u_h^b|_{\Gamma}$ is the jump term across the interface $\Gamma$, whereas $\langle\nabla u_h \cdot \rt{\bm{n}^t} \rangle$ represents the flux through $\Gamma$ \rt{and will be discussed in detail soon}. \rt{Assuming quasi-uniform meshes in $\Omega^t$ and $\Omega^b$, $h_t$ and $h_b$ represent the maximum element sizes of corresponding meshes.} $\beta$ is a penalty parameter. We take $\beta=6p_{\max}^2$ following \cite{ref:johansson19}, where $p_{\max}$ is the maximum degree in the neighboring patches of $\Gamma$.

We discuss two types of fluxes here, the symmetric average flux, 
\begin{equation}
\langle\nabla u_h \cdot \rt{\bm{n}^t} \rangle = \frac{1}{2} ( \nabla u_h^t + \nabla u_h^b ) \cdot \bm{n}^t,
\end{equation}
and the one-sided flux from $\Omega^t$,
\begin{equation}
\langle\nabla u_h \cdot \rt{\bm{n}^t} \rangle = \nabla u_h^t \cdot \bm{n}^t.
\label{eq:topflux}
\end{equation}
It has been shown that it is the flux from a certain \emph{trimmed} patch that causes the instability issue~\cite{ref:wei19u,ref:puppi19}. Therefore, if the symmetric average flux is adopted in Eq.~\eqref{eq:bilinear}, $\nabla u_h^b\cdot \bm{n}^t$ is a flux from the trimmed patch $\Omega_a^b$, and thus stabilization is needed for $\nabla u_h^b\cdot \bm{n}^t$. In contrast, no further treatment is needed to ensure stability if the one-sided flux from the non-trimmed $\Omega^t$ is used\footnote{This is true only in the two-patch union where $\Omega^t$ is already at the top of the overlapping hierarchy. In the multi-patch union, $\Omega^t$ may be cut by certain patches that are on top of $\Omega^t$, \rt{and thus $\nabla u_h^t\cdot \bm{n}^t$ may also need to be stabilized according to, for instance, the minimal stabilization~\cite{ref:wei19u}.}}. It has been shown that the two types of fluxes yield very similar results in terms of solution accuracy and matrix conditioning for linear elliptic problems~\cite{ref:wei19u}. Therefore, we choose the one-sided flux in this work to obtain a formulation without further treatment regarding stabilization.

\rt{
\begin{myremark}
The stabilization through Eq. \eqref{eq:topflux} leverages the adjacency of trimmed and untrimmed meshes. Essentially, Eq. \eqref{eq:topflux} can be viewed as imposing the third equation in Eq. \eqref{eq:poisson_strong} as a Neumann boundary condition on the trimmed domain. Since trimming at Neumann boundaries does not lead to instability, the entire need for additional stabilization vanishes.
\end{myremark}
}

\vspace{+2mm}
\textbf{Linear elasticity}. The second model problem is linear elasticity under assumption of homogeneous and isotropic material, small strains, and small displacements. The strong form is stated as follows. Given $f_i:\,\Omega\to\mathbb{R}$, $g_{D_i}:\,\Gamma_{D_i}\to\mathbb{R}$, and $g_{N_i}:\,\Gamma_{N_i}\to\mathbb{R}$, find the displacement $u_i:\,\Omega\to\mathbb{R}$ such that
\begin{equation}
\left\{
\begin{aligned}
\sigma_{ij,j} +f_i =0 \quad & \text{in} \quad \Omega = \Omega^t \cup \Omega_a^b,\\
u_i^t - u_i^b = 0 \quad & \text{on} \quad \Gamma= \partial \Omega^t \cap \Omega^b, \\
\sigma_{ij}^t n_j^t + \sigma_{ij}^b n_j^b = 0 \quad & \text{on} \quad \Gamma= \partial \Omega^t \cap \Omega^b, \\
u_i = g_{D_i} \quad & \text{on} \quad \Gamma_{D_i}, \\
\sigma_{ij} n_j = g_{N_i} \quad & \text{on} \quad \Gamma_{N_i}, \\
\end{aligned}
\right.
\label{eq:equil}
\end{equation}
where
\begin{equation}
\left\{
\begin{aligned}
\sigma_{ij} &= \lambda \delta_{ij} \delta_{kl} \epsilon_{kl} + \mu (\delta_{ik} \delta_{jl} + \delta_{il} \delta_{jk}) \epsilon_{kl},\\
\epsilon_{kl} &= \frac{1}{2} (u_{k,l} + u_{l,k}).
\end{aligned}
\right.
\label{eq:ss}
\end{equation}
In Eqs. (\ref{eq:equil},\ref{eq:ss}), indices $i$, $j$, $k$, and $l$ take on values $1,\ldots,d$, where $d\in\{2,3\}$ is the number of the spatial dimensions. The displacement field $\bm{u}$ is vector-valued with $u_i$ being the $i$-th component; the same notation applies to $f_i$, $g_i$, $h_i$, and $n_j$ (the $j$-th component of the normal $\bm{n}$). $\sigma_{ij}$ and $\epsilon_{kl}$ are Cartesian components of the Cauchy stress tensor and the infinitesimal strain tensor, respectively. The comma in $\sigma_{ij,j}$ and $u_{i,j}$ denotes differentiation with respect to the spatial coordinate, e.g., $u_{i,j} = \partial u_i / \partial x_j$. Moreover, the summation convention is applied to repeated indices, for instances, $\sigma_{ij,j} = \sigma_{i1,1}+\sigma_{i2,2}$ and $\sigma_{ij}n_j = \sigma_{i1}n_1+\sigma_{i2}n_2$ in 2D. The superscripts ``$t$" and ``$b$" again denote restrictions to the top and bottom patches, respectively. The Dirichlet and Neumann boundary conditions are applied independently in each direction and thus $\Gamma_{D_i} \cap \Gamma_{N_i} = \emptyset$ and $\overline{\Gamma_{D_i} \cup \Gamma_{N_i}} = \partial \Omega$ for $i=1,\ldots,d$. $\delta_{ij}$ is the Kronecker delta, i.e.,
\begin{equation}
\delta_{ij} = \left\{
\begin{aligned}
1 \quad &i=j,\\
0 \quad &\text{otherwise}.
\end{aligned}
\right.
\end{equation}

In Eq. \eqref{eq:ss}, constants $\lambda$ and $\mu$ are material parameters called Lam\'{e} parameters, which are often expressed in terms of the Young's modulus $E$ and and Poisson's ratio $\nu$,
\begin{equation}
\begin{aligned}
\lambda &= \frac{\nu E}{(1+\nu)(1-2\nu)}, \\
\mu &= \frac{E}{2(1+\nu)}. \\
\end{aligned}
\label{eq:lame}
\end{equation} 
In 2D, Eq. \eqref{eq:lame} falls into the plane strain assumption. Also note that different domains may be occupied by different materials, and thus the values of $\lambda$ and $\mu$ can vary from domain to domain.

The corresponding discrete weak formulation with Nitsche's method to deal with the interface is stated as follows: Find $\bm{u}_h =(u_{1,h},u_{2,h}) \in\mathcal{V}_h^{g_{D_1}}\times \mathcal{V}_h^{g_{D_2}}$ (see Eq. \eqref{eq:approx_space}) such that
\begin{equation}
a_h(\bm{u}_h,\bm{v}_h) = l(\bm{v}_h), \quad \forall \bm{v}_h=(v_{1,h},v_{2,h}) \in \mathcal{V}_h^0 \times \mathcal{V}_h^0,
\end{equation}
where
\begin{equation}
\begin{aligned}
a_h(\bm{u}_h,\bm{v}_h) &= \int_{\Omega} \sigma_{ij}(\bm{u}_h) \, \epsilon_{ij}(\bm{v}_h) \\
& -\int_{\Gamma} \langle \sigma_{ij}(\bm{u}_h) \, \rt{n_j^t} \rangle  [v_{i,h}]  
-\int_{\Gamma} \langle \sigma_{ij}(\bm{v}_h) \, \rt{n_j^t} \rangle  [u_{i,h}] \\
& +\beta \left( h_t^{-1}+ h_b^{-1} \right) \int_{\Gamma}  [u_{i,h}]  [v_{i,h}] ,\\ 
\end{aligned}
\label{eq:bilinear_elas}
\end{equation}
and
\begin{equation}
l(\bm{v}_h) = \int_{\Omega} f_i \, v_{i,h} + \sum_{j=1}^d \left( \int_{\Gamma_{N_j}} g_{N_i} \, v_{i,h} \right).
\end{equation}
In Eq. \eqref{eq:bilinear_elas}, $[u_{i,h}]:=u_{i,h}^t - u_{i,h}^b$ and $\langle \sigma_{ij}(\bm{u}_h) \rt{n_j^t}  \rangle$ denote the displacement jump and the stress flux across $\Gamma$, respectively. We again adopt the one-sided flux, i.e.,
\begin{equation}
\langle \sigma_{ij}(\bm{u}_h) \, \rt{n_j^t}  \rangle = \sigma_{ij}^t (\bm{u}_h^t) \, n_j^t,
\end{equation}
which does not need further treatment for stabilization. The penalty parameter $\beta$ depends on both spline degrees and material properties,
\begin{equation}
\beta = 6 p_{\max}^2 \times 8(3\lambda_{\max}+2\mu_{\max}),
\end{equation}
where $\lambda_{\max}$ and $\mu_{\max}$ are the maximum Lam\'{e} constants of neighboring patches. The choice is inspired by \cite{ref:johansson19,ref:hansbo05}. 

\begin{figure*}[htb]
\centering
\includegraphics[width=\textwidth]{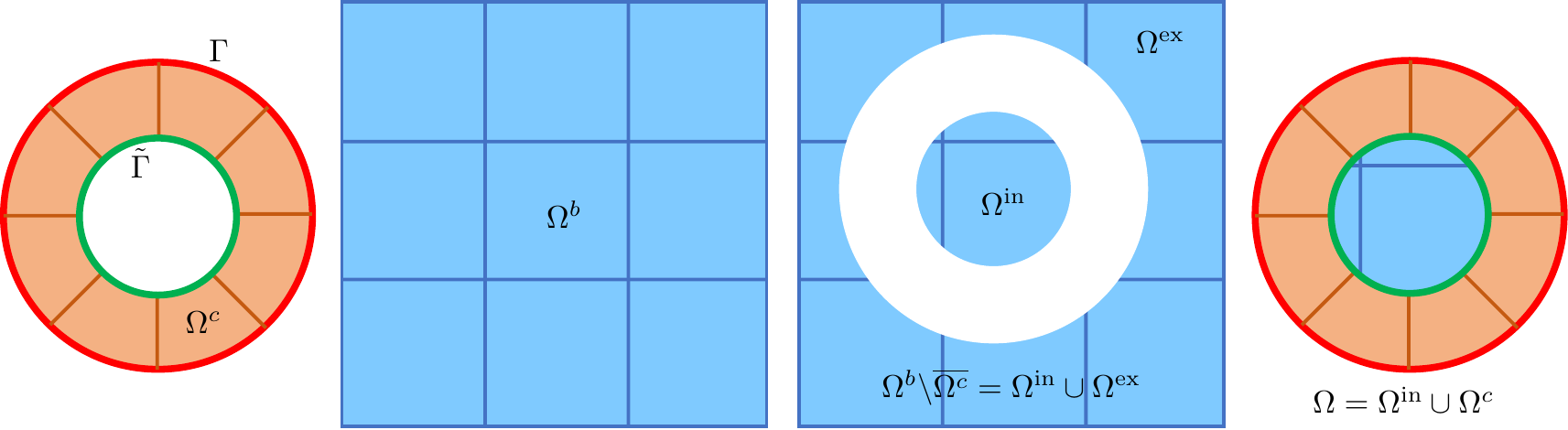}
\begin{tabular}{cccc}
(a) \hspace{+0.2\textwidth}  & (b) \hspace{+0.23\textwidth} & (c) \hspace{+0.22\textwidth} & (d) \\
\end{tabular}
\caption{Geometric construction of IBCM for the boundary type. (a) The input boundary curve $\Gamma$ (the red circle) and the resulting conformal layer $\Omega^c$ by extruding $\Gamma$ inwards to the target curve $\tilde{\Gamma}$, (b) the background patch $\Omega^b$, (c) two disconnected regions ($\Omega^{\mathrm{ex}}$ and $\Omega^{\mathrm{in}}$) by cutting $\Omega^b$ with $\Omega^c$, and (d) the computational domain $\Omega$ composed of $\Omega^{\mathrm{in}}$ and $\Omega^c$.}
\label{fig:boundary_type}
\end{figure*}

\begin{figure*}[htb]
\centering
\begin{tabular}{cccc}
\includegraphics[width=0.23\textwidth]{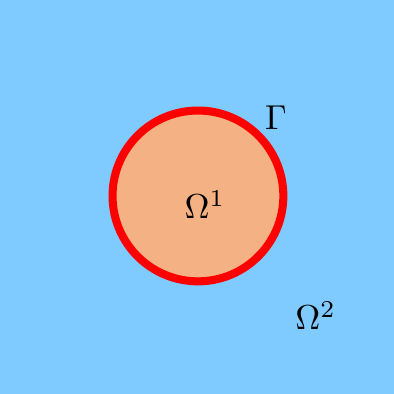}&\hspace{-2mm}
\includegraphics[width=0.23\textwidth]{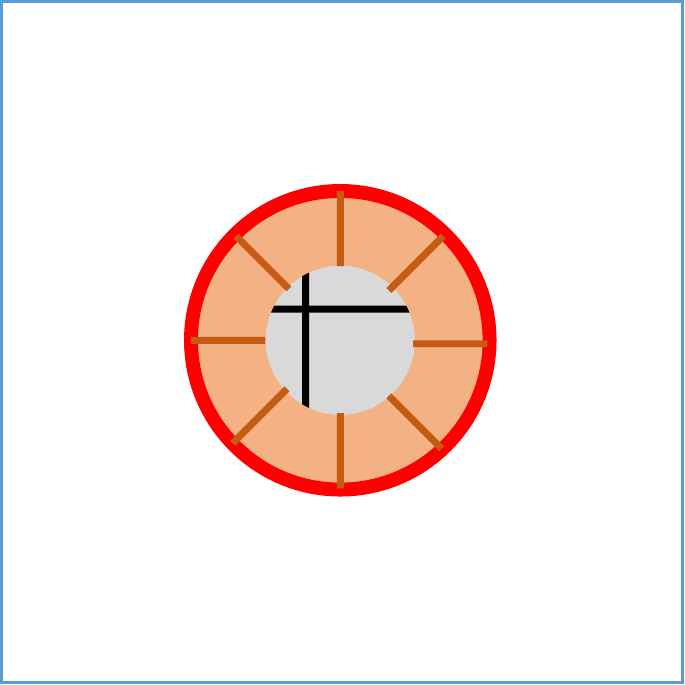}&\hspace{-2mm}
\includegraphics[width=0.23\textwidth]{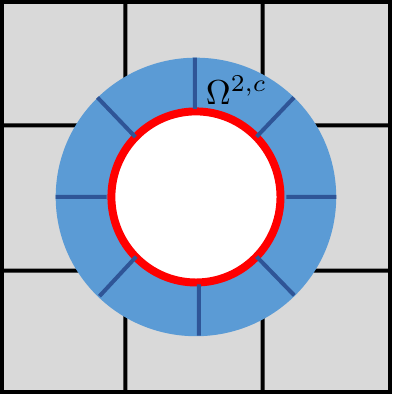}&\hspace{-2mm}
\includegraphics[width=0.23\textwidth]{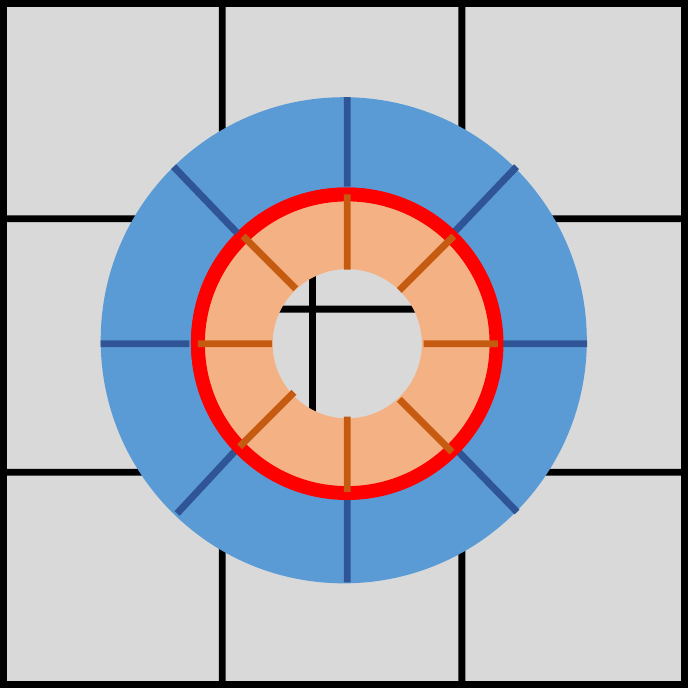}\\
(a)  & (b) & (c) & (d) \\
\end{tabular}
\caption{Geometric construction of IBCM for the interface type. (a) Two domains with different material properties and the interface $\Gamma$, (b, c) the IBCM representation of $\Omega^1$ and $\Omega^2$, respectively, and (d) the entire IBCM representation of $\Omega^1\cup \Omega^2$.}
\label{fig:interface_type}
\end{figure*}

\section{Immersed Boundary-Conformal Isogeometric Method}
\label{sec:IBCM}

In this section, we introduce the proposed method, namely the \emph{Immersed Boundary-Conformal Method} (IBCM). We start with two types of geometric constructions using IBCM. We then proceed to discuss the key technologies that guarantee IBCM to work properly, including parameterization-consistent extrusion and analysis-aware treatment of Boolean operations in IBCM.

\subsection{Geometric construction}
\label{sec:IBCM_geom}

The geometric construction of IBCM is conceptually simple and can be divided into two types depending on the features of interest: the boundary type and the interface type. The \emph{boundary type} aims to capture boundary features, such as boundary shapes and local solution features near boundaries. The \emph{interface type}, on the other hand, deals with interfaces between different materials such as inclusions and fiber-reinforced materials, where stresses exhibit discontinuity across material interfaces. The major steps in constructing an IBCM representation are the same in both types, including extrusion, trimming and union. 

\textbf{Boundary-type construction}. We first explain the boundary type. In the extrusion step, we start with the representation of a boundary $\Gamma$, which is a closed loop formed by a set of oriented and connected B-spline/NURBS curves. We extrude the loop inwards to yield a ring-like layer $\Omega^c$ that is obviously conformal to the given boundary, so we call it a \emph{conformal layer}. $\Omega^c$ is generally represented by multiple B-spline/NURBS patches. Details of constructing $\Omega^c$ will be discussed in Section \ref{sec:extrusion}.

We explain the remaining steps with the reference to Fig.~\ref{fig:boundary_type}. In the trimming step, we embed $\Omega^c$ into a sufficiently large background domain $\Omega^b$ such that $\Omega^c \subset \Omega^b$. $\Omega^b$ is usually represented by a B-spline mesh defined on a Cartesian grid. As a result, $\Omega^b$ is cut by $\Omega^c$ into two disconnected regions: an exterior region $\Omega^{\mathrm{ex}}$ and an interior region $\Omega^{\mathrm{in}}$. In other words, we have $\Omega^b \backslash \overline{\Omega^c}=\Omega^{\mathrm{ex}} \cup \Omega^{\mathrm{in}}$. 

Finally in the union step, $\Omega^{\mathrm{in}}$ is coupled with $\Omega^c$ to constitute the computational domain $\Omega$, i.e., $\Omega=\Omega^{\mathrm{in}} \cup \Omega^c$. This way, we obtain the IBCM representation of $\Omega$. 

\textbf{Interface-type construction}. We next discuss the interface-type construction, which can be obtained by repeatedly applying the boundary-type construction. Let $\Gamma$ denote an interface of two domains $\Omega^1$ and $\Omega^2$, each of which has its own material properties. We explain details in the following with the help of Fig. \ref{fig:interface_type}.

We first create an IBCM representation for $\Omega^1$, which in fact follows the same procedure as the boundary-type construction; see Fig. \ref{fig:interface_type}(b). The IBCM construction for $\Omega^2$ is almost the same. Now, $\Gamma$ is extruded towards the interior of $\Omega^2$ to obtain a conformal layer $\Omega^{2,c}$. Clearly, the extrusion direction is opposite to that for $\Omega^1$. After the corresponding background mesh is cut, the exterior region is coupled with $\Omega^{2,c}$; see Fig. \ref{fig:interface_type}(c). Finally, the two IBCM representations are combined together on $\Gamma$ in a conformal manner, leading to an entire IBCM representation of $\Omega^1\cup \Omega^2$, as shown in Fig. \ref{fig:interface_type}(d).

In general, when more than one boundary/interface feature appears in geometric modeling, we need to construct an IBCM representation for each of them. Combining all the resulting IBCM representations together is straightforward through conformal interfaces. We will present such an example in Section~\ref{sec:result}.

\begin{myremark}
IBCM aims to leverage the geometric flexibility of immersed methods with the advantages of boundary-fitted methods. While the large portion of a geometry is represented following the immersed manner, its boundary (or interface), as the key geometric feature, has a boundary-fitted representation through the conformal layer. In other words, the conformal mesh is placed where it is needed most. On the other hand, meshing is still needed in IBCM to obtain a desired conformal layer through extrusion. Although it remains a challenge in general cases, extrusion is much easier to manage than finding a boundary-conformal parameterization for the entire domain. From this perspective, IBCM helps alleviate the meshing difficulties encountered in IGA while retaining geometry-aligned discretization around key geometric features.
\end{myremark}

\begin{myremark}
Several benefits of conformal discretization are immediately available in IBCM thanks to the conformal layer. First, it is possible to strongly impose Dirichlet boundary conditions, which is preferable when point-wise satisfaction is desired (e.g., the clamped boundary in solid mechanics). In contrast, imposing Dirichlet boundary conditions often poses as one of the biggest challenges in immersed methods. Second, IBCM is also suitable to model material interfaces when perfect bonding is the case, as the kinematic constraints are naturally met due to the conformal representation of the interface.
\end{myremark}

\begin{myremark}
In immersed boundary methods, the geometry boundary always serves as the trimming loop on the background mesh. In contrast, IBCM moves trimming curves away from the boundary/interface. As such geometric features are often critical to solution accuracy and trimming is the origin of various issues in analysis, we expect that by separating the two, IBCM can benefit solution accuracy. Moreover, as conformal layers align with geometric features, it is much more intuitive and convenient for IBCM to control mesh resolutions there than by using immersed methods.
\end{myremark}

\subsection{Parameterization-consistent extrusion}
\label{sec:extrusion}

We aim for a conformal layer $\conformalLayer$ that (1) represents the boundary or interface $\interface$ exactly and (2) is consistent with the parameterization of $\interface$.
First, we define a target curve $\targetCurve$, which specifies the interface of $\conformalLayer$ to the interior region of the background domain $\backDomainIn$; see Fig.~\ref{fig:boundary_type}(a).

From a geometric point of view, the most natural choice for $\targetCurve$ is an offset to $\interface$.
Offset curves are defined as the locus of points that are at a constant distance $d$ along the normal vector from the so-called progenitor curve \cite{Hoschek1992b,patrikalakis2009b},
and there are several techniques in the literature to generate them; see e.g.~\cite{Maekawa1999a,Pham1992a}.
\rt{These schemes are well-suited to preserve geometric features such as kinks and cusps.
In general, 
the resulting offset curves} are more complex than the progenitor curve, and they may have \rt{additional} cusps and may self-intersect locally or globally. 
\begin{figure}[t]
    \centering
    \includegraphics[width=0.8\columnwidth]{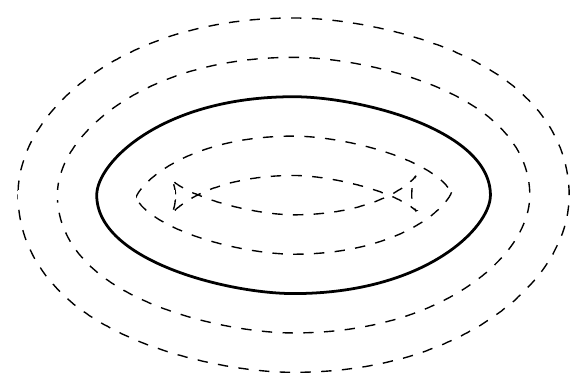}
    \caption{A progenitor curve (thick) and various offset curves (dashed) with different distance $d$.}
    \label{fig:offsetExample}
\end{figure}
These self-intersections are linked to the distance $d$ of the offset curve (see Fig.~\ref{fig:offsetExample}). 
In particular, local self-intersections occur in concave regions when the absolute value of $d$ exceeds the minimum radius of curvature, and global self-intersections occur when the distance between two distinct points on the progenitor curve is smaller than $2d$ \cite{patrikalakis2009b}.
Wallner et al.~\cite{WALLNER2001a} presented an approach to determine the maximal offset distance to avoid self-intersections.

To sum up, using offset curves for the construction of the conformal layer $\conformalLayer$ may imply that the parameterization from the target curve $\targetCurve$ does not match with the one from the interface $\interface$, especially when $d$ is large. 
In this case, $\conformalLayer$ may be constructed as a loft surface between $\interface$ and $\targetCurve$, where a superset of the knot vectors of $\interface$ and $\targetCurve$ define the resulting parameterization. 
\rt{A loft surface is created by fitting a series of given curves, with the control over the tangents of the surface. The related function is available in many CAD systems such as Rhinoceros \cite{ref:rhino}.}

This approach has some disadvantages: First, the final parameterization of $\conformalLayer$ is partly determined by the offset curve scheme, and second, $C^0$-continuities may be introduced when the offset curve has  \rt{additional} cusps or self-intersections that have been trimmed away.

Fortunately, we are very flexible in the definition of the target curve $\targetCurve$ in terms of the distance $d$ and its shape, as will be shown in Section \ref{sec:result}. Hence, we suggest the following procedure for the construction of $\conformalLayer$:
\begin{enumerate}
    \item Define $\targetCurve$ either by an offset curve or any other curve that does not intersect $\interface$ and is at least as smooth as $\interface$.
    \item Project the Greville points of $\interface$ onto $\targetCurve$.
    \item Use the projected Greville points to construct an approximation of $\targetCurve$ that has the same \rt{knots} and degree as $\interface$.
\end{enumerate}
%
\rt{Greville points are defined as follows. We consider a degree-$p$ NURBS curve whose knot vector is $\{\xi_1,\xi_2,\ldots,\xi_{n+p+1}\}$, where $n$ is the number of control points. The Greville abscissae are defined as an average of certain knots: $g_i=(\xi_{i+1}+\cdots +\xi_{i+p})/p$, $i=1,\ldots,n$. Letting $\bm{C}(\xi)$ be the geometric mapping of the NURBS curve, we call $\bm{C}(g_i)$ its Greville points.}

This Greville point projection allows a straightforward realization of the conformal layer $\conformalLayer$ as a parameterization-consistent extrusion of $\interface$. 
\rt{
Note that the knot vector of the boundary or interface $\interface$ determines the continuity of the final approximation of $\targetCurve$, which may be less smooth than the target curve $\targetCurve$. In general, this does not lead to any complications in the construction. However, for sharp features of $\interface$ such as cusps, it is beneficial to incorporate the present $C^0$ continuity already in the target curve. To do so, an offset with a relatively small distance $d$ can be employed to capture the corresponding portion of $\targetCurve$.
}

\begin{myremark}
\rt{When} the target curve $\targetCurve$ is given by an offset curve, we first check if it has \rt{additional} cusps or regions with high curvature. \rt{If so, we} replace them with rounded fillets to avoid that projected Greville points coincide.
\end{myremark}

\rt{
\begin{myremark}
    The proposed construction of the conformal layer $\conformalLayer$ does not guarantee that the resulting parameterization is bijective.
    Thus, we check if the Jacobian determinant is greater than zero at the integration points.
    If this is not the case, $\conformalLayer$ may be reconstructed either by simply decreasing the offset distance $d$ or using more advanced techniques such as elliptic grid generation~\cite{ref:hinz18}, the introduction of internal guiding curves~\cite{Randrianarivony2006phd} or integer-grid maps~\cite{Bommes2013a}. 
    It is pointed out that all examples in this paper do not need any correction of the conformal layer because the parameterizations obtained are bijective; see Section \ref{sec:result}.
    \end{myremark}
}

\subsection{Boolean operations in IBCM}
\label{sec:IBCM_bool}

Recall that the geometric construction of IBCM involves trimming and union. Trimming applies to the background B-spline mesh, where we follow Section \ref{sec:trim} to reparameterize cut elements for numerical integration. Note that the geometric mapping of the background patch is usually an identity map, which eliminates the need to find an approximate trimming loop in the parametric domain through the inversion algorithm.

In the union operation, part of the trimmed background mesh is coupled with the conformal layer using Nitsche's method. From the overlapping perspective, the background patch lies on the bottom whereas the conformal layer is on top. We follow Section \ref{sec:union} to generate an interface quadrature mesh. Moreover, as discussed in Section \ref{sec:stabf}, with the one-sided flux from the (non-trimmed) conformal layer, Nitsche's formulation is stable and needs no further treatment. To this end, Boolean operations are properly handled in IBCM such that it can be readily applied to linear elliptic problems.

\begin{myremark}
As weak coupling depends on the problem of interest, extending the method to other problems, such as shells, nonlinear elasticity, and Stokes or Navier-Stokes problems, requires further investigation particularly in the analysis side. A different coupling method may also be needed when it is too cumbersome to use Nitsche's formulation. There are indeed many challenging problems for IBCM to accommodate a larger class of problems. Among them, our priority lies in extending the method to 3D, which needs robust implementation to support trimming, union and extrusion.
\end{myremark}

\section{Numerical Examples}
\label{sec:result}

In this section, we first study several 2D benchmark problems to evaluate the accuracy and convergence of the proposed method. They include problems concerning boundary features as well as modeling interfaces between different materials. The former includes examples of a plate with a hole, an L-shaped domain, and a flower geometry, whereas the latter studies a bimaterial disk. Moreover, we present two examples to show the capability of IBCM in representing complex geometric features, which includes a spanner model and a fiber-reinforced composite. \rt{In the end, we push a step forward to apply IBCM to an advection-diffusion problem, where we will observe that IBCM can help efficiently resolve the boundary-layer phenomenon in an intuitive manner.} In all the tests, biquadratic spline bases are used. In addition, Dirichlet boundary conditions are always strongly imposed \rt{except for the advection-diffusion problem}, but as splines are not interpolatory, non-homogeneous Dirichlet data needs to be projected to the involved spline spaces.

\subsection{Plate with a hole}

\begin{figure}[htb]
\centering
\includegraphics[width=0.7\columnwidth]{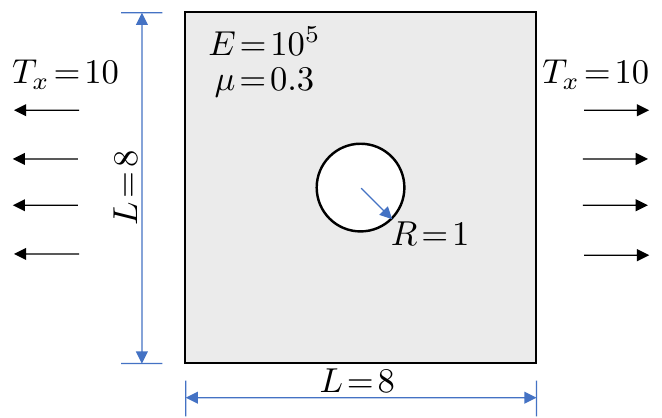}
\caption{Problem settings of the plate-with-a-hole example. $T_x$ is the applied traction at infinity in the $x$ direction.}
\label{fig:platehole_set}
\end{figure}

We start with the plate-with-a-hole test. It is a linear elasticity problem where an infinite plate with a circular hole is under constant in-plane tension. A finite portion with the hole being at the center is taken for the numerical test; see Fig. \ref{fig:platehole_set} for the problem settings.

\begin{figure}[htb]
\centering
\begin{tabular}{cc}
\includegraphics[width=0.46\columnwidth]{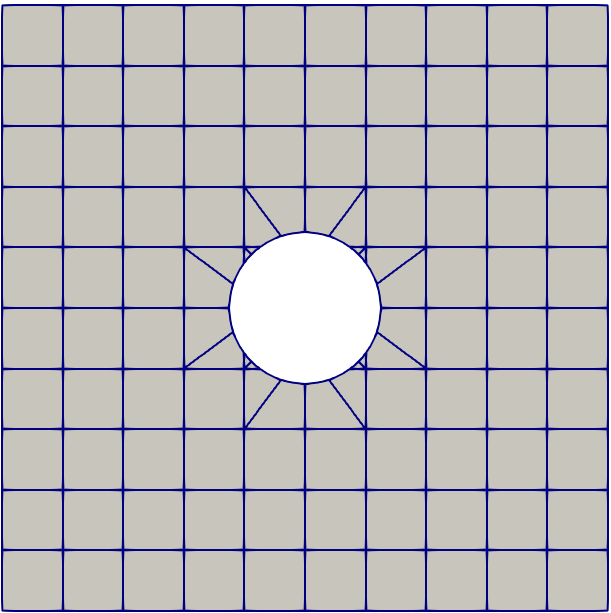}&\hspace{-2mm}
\includegraphics[width=0.46\columnwidth]{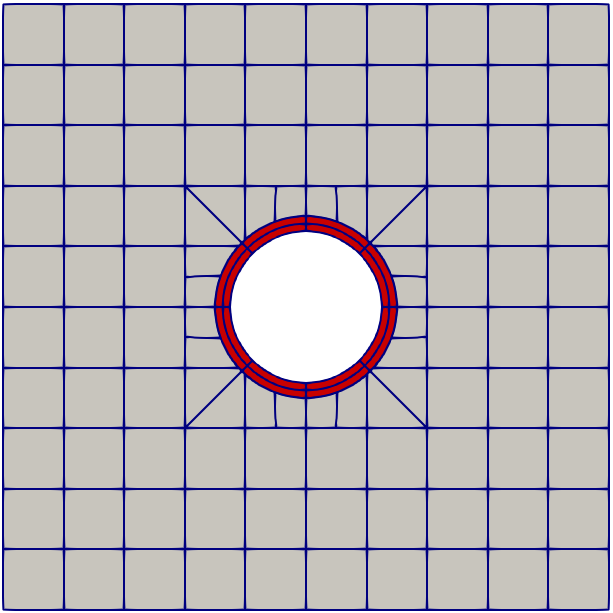}\\
(a) & (b) \\
\end{tabular}
\caption{Initial B\'ezier meshes of the plate-with-a-hole problem via trimming (a) and IBCM (b). In (b), an extra annulus is added around the hole as the conformal layer. \rt{In the background mesh, integration cells of trimmed elements are used for visualization.}}
\label{fig:platehole_geom}
\end{figure}

Clearly, the hole is a boundary geometric feature of interest. We study two kinds of geometric constructions for comparison, trimming versus IBCM; see Fig.~\ref{fig:platehole_geom}. In the trimming construction, the hole is represented by a NURBS curve, which is also the trimming loop of the background mesh. In this case, we only need to reparameterize cut elements for numerical integration. Alternatively in an IBCM construction, we observe that a conformal layer (marked in red), represented by a NURBS patch, is an annulus sitting on top of the background mesh. The geometric feature is represented by a conformal discretization in the computation domain. The conformal layer is then coupled with the cut background mesh through the union operation. Moreover, the trimming loop of the background mesh becomes the outer circle of the annulus, which, as opposed to the trimming construction, is no longer the geometric feature itself.

\begin{figure}[htb]
\centering
\begin{tabular}{cc}
\includegraphics[width=0.46\columnwidth]{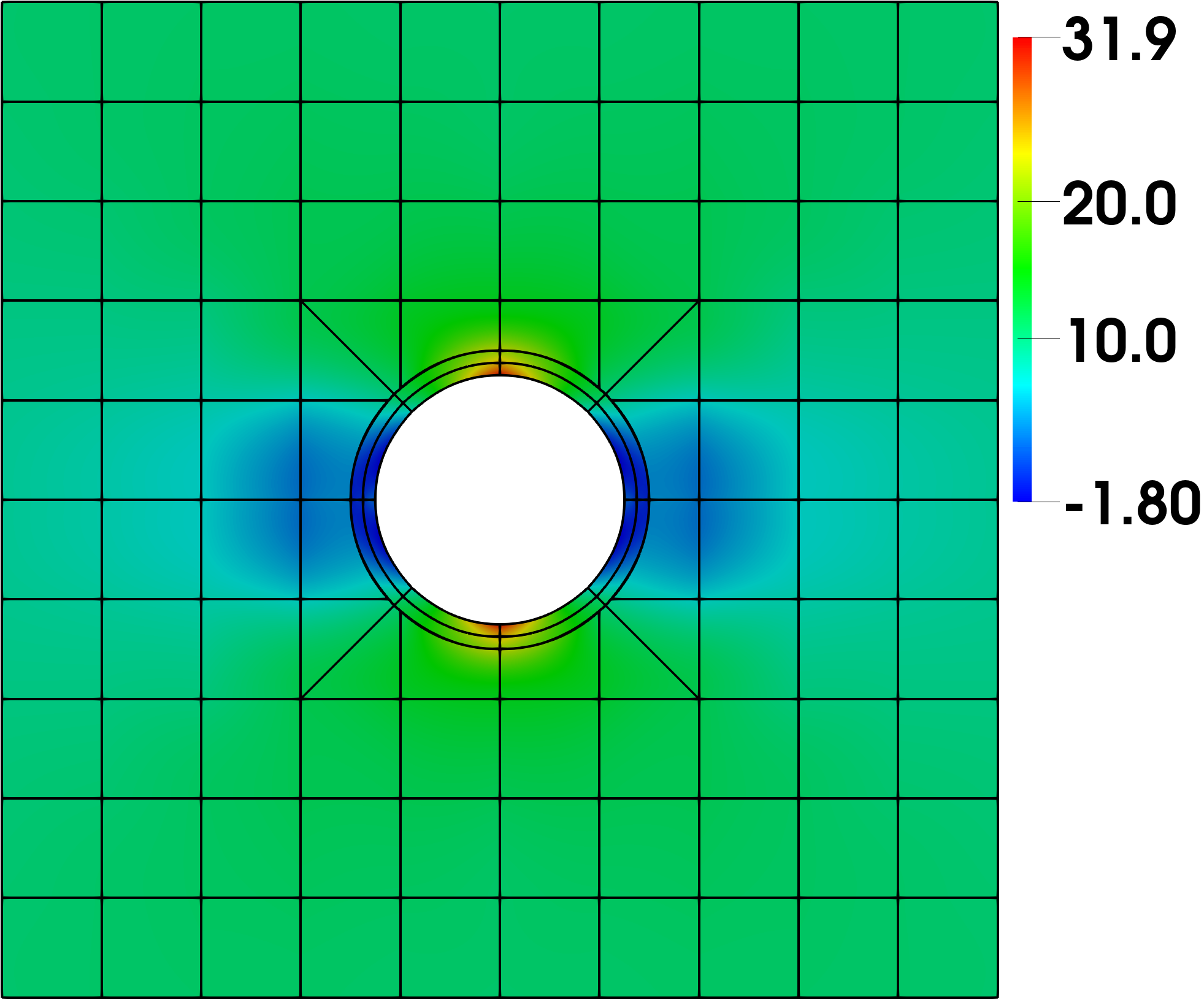}&\hspace{-2mm}
\includegraphics[width=0.46\columnwidth]{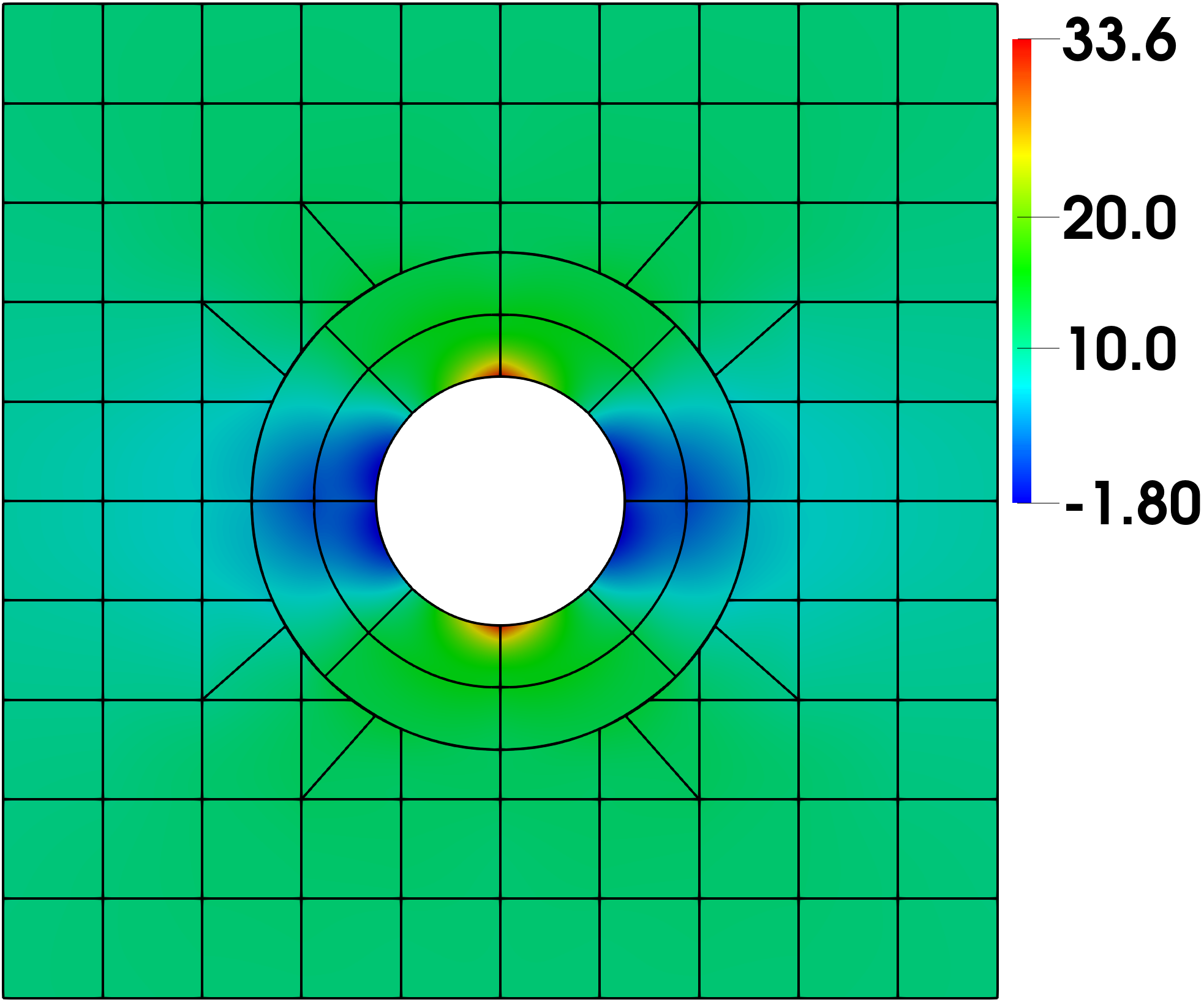}\\
(a) $t=0.2R$ & (b) $t=1.0R$ \\
\end{tabular}
\caption{Different thicknesses of the conformal layer. The color map represents numerical solution of $\sigma_{xx}$ on the initial mesh. The exact maximum $\sigma_{xx}$ is 30. \rt{In the background mesh, integration cells of trimmed elements are used for visualization.}}
\label{fig:platehole_t}
\end{figure}

\begin{figure}[htb]
\centering
\includegraphics[width=0.48\columnwidth]{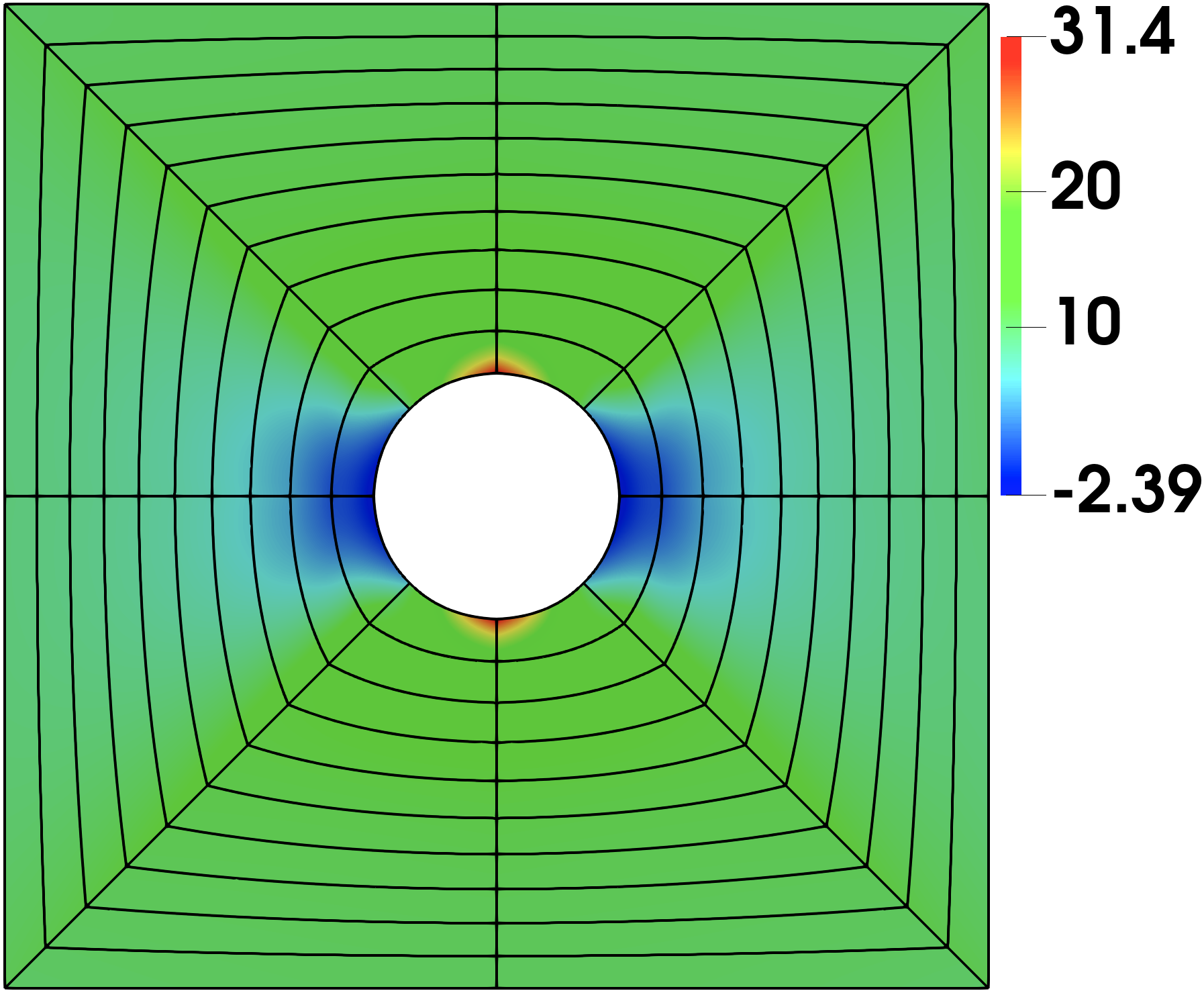}
\caption{Conformal discretization for the whole plate-with-a-hole geometry. The color map represents numerical solution of $\sigma_{xx}$ on the initial mesh. The exact maximum $\sigma_{xx}$ is 30.}
\label{fig:platehole_conform}
\end{figure}

\begin{figure}[htb]
\centering
\begin{tabular}{cc}
\includegraphics[width=0.8\columnwidth]{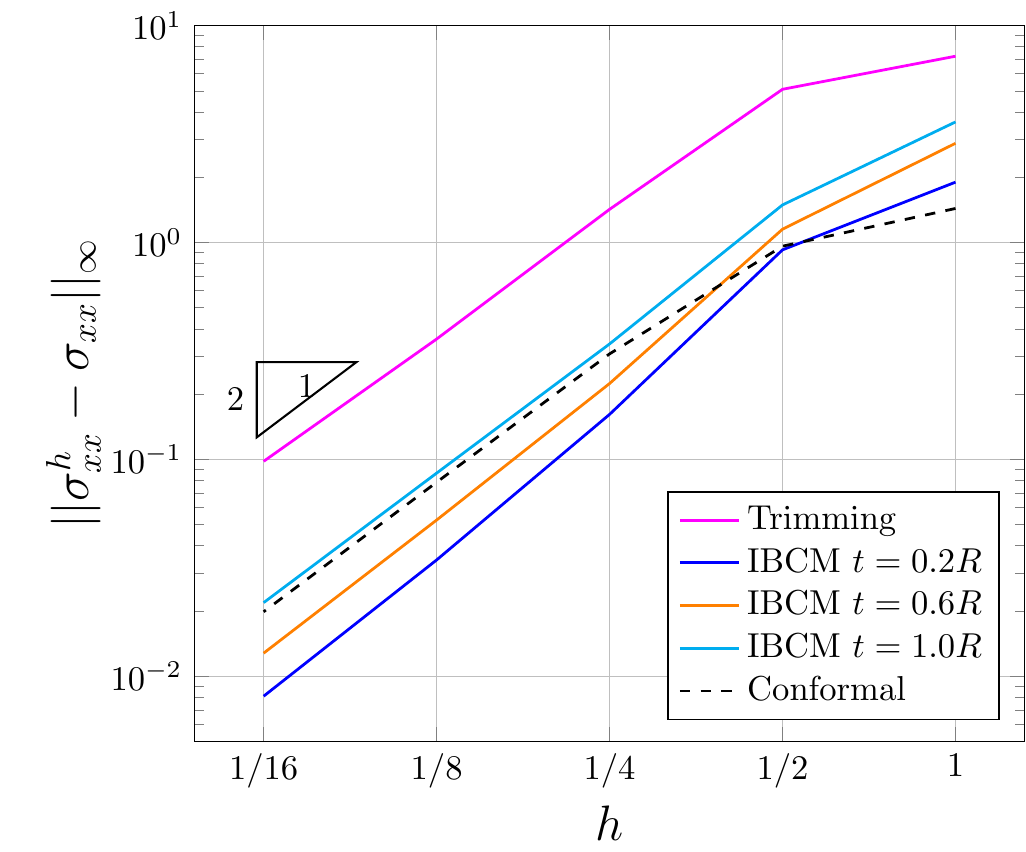}\\
(a)\\
\includegraphics[width=0.8\columnwidth]{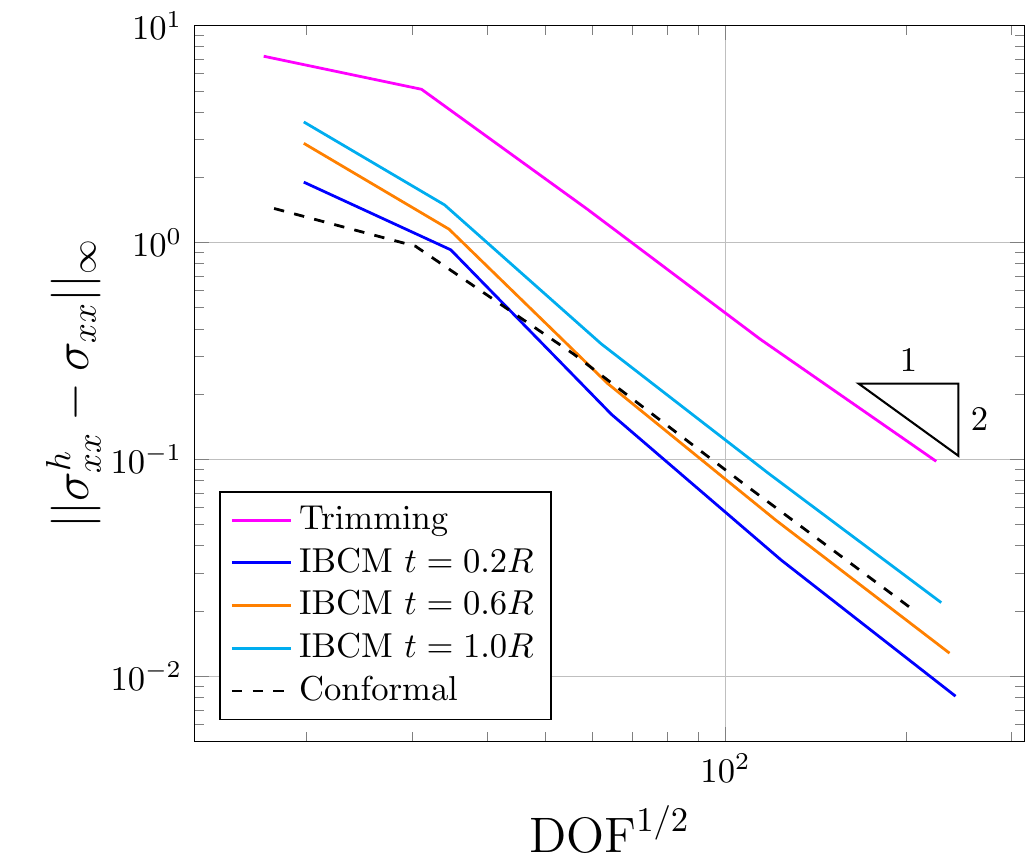}\\
(b)\\
\end{tabular}
\caption{Convergence plots with respect to $h$ (a) and $\mathrm{DOF}^{1/2}$ (b) in the plate-with-a-hole problem.}
\label{fig:platehole_conv}
\end{figure}

We also study the influence of the thickness $t$ of the conformal layer. We consider $t=0.2R$, $t=0.6R$ and $t=1.0R$, with $R$ being the radius of the hole. In each case, the conformal layer (i.e., the annulus) is represented by a $8\times 2$ mesh; see Fig. \ref{fig:platehole_t}.

\rt{As a reference, we further take a conformal discretization for the whole geometry, which in this case is easy to obtain via a single biquadratic NURBS patch; see the input B\'ezier mesh in Fig. \ref{fig:platehole_conform}. The element size around the hole is comparable to that of IBCM with $t=0.6R$.}

We summarize the convergence plots in Fig. \ref{fig:platehole_conv}, with respect to both the element size indicator $h$ and the square root of degrees of freedom $\mathrm{DOF}^{1/2}$. We observe that compared to the trimming construction, all IBCM constructions yields much more accurate results. As expected, the case of $t=0.2R$ has the smallest error because its conformal layer has the smallest element size around the hole. The other two cases ($t=0.6R$ and $t=1.0R$), on the other hand, also achieve comparable accuracy and still improve a lot compared to the trimming construction. This implies that the solution may not be sensitive to the thickness of the conformal layer, which, however, needs further study to conclude. \rt{We also notice in Fig. \ref{fig:platehole_conv}(b) that the result of the full conformal discretization is very close to that of IBCM with $t=0.6R$, which is expected because they have a comparable element size around the hole. In other words, IBCM can achieve the same level of accuracy as a full conformal discretization when critical geometric features are resolved similarly.}

\begin{figure}[htb]
\centering
\begin{tabular}{cc}
\includegraphics[width=0.46\columnwidth]{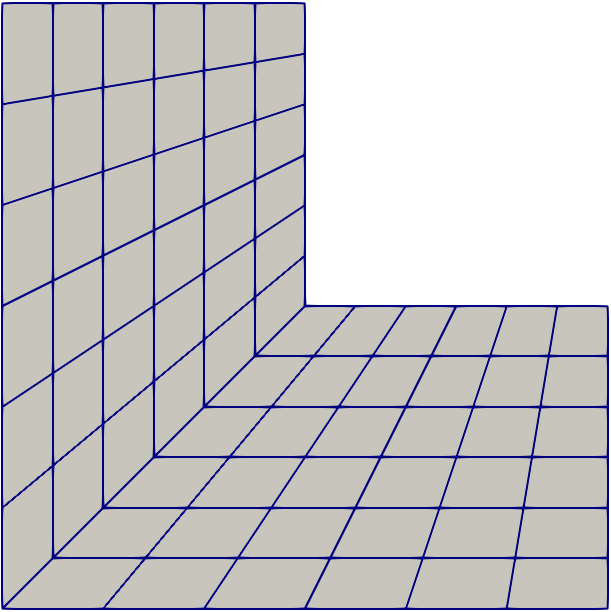}&\hspace{-2mm}
\includegraphics[width=0.46\columnwidth]{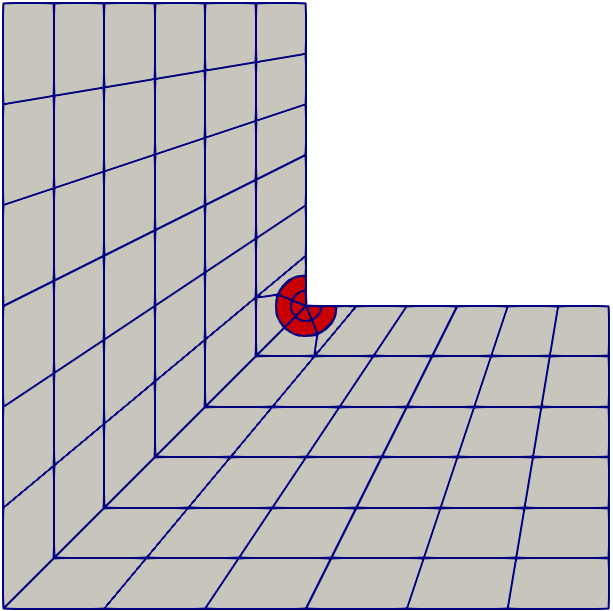}\\
(a) & (b)\\
\end{tabular}
\caption{Initial B\'ezier meshes of the L-shaped domain by a single B-spline patch (a) and via IBCM (b). In (b), a conformal layer (shaded red) is added on top of the interior material to yield a conformal discretization around the material interface.}
\label{fig:Lshape_geom}
\end{figure}

\subsection{L-shaped domain}

This example is aimed to show how IBCM \rt{improves} solution accuracy when the solution has a local feature, such as a large gradient or even singularity. Such local phenomena are often closely related to geometric features. For example, the solution gradient may exhibit singularity at certain corners. The key idea is to add a layer of mesh to where local features of the solution are expected. This indeed needs ``rough" a priori knowledge about the solution field.

\begin{figure}[htb]
\centering
\begin{tabular}{c}
\includegraphics[width=0.8\columnwidth]{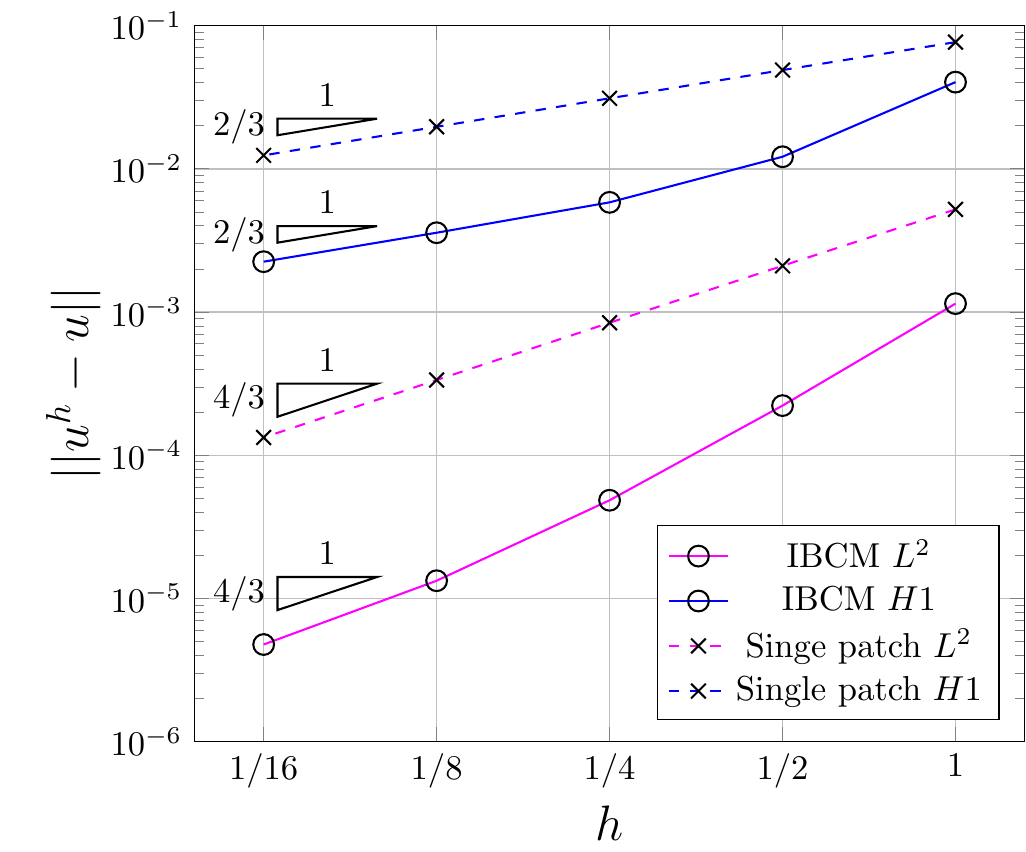}\\
(a)\\
\includegraphics[width=0.8\columnwidth]{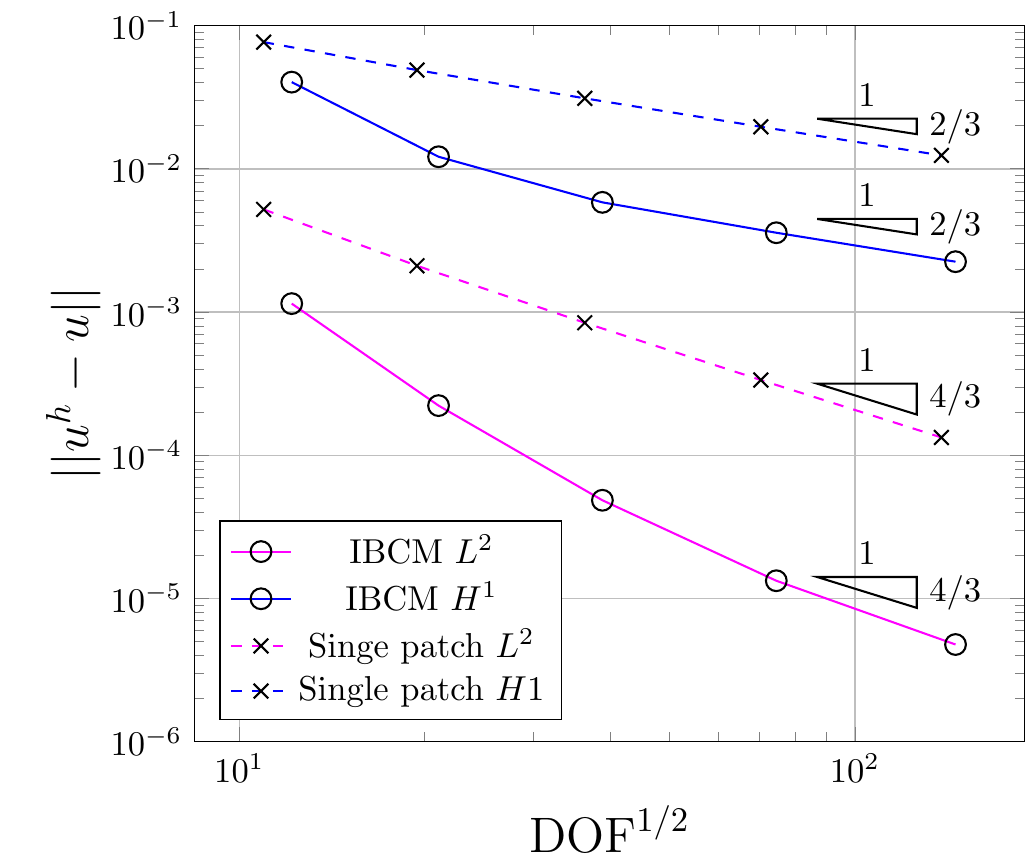}\\
(b)\\
\end{tabular}
\caption{Convergence plots with respect to $h$ (a) and $\mathrm{DOF}^{1/2}$ (b) in the L-shaped domain problem. \rt{In the background mesh of (b), integration cells of trimmed elements are used for visualization.}}
\label{fig:Lshape_conv}
\end{figure}

As an example, we solve the Laplace equation on a L-shaped domain, where the solution at the reentrant corner exhibits a singularity in its first order derivative. As a result, the convergence rates are governed by the solution regularity rather than the degree of the basis, which are expected to be $\frac{4}{3}$ and $\frac{2}{3}$ for $L^2$- and $H^1$-norm errors, respectively. As a reference test shown in Fig. \ref{fig:Lshape_geom}(a), the domain is represented by a single B-spline patch that is globally $C^0$-continuous due to sharp corners. Alternatively, we put an extra layer of mesh (three quarters of a disk) on top of the B-spline patch; see Fig. \ref{fig:Lshape_geom}(b). The extra layer is represented by a degenerated NURBS patch. Note that such a construction is merely a union of overlapping patches rather than an IBCM construction, but the idea coincides with IBCM in the spirit of better capturing local geometry/solution features by adding extra layers. Therefore, we treat it as a special case of IBCM. We compare and summarize the convergence plots of the two geometric constructions in Fig. \ref{fig:Lshape_conv}. We observe that both constructions eventually yield the same convergence rates but IBCM, as expected, greatly improves accuracy with the same DOF (or mesh size).

\rt{On the other hand, it is worth mentioning that adaptive mesh refinement is usually the method of choice to recover optimal convergence for problems with irregular solutions. In IGA, T-splines \cite{ref:sederberg03,ref:veiga12,ref:xli18,ref:casquero20}, hierarchical B-splines \cite{ref:vuong11,ref:giannelli12,ref:wei15,ref:buffa16,ref:bracco19}, and locally-refinable B-splines \cite{ref:johannessen14,ref:patrizi20} are typical examples in this family of methods.}

\subsection{Flower}

\begin{figure}[htb]
\centering
\begin{tabular}{cc}
\includegraphics[width=0.46\columnwidth]{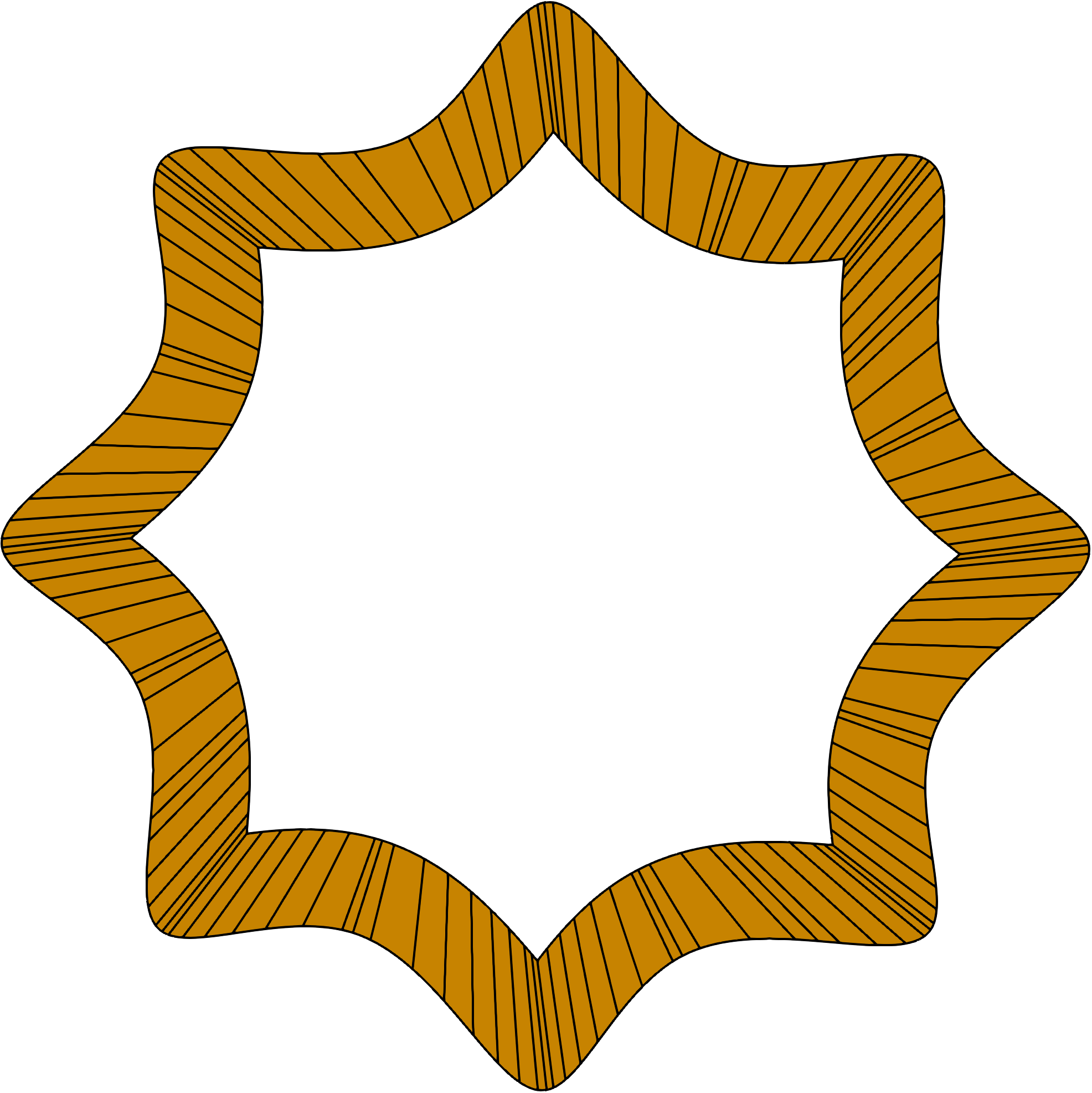}&\hspace{-2mm}
\includegraphics[width=0.46\columnwidth]{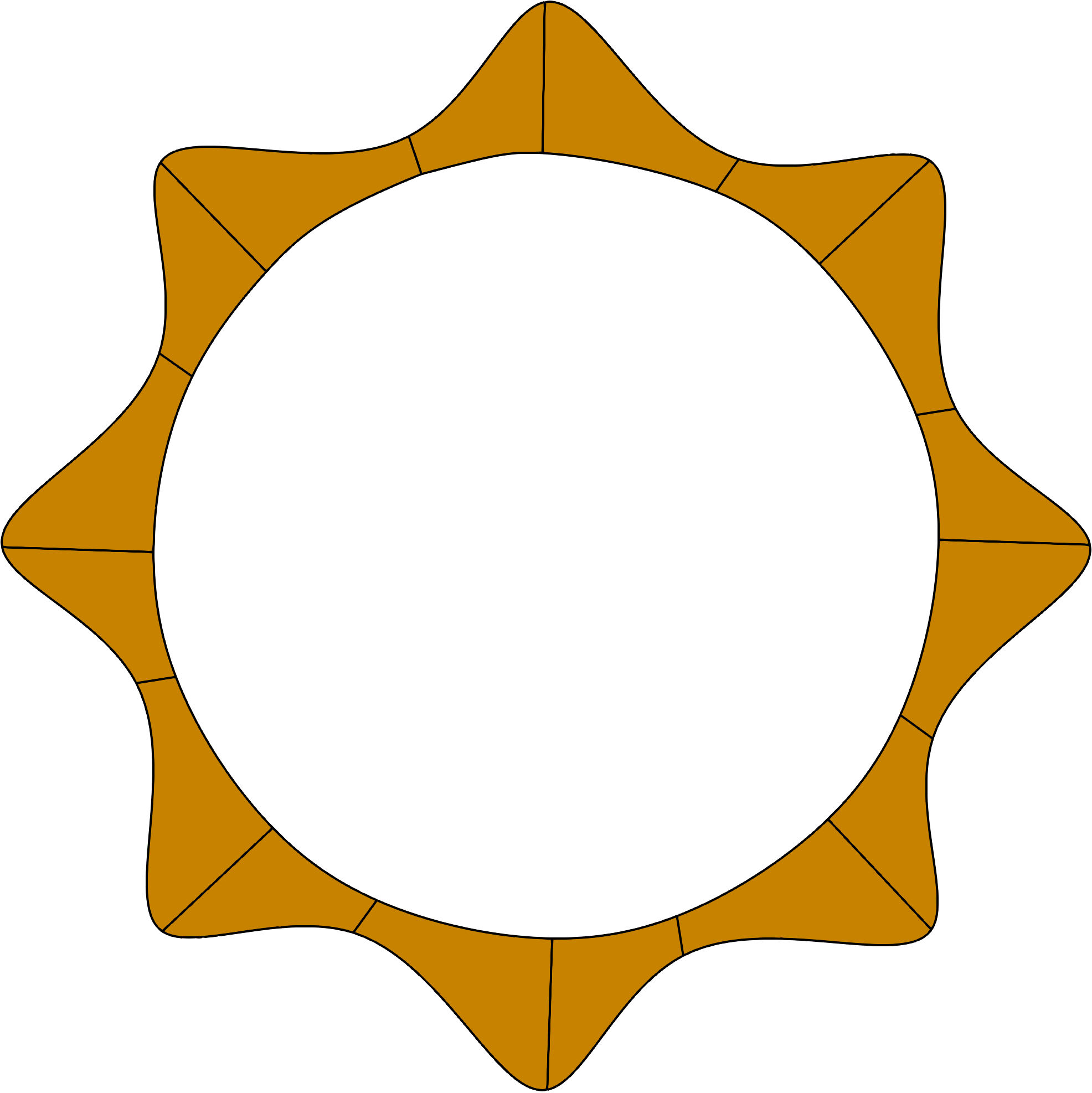}\\
(a) & (b)\\
\end{tabular}
\caption{Conformal layers obtained through an offset curve (a) and a circular target curve (b).}
\label{fig:flower_bl}
\end{figure}

Next, we solve Poisson's equation on a flower-shaped domain to test how the shape of a conformal layer influences the numerical solution. The input is a B-spline curve $\Gamma$ representing the flower boundary. We study two different constructions for the target curve $\tilde{\Gamma}$: (1) $\tilde{\Gamma}$ is constructed as an offset curve; and (2) $\tilde{\Gamma}$ is simply a circle; see Fig. \ref{fig:flower_bl}. In Case 1, the offset curve has a similar shape to $\Gamma$, but it generally has a different knot vector from $\Gamma$. Therefore, a loft surface is constructed as the conformal layer $\Omega^c$. As a result, its knot vector is a superset of those of both $\Gamma$ and $\tilde{\Gamma}$, leading to a dense mesh in $\Omega^c$. Moreover, $C^0$ continuities may be introduced to $\Omega^c$ due to the presence of cusps and removal of self-intersections in $\tilde{\Gamma}$. 

On the other hand, due to the flexibility of Case 2, $\tilde{\Gamma}$ can maintain the same knot vector as $\Gamma$ through Greville point projection. Thus, $\Omega^c$ is constructed simply as ruled surface between $\Gamma$ and $\tilde{\Gamma}$. This way, the mesh resolution of $\Omega^c$ is controlled by the input boundary curve.

\begin{figure}[htb]
\centering
\begin{tabular}{cc}
\includegraphics[width=0.46\columnwidth]{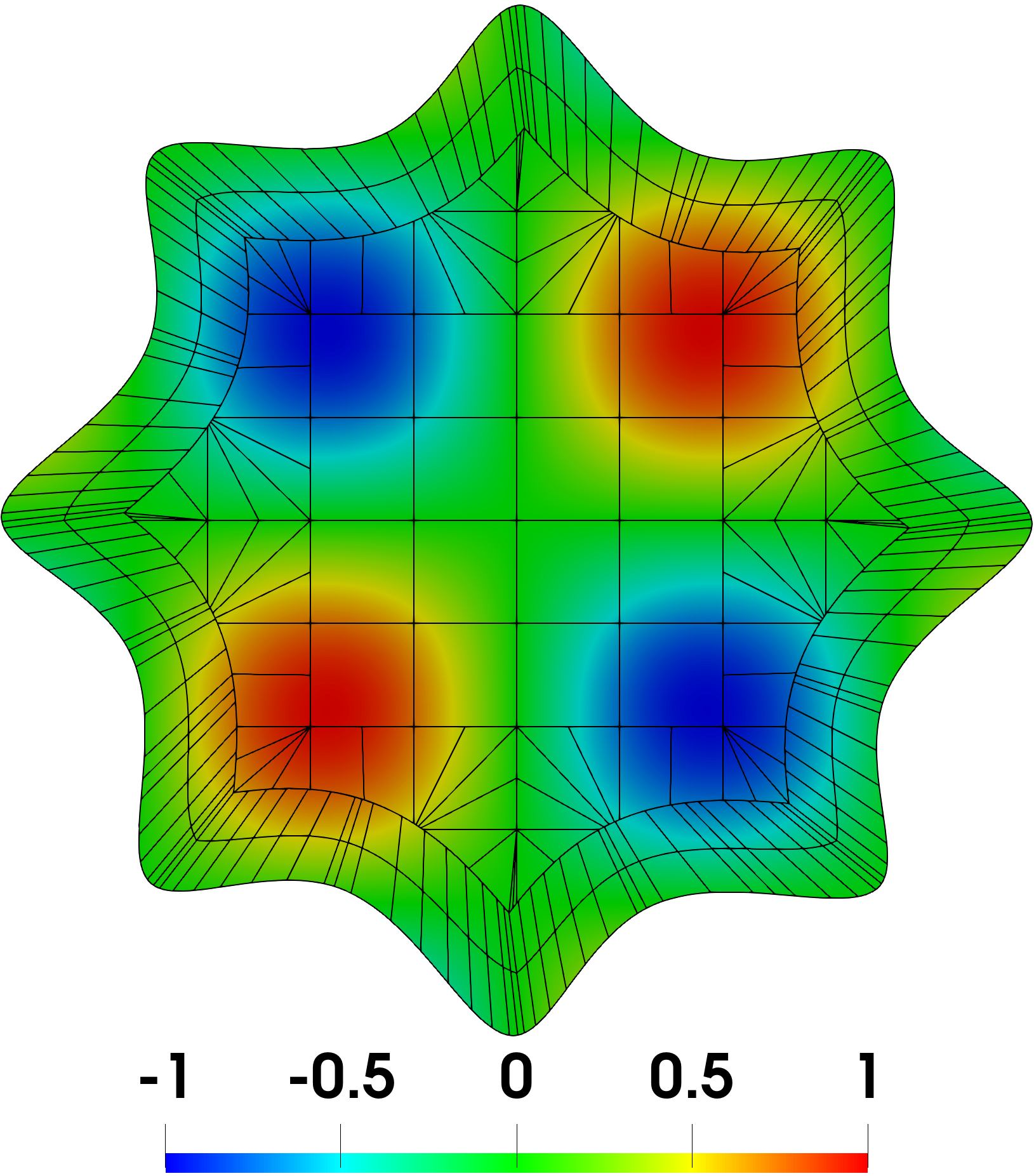}&\hspace{-2mm}
\includegraphics[width=0.46\columnwidth]{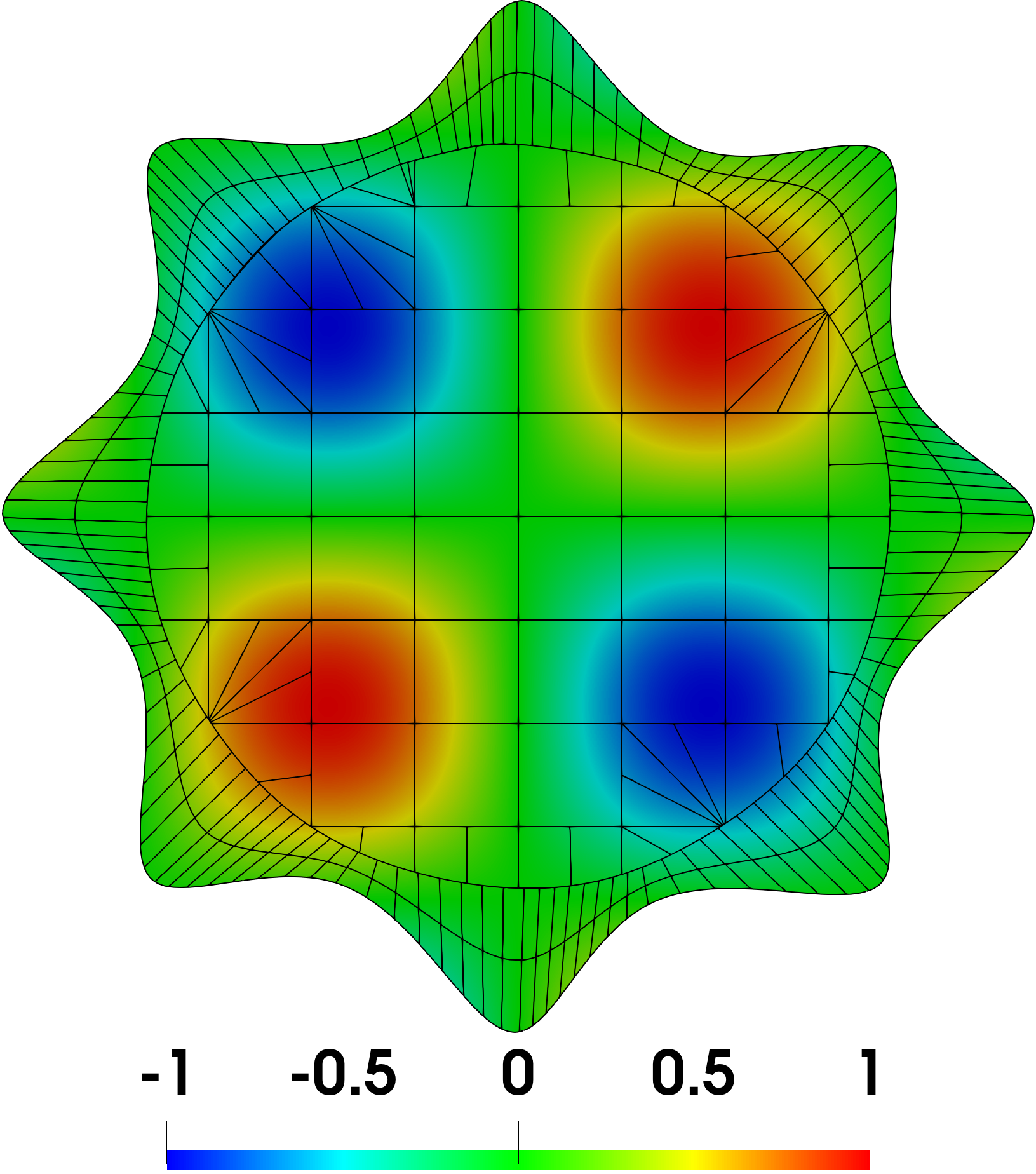}\\
\end{tabular}
\caption{Solutions on the initial B\'{e}zier mesh that has an offset-based conformal layer (a), and on the mesh that has a circle-based conformal layer (b). \rt{In the background mesh, integration cells of trimmed elements are used for visualization.}}
\label{fig:flower_u}
\end{figure}

\begin{figure}[htb]
\centering
\includegraphics[width=0.8\columnwidth]{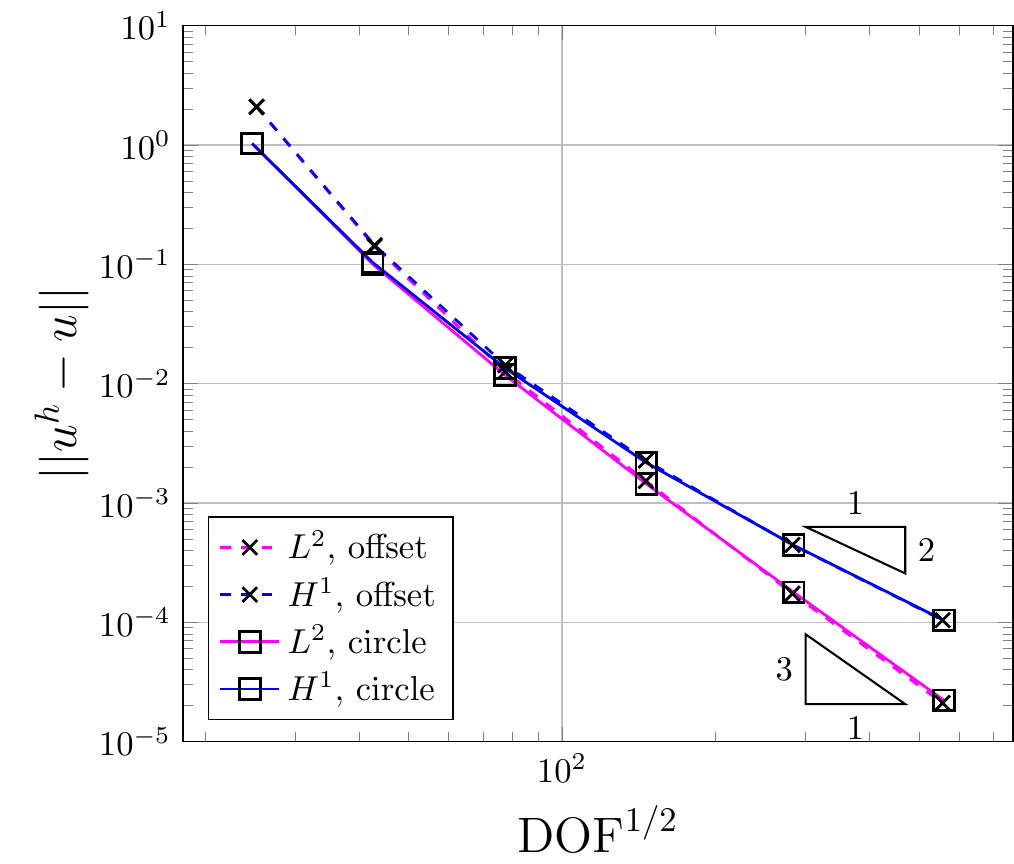}
\caption{Convergence plots with respect to $\mathrm{DOF}^{1/2}$ in the flower example. Note that ``offset" indicates the results using the offset target curve, whereas ``circle" corresponds to results using the circular target curve.}
\label{fig:flower_conv}
\end{figure}

With the conformal layers, we construct their corresponding IBCM representations and perform a convergence study using the following manufactured solution,
\begin{equation}
u(x,y) = \sin \left(\frac{\pi x}{R} \right) \sin \left( \frac{\pi y}{R} \right),
\end{equation}
where $R$ is a constant related to the ``radius" of the flower. The Dirichlet boundary condition is strongly imposed on the entire boundary. The solutions on the initial B\'{e}zier meshes are shown in Fig. \ref{fig:flower_u}, where we have refined the conformal layers to make them have a similar mesh resolution. Moreover,  we summarize the convergence plots in Fig. \ref{fig:flower_conv}. We observe that the results are almost identical in both geometric representations. In other words, numerical results are almost not influenced by the shape of the target curve, although further study is needed to come to a conclusion. Nonetheless, it indicates that we are not restricted to a specific choice for the target curve, which in turn provides further flexibility to the IBCM construction.

\subsection{Bimaterial disk}

\begin{figure}[htb]
\centering
\includegraphics[width=0.5\columnwidth]{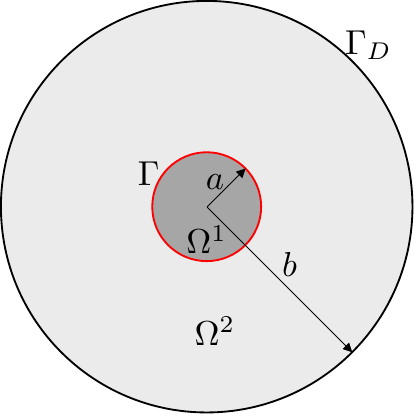}
\caption{A disk composed of two domains with different materials $\Omega^1$ and $\Omega^2$, where $\Gamma$ (the red circle) is the material interface.}
\label{fig:disk_set}
\end{figure}

This test is motivated by \cite{ref:sukumar01} and is aimed to test how IBCM resolves interface features when modeling multiple materials. We solve the linear elasticity problem on a disk composed of two different materials, the interior material $\Omega^1$ and the exterior material $\Omega^2$; see Fig. \ref{fig:disk_set}. Their corresponding material properties are given as $E_1=1$, $\nu_1=0.25$, and $E_2=10$, $\nu_2=0.3$, respectively. The Lam\'{e} constants ($\lambda_1$, $\mu_1$ and $\lambda_2$, $\mu_2$) are obtained according to Eq. \eqref{eq:lame}. The exact solutions are given in polar coordinates $(r,\theta)$. The displacement field is written as
\begin{equation}
\begin{aligned}
&u_r(r) = \left\{
\begin{array}{lll}
&\left[ \left( 1 - \frac{b^2}{a^2} \right) \alpha + \frac{b^2}{a^2}  \right] r , \quad & 0 \leq r \leq a, \\
&\left( r- \frac{b^2}{r} \right) \alpha + \frac{b^2}{r}, \quad & a < r \leq b,
\end{array}
\right. \\
&u_{\theta} = 0 ,\\
\end{aligned}
\label{eq:bidisk_u}
\end{equation}
where $a$ and $b$ are the radius of $\Omega^1$ and the outer radius of $\Omega^2$, respectively, and
\begin{equation}
\alpha = \frac{(\lambda_1 + \mu_1 + \mu_2) b^2}{(\lambda_2 +\mu_2)a^2 + (\lambda_1 + \mu_1)(b^2-a^2) + \mu_2 b^2} .
\end{equation}
We set $a=0.4$ and $b=2.0$ in the test. The radial ($\epsilon_{rr}$) and hoop ($\epsilon_{\theta \theta}$) strains are given by
\begin{equation}
\epsilon_{rr} (r) = \left\{
\begin{array}{lll}
&\left( 1- \frac{b^2}{a^2} \right) \alpha + \frac{b^2}{a^2}, \quad & 0 \leq r \leq a, \\
&\left( 1+ \frac{b^2}{r^2} \right) \alpha - \frac{b^2}{r^2}, \quad & a < r \leq b, \\
\end{array} 
\right.
\label{eq:bidisk_srr}
\end{equation}
and
\begin{equation}
\epsilon_{\theta\theta} (r) = \left\{
\begin{array}{lll}
&\left( 1- \frac{b^2}{a^2} \right) \alpha + \frac{b^2}{a^2}, \quad & 0 \leq r \leq a, \\
&\left( 1- \frac{b^2}{r^2} \right) \alpha + \frac{b^2}{r^2}, \quad & a < r \leq b. \\
\end{array} 
\right.
\label{eq:bidisk_stt}
\end{equation}
The radial ($\sigma_{rr}$) and hoop ($\sigma_{\theta \theta}$) stresses are
\begin{equation}
\begin{array}{lll}
&\sigma_{rr}(r) &= \lambda (\epsilon_{rr} + \epsilon_{\theta \theta}) + 2 \mu \epsilon_{rr}, \\
&\sigma_{\theta \theta}(r) &= \lambda (\epsilon_{rr} + \epsilon_{\theta \theta}) + 2 \mu \epsilon_{\theta\theta}, \\
\end{array}
\end{equation}
where \rt{$(\lambda,\mu)\in\{(\lambda_1,\mu_1),(\lambda_2,\mu_2)\}$}. The shear stress is zero everywhere. All solutions are axisymmetric in the sense that they are independent of $\theta$. We observe that in Eqs.~(\ref{eq:bidisk_u}, \ref{eq:bidisk_srr}, \ref{eq:bidisk_stt}), $u_r$ and $\epsilon_{\theta\theta}$ are continuous across the material interface, but $\epsilon_{rr}$ experiences a discontinuity. Moreover, both stresses ($\sigma_{rr}$ and $\sigma_{\theta\theta}$) are discontinuous due to different material parameters.

\begin{figure}[htb]
\centering
\begin{tabular}{cc}
\includegraphics[width=0.46\columnwidth]{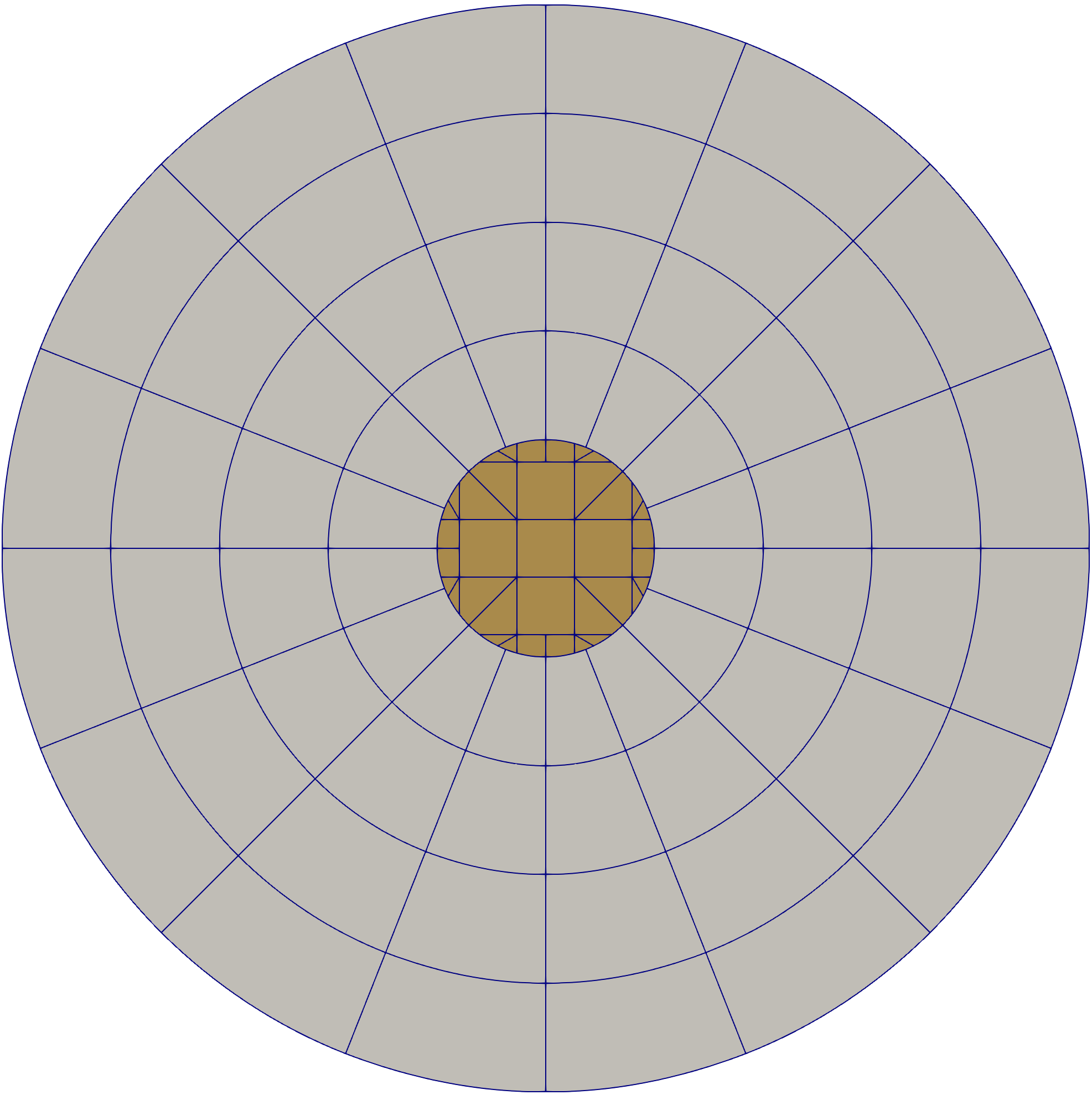}&\hspace{-2mm}
\includegraphics[width=0.46\columnwidth]{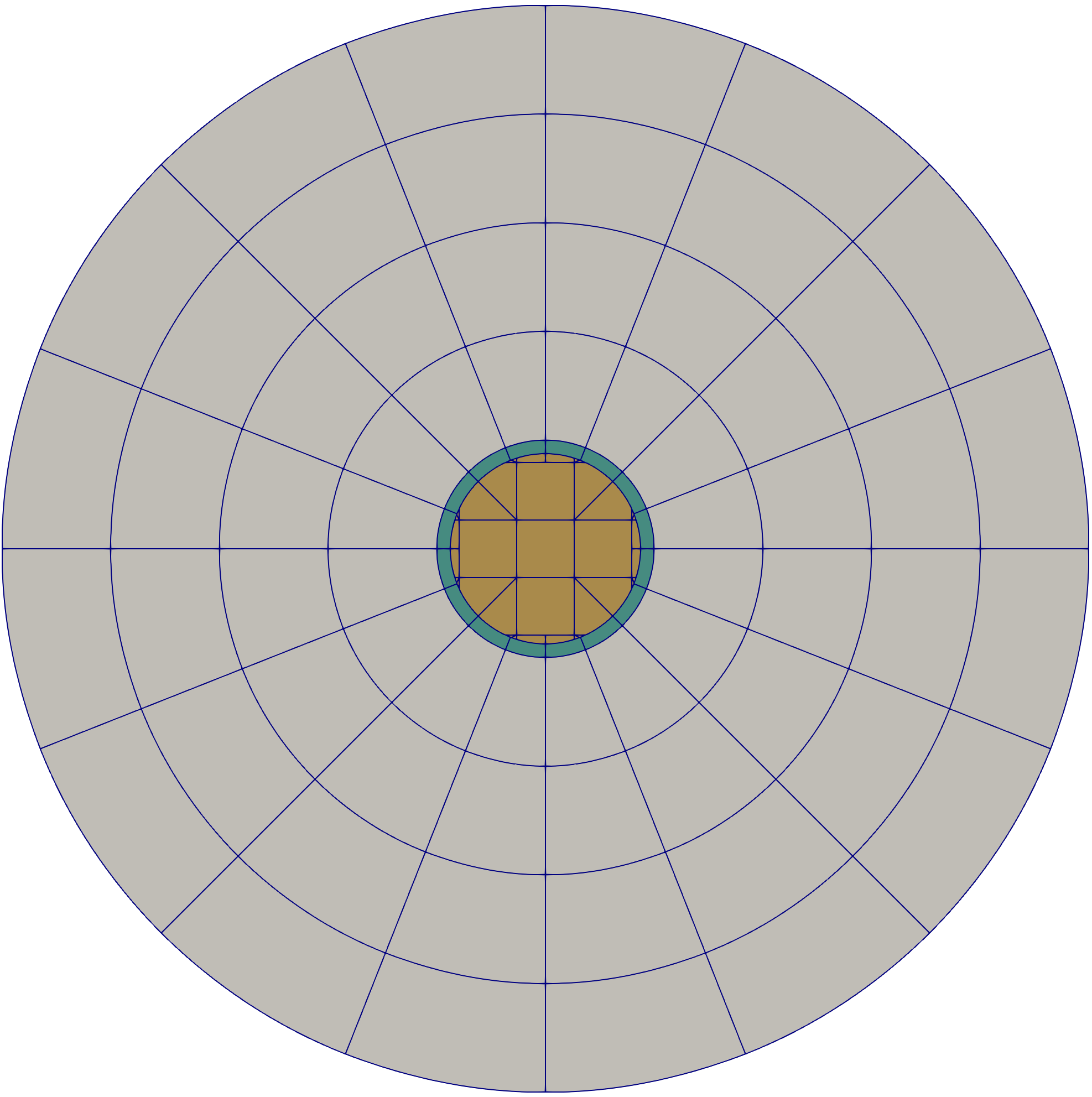}\\
(a) & (b)\\
\end{tabular}
\caption{Initial B\'ezier meshes of the bimaterial disk via union (a) and IBCM (b), where each patch is represented by a different color. Particularly in (b), the interior material $\Omega^1$ is represented by two patches: the green annulus and the yellow background mesh. \rt{In the background mesh, integration cells of trimmed elements are used for visualization.}}
\label{fig:disk_geom}
\end{figure}

\begin{figure}[htb]
\centering
\begin{tabular}{c}
\includegraphics[width=0.8\columnwidth]{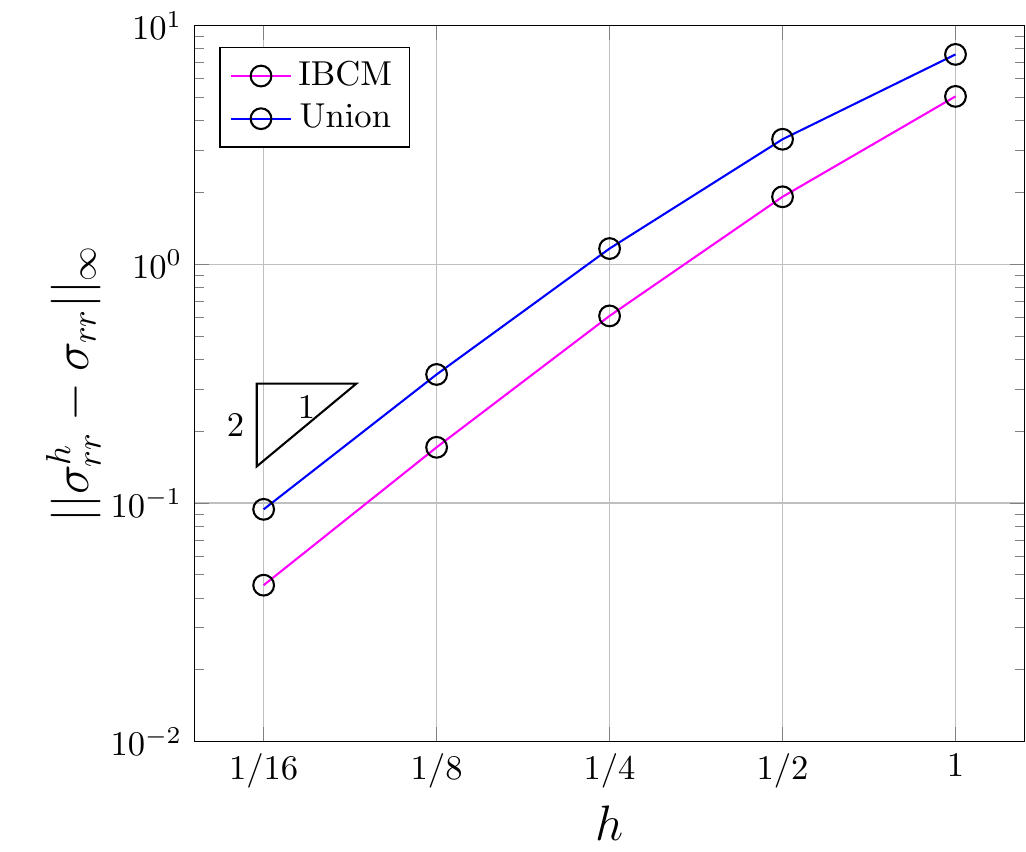}\\
(a)\\
\includegraphics[width=0.8\columnwidth]{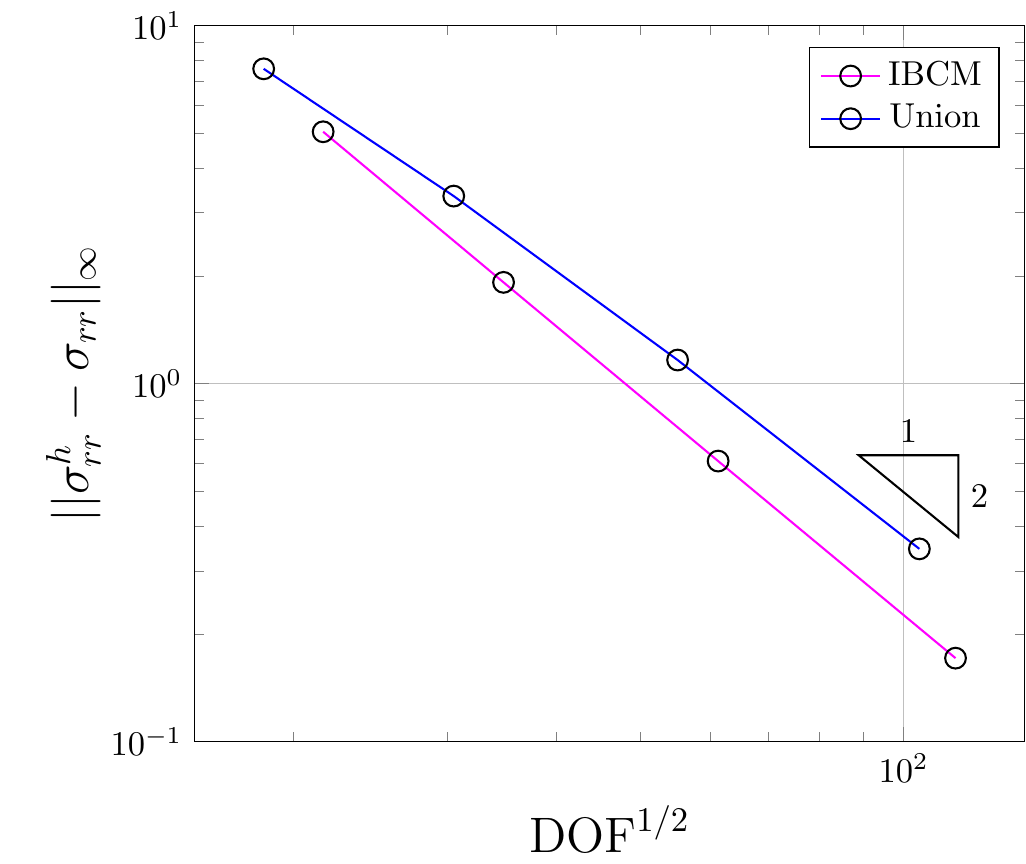}\\
(b)\\
\end{tabular}
\caption{Convergence plots with respect to $h$ (a) and $\mathrm{DOF}^{1/2}$ (b) in the bimaterial disk test.}
\label{fig:disk_conv}
\end{figure}

Two geometric constructions are studied: a union of two overlapping patches (as in~\cite{ref:wei19u}) and an IBCM construction; see Fig. \ref{fig:disk_geom}. In the union construction, the grey annulus sits on top of the yellow background patch, which represent the material domains $\Omega^2$ and $\Omega^1$, respectively. On the other hand, the interior material $\Omega^1$ in an IBCM representation is represented by two patches, the green conformal layer and the yellow background mesh. The mesh of $\Omega^2$ remains the same. According to Eq. \eqref{eq:bidisk_u}, the exact displacement, $u_r(b)=b$, is applied on the outer boundary $\Gamma_D$; see Fig. \ref{fig:disk_set}. We summarize the convergence plots in Fig. \ref{fig:disk_conv}. Again, we observe the same expected convergence rates but improved accuracy per DOF in IBCM.

\subsection{Spanner}

\begin{figure}[htb]
\centering
\begin{tabular}{c}
\includegraphics[width=0.8\columnwidth]{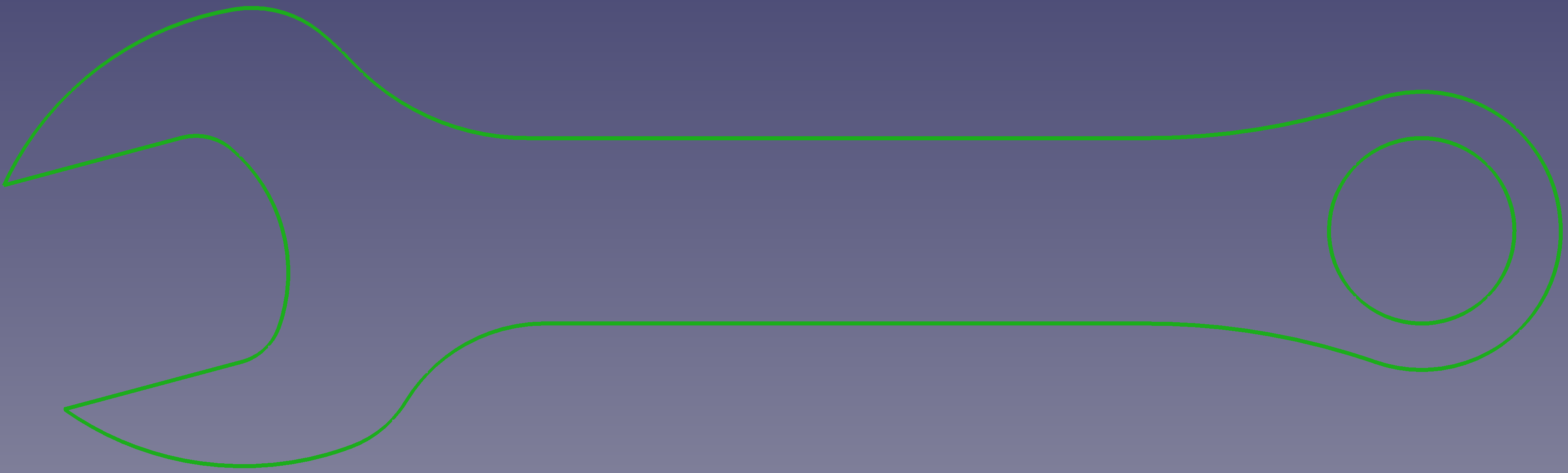}\\
(a)\\
\includegraphics[width=0.8\columnwidth]{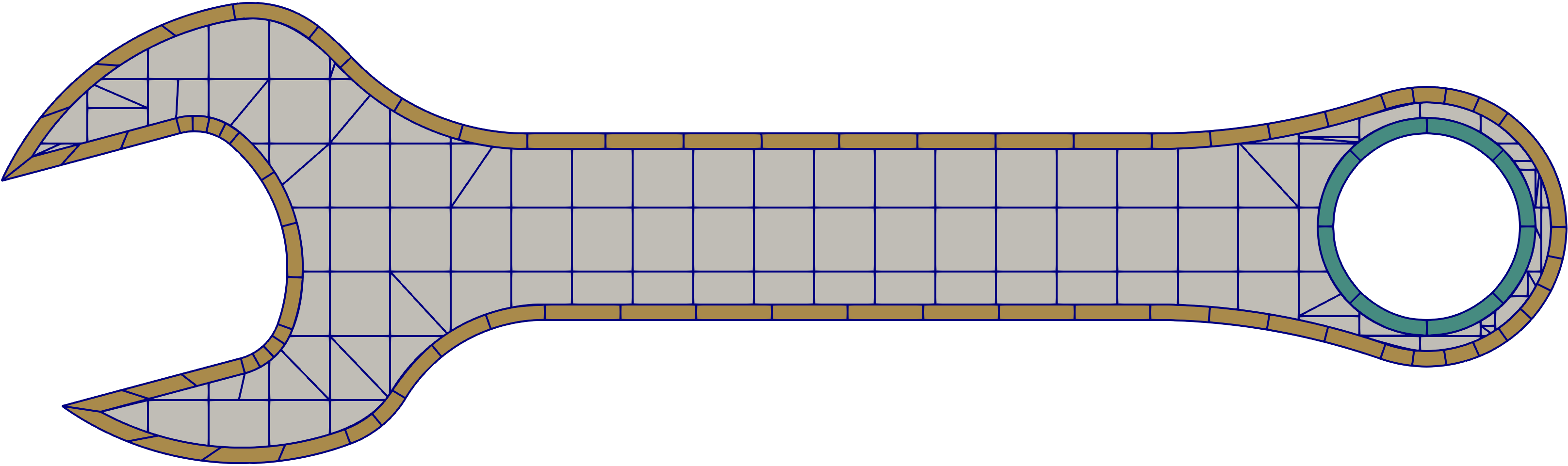}\\
(b)\\
\includegraphics[width=0.8\columnwidth]{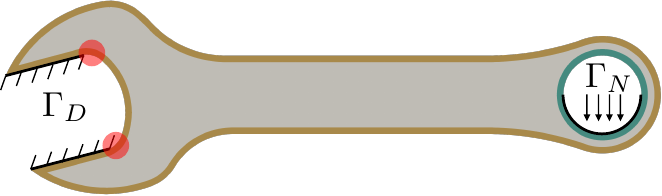}\\
(c)\\
\end{tabular}
\caption{The input B-rep of spanner (a), the corresponding IBCM representation (b), and boundary conditions (c). \rt{In the background mesh of (b), integration cells of trimmed elements are used for visualization.} In (c), stress singularities are expected in the red-spotted regions due to the sudden change from the clamped displacement condition to the traction-free condition.}
\label{fig:spanner_input}
\end{figure}

\begin{figure*}[htb]
\centering
\begin{tabular}{cc}
\includegraphics[width=0.6\textwidth]{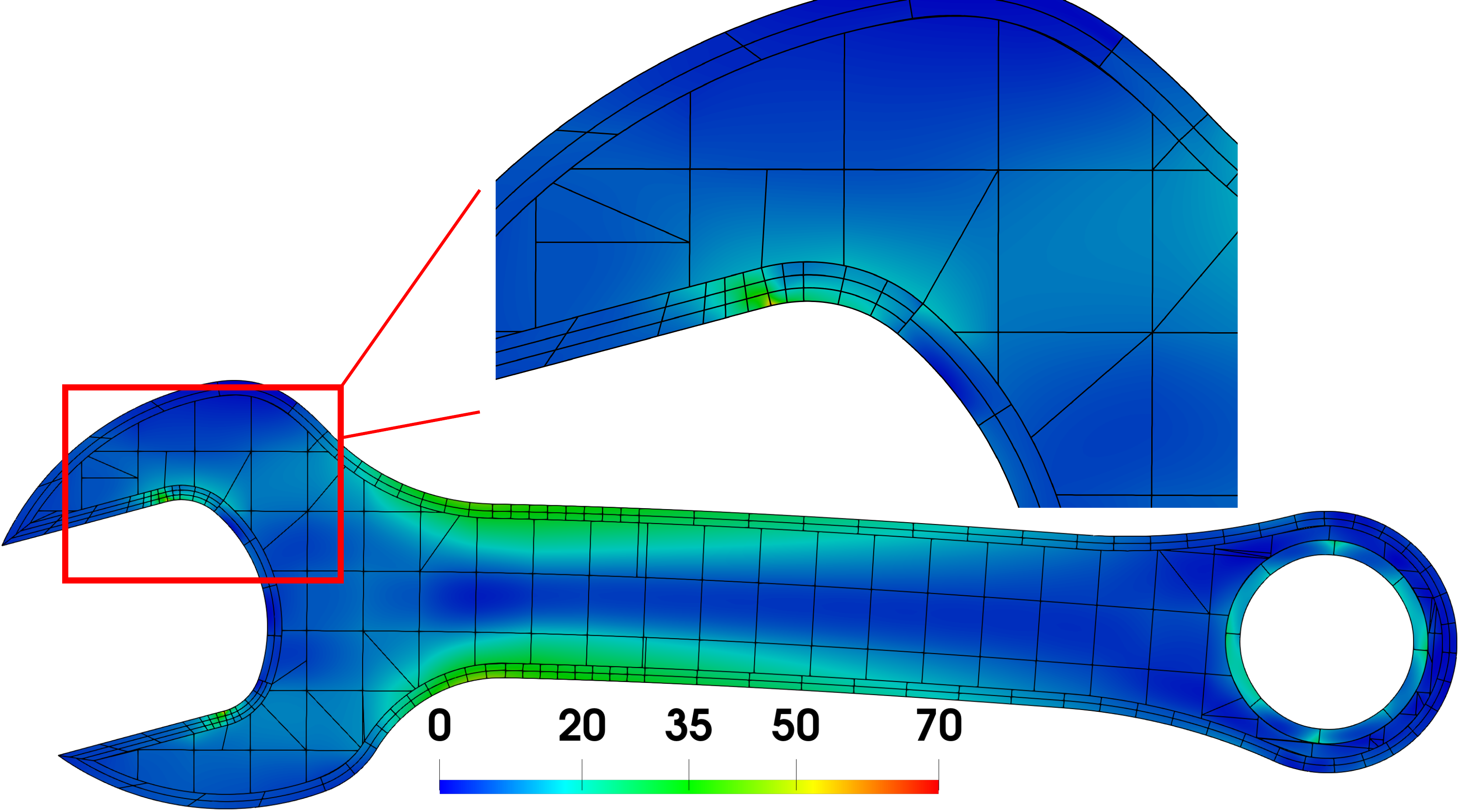} &\hspace{-8mm}
\includegraphics[width=0.4\textwidth]{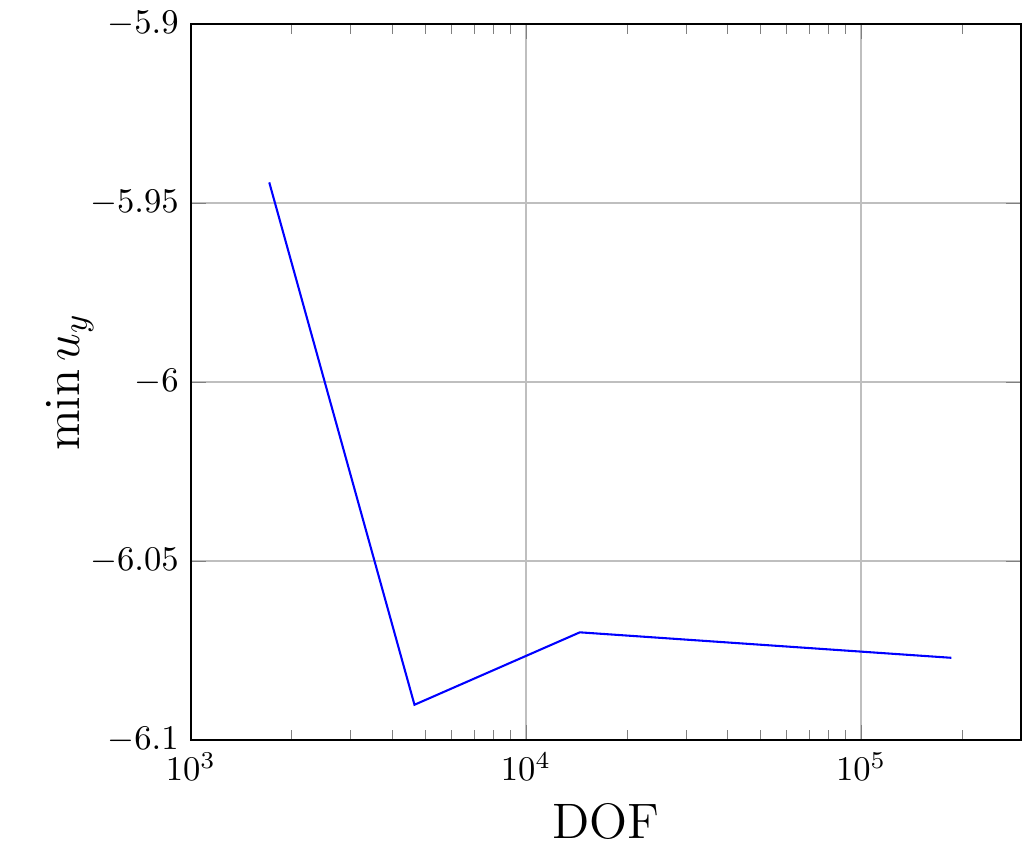}\\
(a) & (b)\\
\end{tabular}
\caption{\rt{Results of the spanner model. (a) The von-Mises stress on the IBCM representation of the spanner, and (b) the convergence plot of the maximum downwards displacement ($-y$ direction). In (a), the mesh resolution is (manually) adapted to large stresses. In the zoomed-in figure of (a), some part of the mesh in the conformal layer is not conformal. In the background mesh of (a), integration cells of trimmed elements are used for visualization.}}
\label{fig:spanner_stress}
\end{figure*}

We next present a spanner model to demonstrate how to use IBCM to represent complex geometries. We will also show that it is flexible for IBCM to deal with Dirichlet boundary conditions as well as to control the mesh resolution around areas of interest. As shown in Fig. \ref{fig:spanner_input}(a), the input B-rep of the spanner is composed of two loops, each represented by a set of NURBS curves. A detailed description of the geometry data is given in \cite{ref:stanford19}. We first generate an offset curve for each loop, which, in this particular case, has a nearly identical parameterization as the input loop. The conformal layer is then immediately available as a loft surface; see Fig. \ref{fig:spanner_input}(b). The remaining procedure of IBCM follows the boundary-type construction.

We proceed to solve the linear elasticity problem on the IBCM representation, with material constants $E=10^{5}$ and $\nu=0.3$. As shown in Fig.~\ref{fig:spanner_input}(c), we clamp part of the jaw ($\bm{u}=\bm{0}$ on $\Gamma_D$), and we apply traction downwards on the lower half of the handle with a magnitude of $1.0$. All the remaining part of the boundary is traction free, i.e., homogeneous Neumann boundary condition. This way, singularities (in stresses) are expected in the transition regions from the clamped displacement condition to the traction-free condition; see the red-spotted regions in Fig.~\ref{fig:spanner_input}(c). We further ``locally" refine the mesh to capture such features. The von Mises stress is shown in Fig.~\ref{fig:spanner_stress}(a). We observe that with the conformal layer, we can easily adapt the mesh resolution to the solution features. Moreover, when a conformal layer has a multi-patch representation, meshes in different patches do not need to be conformal across patch interfaces; see the zoomed picture in Fig.~\ref{fig:spanner_stress}(a). Non-conformal patch interfaces are coupled with Nitsche's method and are treated as a special case of the general overlapping construction. This way, mesh refinement is localized to the patch level.

\rt{Moreover, we check the convergence of the maximum downwards displacement (the $-y$ direction), or $\min\,u_y$, with a series of globally refined meshes. The convergence plot is shown in Fig. \ref{fig:spanner_stress}(b) with respect to DOF, where we observe a convergent result. Here we study the convergence in displacements rather than stresses because in this example, stresses exhibit singularity due to the boundary conditions, which we do on purpose to show how IBCM can be used to flexibly capture such a solution feature.}

\subsection{Fiber-reinforced composite}

\begin{figure}[htb]
\centering
\includegraphics[width=0.7\columnwidth]{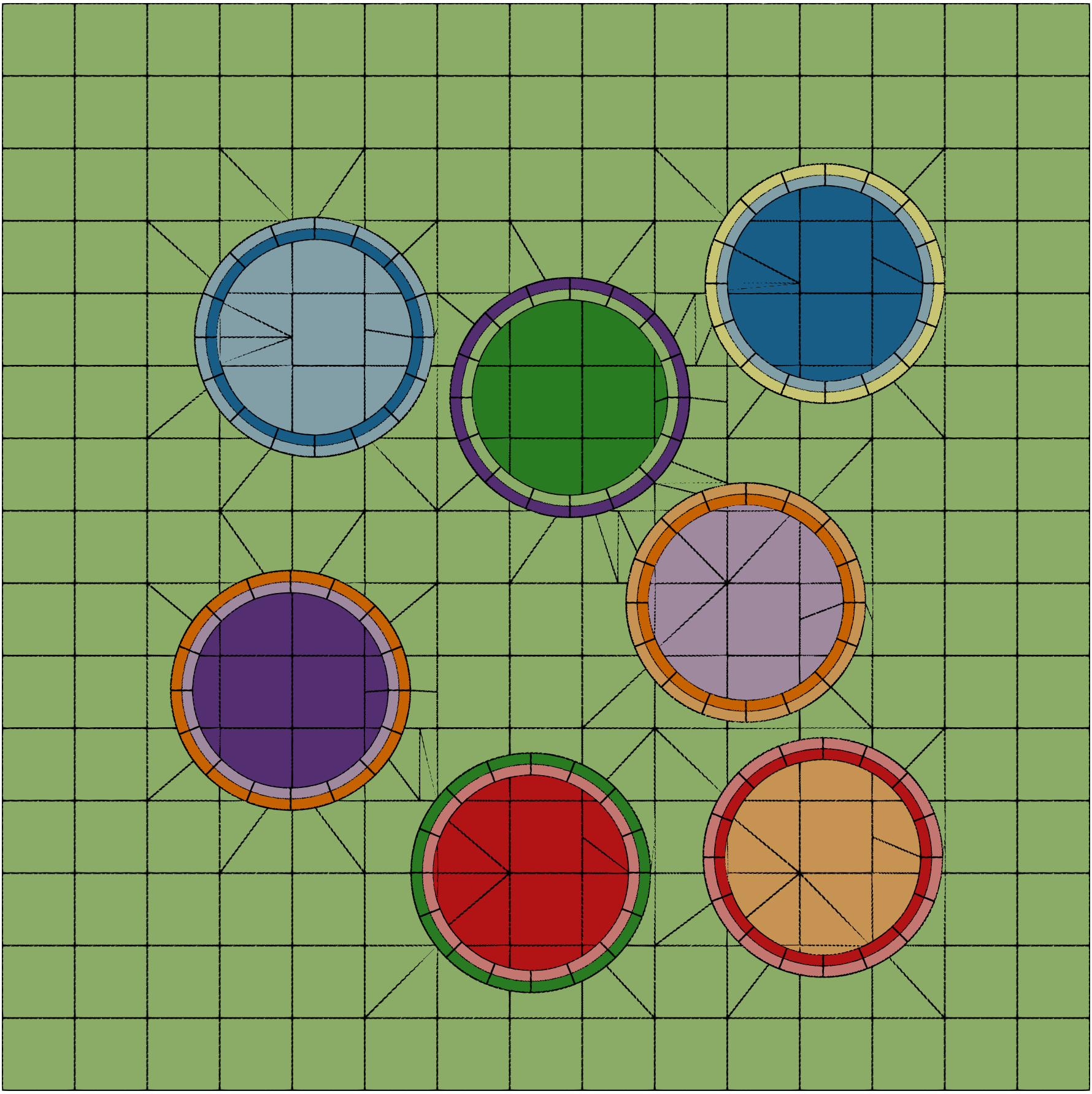}
\caption{The initial IBCM representation of the fiber-reinforced composite, where the circular fibers are randomly positioned and different patches are represented by different colors. \rt{In the background mesh, integration cells of trimmed elements are used for visualization.}}
\label{fig:composite_mesh}
\end{figure}

Inspired by \cite{ref:herraez16}, we perform stress analysis on a fiber-reinforced composite. The test is \rt{mainly} aimed to show the capability of IBCM to represent complex material interfaces with conformal discretizations. A typical cross section is studied under the plane strain assumption; see Fig. \ref{fig:composite_mesh}. The matrix material is modeled as a square that occupies the domain $[-5,5]\times[-5,5]$. Then multiple circular fibers of unit radii are randomly positioned in the matrix, provided that none of them overlaps with another. For each circular material interface, a pair of conformal annuli is obtained following the interface-type construction, where one has the material property of fibers and the other has the property of the matrix. Every fiber has an independent background mesh, whereas the background mesh of the matrix is cut by multiple annuli. Combining all the annuli and the active parts of background meshes yields the IBCM representation for the fiber-reinforced composite.

\begin{figure}[htb]
\centering
\begin{tabular}{c}
\includegraphics[width=0.8\columnwidth]{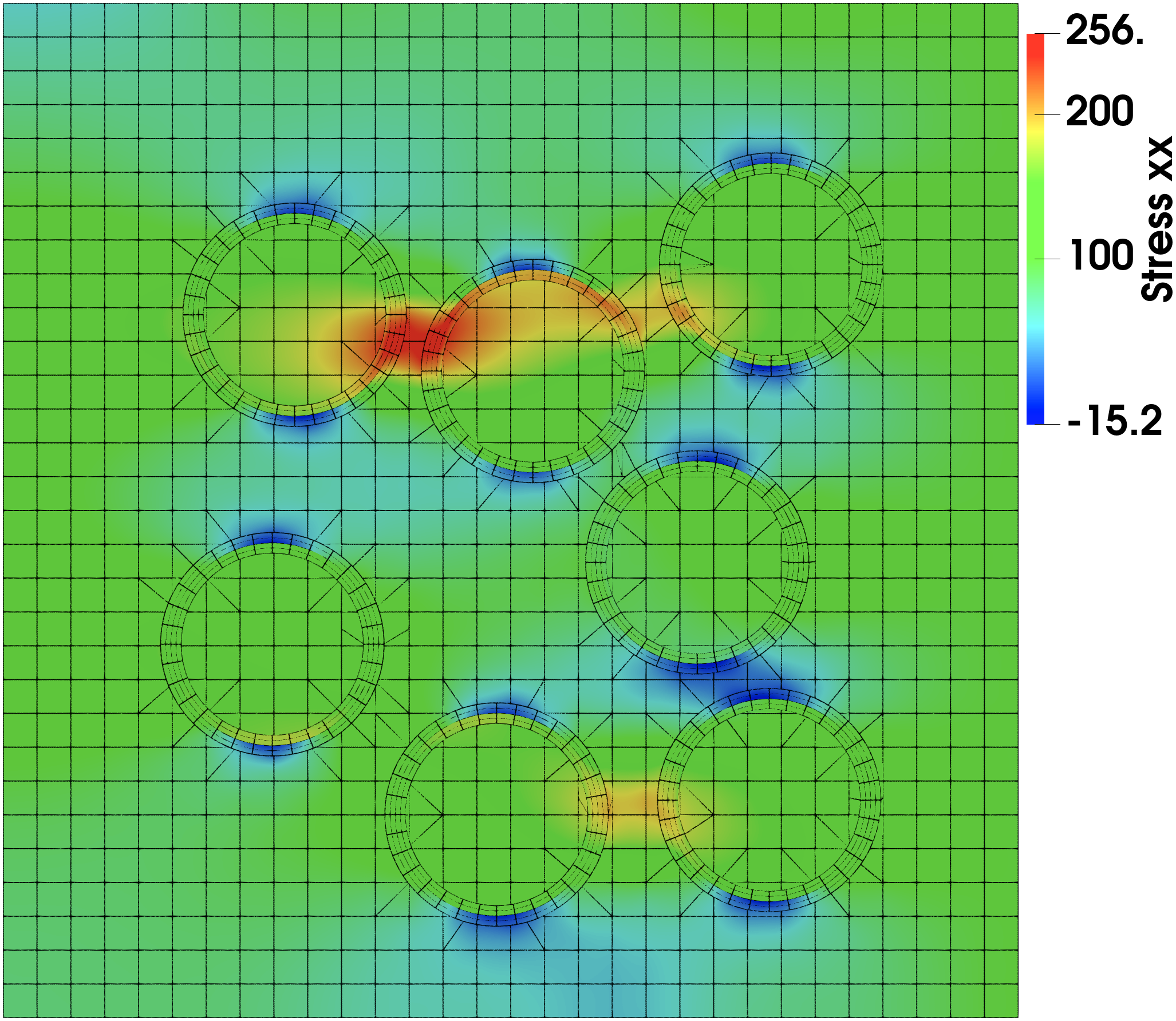}\\
(a)\\
\hspace*{-1cm}\includegraphics[width=0.9\columnwidth]{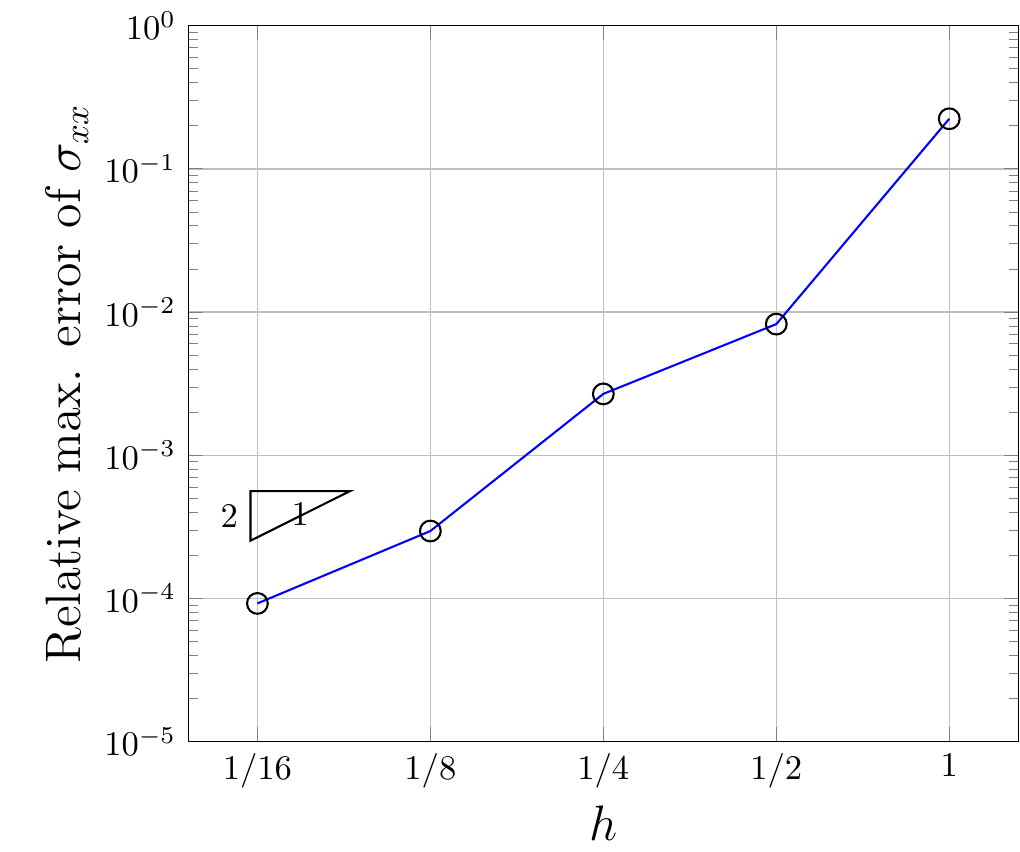}\\
(b)\\
\end{tabular}
\caption{\rt{Results of the fiber-reinforced composite. (a) The stress $\sigma_{xx}$ computed on the mesh after one refinement, and (b) the convergence plot of $\sigma_{xx}$. In the background mesh of (a), integration cells of trimmed elements are used for visualization.}}
\label{fig:composite_stress}
\end{figure}

We solve the linear elasticity problem on the resulting geometric model. Uniform tension $T_x=100$ is imposed on the right boundary of the matrix, whereas $u_x=0$ and $u_y=0$ are imposed on the left and bottom boundaries, respectively. Homogeneous Neumann boundary conditions are applied elsewhere. Material constants are given as follows: $E_{\mathrm{fiber}} =100$, $\nu_{\mathrm{fiber}}=0.33$, and $E_{\mathrm{matrix}} =1$, $\nu_{\mathrm{matrix}}=0.3$. As a result of discontinuous material data, stress discontinuities are expected across material interfaces. 

\rt{Starting from the input mesh in Fig. \ref{fig:composite_mesh}, we obtain a series of globally refined meshes for a convergence study. The distribution of the stress $\sigma_{xx}$ is shown in Fig. \ref{fig:composite_stress}(a), which is computed with the mesh after one refinement. As expected, large stresses occur in fibers as well as in the narrow bands between fibers. The convergence plot is shown in Fig. \ref{fig:composite_stress}(b). In each mesh, we compute the maximum $\sigma_{xx}$ (the $x-x$ stress) among all the sampled points, and we check the convergence of this quantity. A relative error is evaluated using an overkill solution, which is obtained with the mesh after five times of refinement. We observe that the error converges roughly at an expected rate of two.}

\rt{
\subsection{Advection-diffusion problem}

\newcommand{\po}{\partial\Omega}

As the last example, we push a step forward to study IBCM beyond the elliptic problems. It is motivated by the singularly perturbed model problem \cite{ref:hemker96}, but here we will adopt a simplified boundary condition. The underlying partial differential equation (PDE) is an advection-diffusion equation, which is fundamentally different from those in the previous examples because of the advection term. It states as follows. Given $f:\,\Omega\to\mathbb{R}$ and $g_D:\,\Gamma_D\to\mathbb{R}$, find $u:\,\Omega\to\mathbb{R}$ such that
\begin{equation}
\left\{
\begin{aligned}
-\epsilon \Delta u + \bm{a}\cdot \nabla u = f \quad & \text{in} \quad \Omega=\Omega^t \cup \Omega_a^b, \\
u^t - u^b = 0 \quad & \text{on} \quad \Gamma = \partial \Omega^t \cap \Omega^b, \\
\nabla u^t \cdot \bm{n}^t + \nabla u^b \cdot \bm{n}^b = 0 \quad & \text{on} \quad \Gamma= \partial \Omega^t \cap \Omega^b, \\
u = g_D \quad & \text{on} \quad \Gamma_D\equiv \po, \\
\end{aligned}
\right.
\label{eq:ad_strong}
\end{equation}
where $\epsilon\in\mathbb{R}$ and $\bm{a}\in\mathbb{R}^2$ are a diffusivity constant and a velocity constant, respectively, and the other notations (including those in the following) are the same as those in Eq. \eqref{eq:poisson_strong}. Note that we only consider the Dirichlet boundary condition in this problem. 

We follow \cite{ref:buffa06} to obtain the corresponding weak formulation. $\po$ is divided into the outflow boundary $\po^+ = \{\bm{x}\in\po: \bm{a}\cdot \bm{n} \geq 0\}$ and the inflow boundary $\po^- = \{\bm{x}\in\po: \bm{a}\cdot \bm{n} < 0\}$, where $\bm{n}$ is the outward normal of $\po$. The discrete weak formulation states as follows: Find $u_h \in \mathcal{V}_h^{g_D}$ such that 
\begin{equation}
a_h(u_h,v_h) = l (v_h), \quad \forall v_h \in \mathcal{V}_h^0,
\end{equation}
where
\begin{equation}
\begin{aligned}
a_h(u_h,v_h) &= \epsilon \int_{\Omega} \nabla u_h \cdot \nabla v_h - \int_{\Omega} u_h \, \bm{a} \cdot \nabla v_h \\
&- \epsilon \int_{\Gamma} \langle\nabla u_h \cdot \rt{\bm{n}^t} \rangle [v_h] - \epsilon \int_{\Gamma} \langle\nabla v_h \cdot \bm{n}^t \rangle [u_h] \\
&+ \epsilon \beta (h_t^{-1} + h_b^{-1}) \int_{\Gamma} [u_h] [v_h], \\
&+ \int_{\Gamma} \mathrm{up} (u_h) \, [v_h] \, \bm{a} \cdot \bm{n}^t \\
&+ \int_{\po^+} u_h \, v_h \,  \bm{a} \cdot \bm{n} +  \epsilon \beta h^{-1} \int_{\po} u_h \, v_h
\end{aligned}
\label{eq:ad_weak}
\end{equation}
and
\begin{equation}
l(v_h) = \int_{\Omega} f \, v_h + \epsilon \beta h^{-1} \int_{\partial\Omega} g_D \, v_h - \int_{\partial\Omega^-} g_D \, v_h \, \bm{a}\cdot\bm{n},
\end{equation}
where
\begin{equation}
h(\bm{x}) =\left\{
\begin{aligned}
&h_t \quad \text{if } \bm{x} \in \po \cap \po^t,\\
&h_b \quad \text{if } \bm{x} \in \po \cap \po^b,
\end{aligned}
\right.
\end{equation}
and
\begin{equation}
\mathrm{up}(u_h(\bm{x})) =\left\{
\begin{aligned}
&u^t(\bm{x}) \quad \text{if } \bm{a} \cdot \bm{n^t}(\bm{x}) \geq 0,\\
&u^b(\bm{x}) \quad \text{otherwise}.
\end{aligned}
\right.
\label{eq:upwind}
\end{equation}
Eq. \eqref{eq:upwind} represents an \emph{upwind} scheme that takes certain quantity of interest from the upwind side. Note that the normal of the interface $\Gamma$ is $\bm{n}^t$, i.e., the outward normal of the top domain $\Omega^t$.

\begin{figure}[htb]
\centering
\includegraphics[width=0.8\columnwidth]{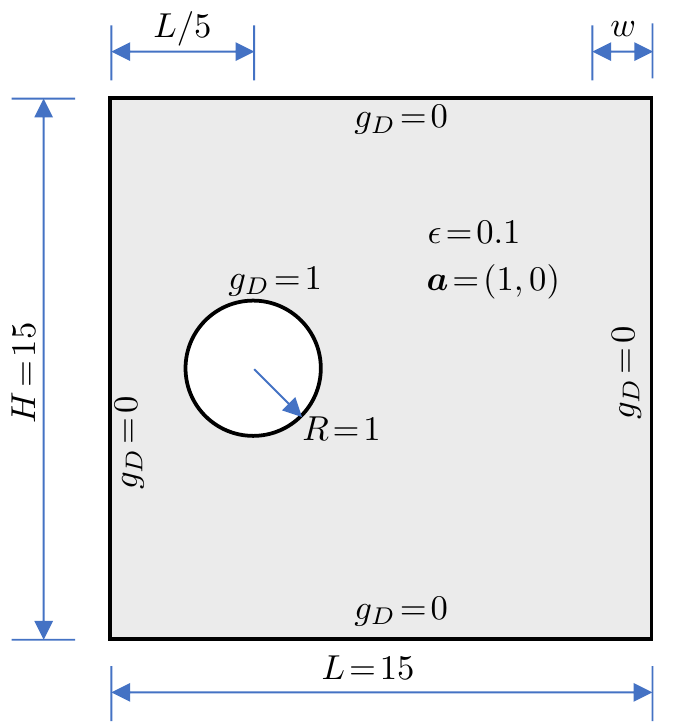}
\caption{The geometry and problem setting for the advection-diffusion problem.}
\label{fig:ad_plate_hole_geom}
\end{figure}

As shown in Fig. \ref{fig:ad_plate_hole_geom}, the geometry of interest is a plate with a hole. The plate has the dimension $L\times H$ where $L=H=15$. The hole is centered at $(L/5,H/2)$ with a radius of $R=1$. The diffusivity is set to be $\epsilon=0.1$ to make the problem advection-dominated, and a constant rightward velocity of $\bm{a}=(1,0)$ is adopted. According to Eq. \eqref{eq:ad_weak}, we \emph{weakly} impose the Dirichlet boundary condition on all the boundaries. More specifically, we impose $g_D=1$ on the hole and $g_D=0$ on the other boundaries. The boundary-layer phenomenon is expected around the hole and the right boundary, where dense meshes are needed to resolve it. It is straightforward to ``locally" refine the mesh near the right boundary by adding more knot lines. Therefore, we focus on the treatment of the hole, where we compare the performance of trimming and IBCM.

\begin{figure}[htb]
\centering
\begin{tabular}{cc}
\includegraphics[width=0.47\columnwidth]{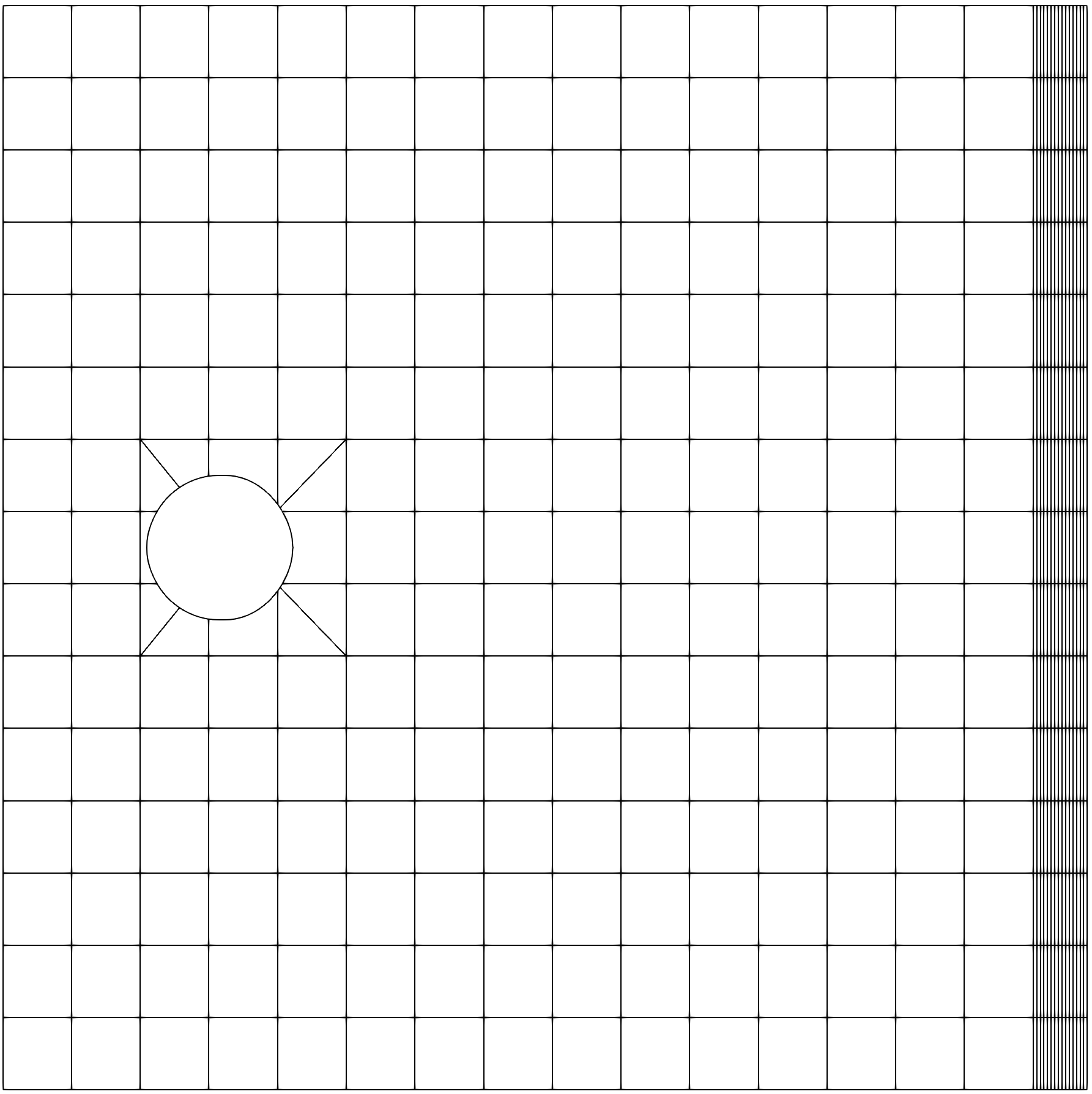}&\hspace{-2mm}
\includegraphics[width=0.47\columnwidth]{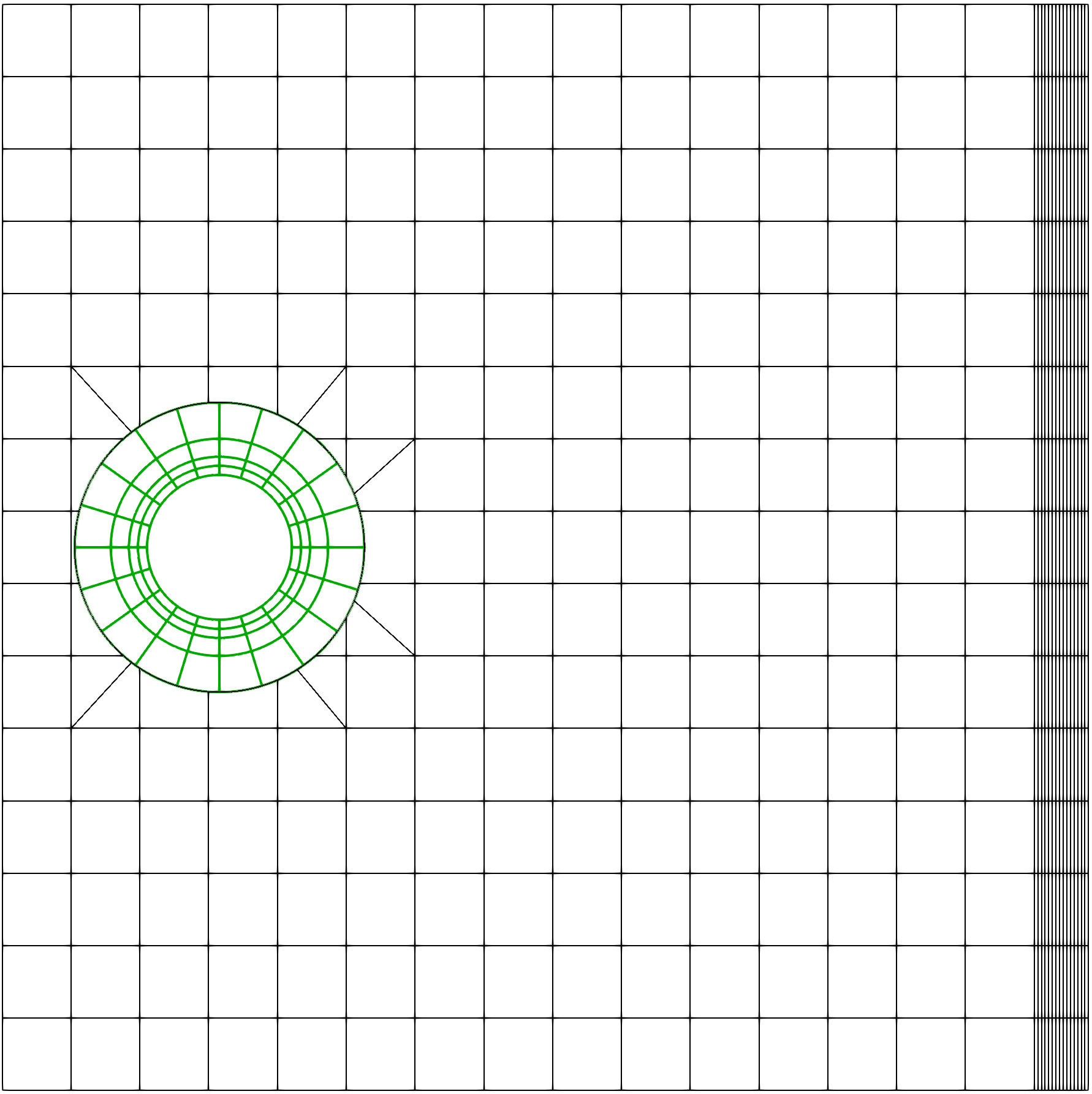}\\
(a) & (b)\\
\end{tabular}
\caption{The initial B\'ezier meshes of the plate with a hole via trimming (a) and IBCM (b). In the background mesh, integration cells of trimmed elements are used for visualization.}
\label{fig:ad_mesh}
\end{figure}

We first study the geometric representation obtained via trimming. The initial B\'ezier mesh is shown in Fig. \ref{fig:ad_mesh}(a), where a dense mesh of width $w=0.05L$ (see Fig. \ref{fig:ad_plate_hole_geom}) is adopted near the right boundary. We solve the advection-diffusion problem on a series of globally refined meshes. The results are shown in Fig. \ref{fig:ad_trim_result}. We observe that the high-quality solution is obtained only when the mesh size is sufficiently small with respect to the boundary layer; see Fig. \ref{fig:ad_trim_result}(f).

\begin{figure}[htb]
\centering
\begin{tabular}{cc}
\includegraphics[width=0.47\columnwidth]{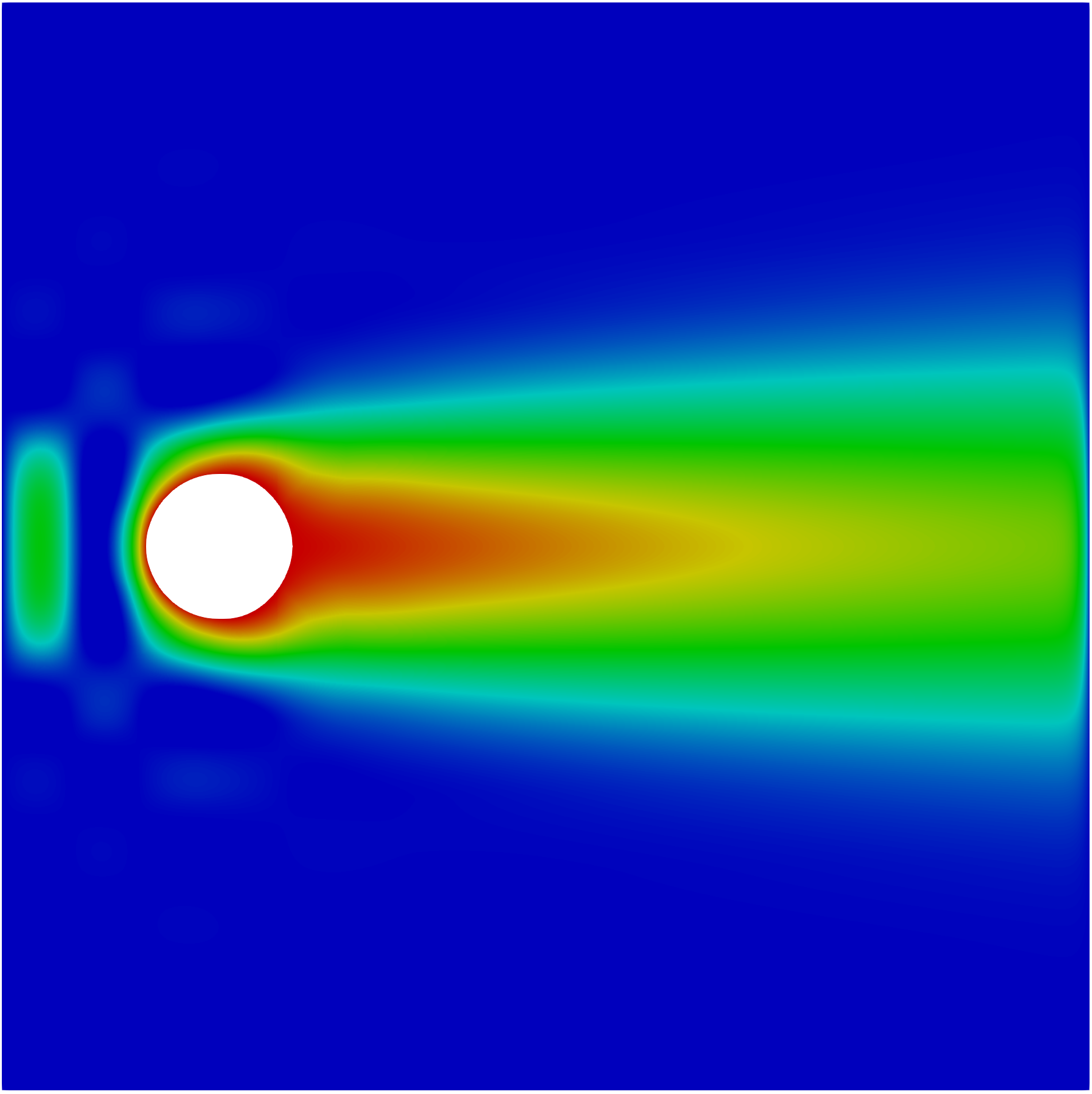}&\hspace{-2mm}
\includegraphics[width=0.47\columnwidth]{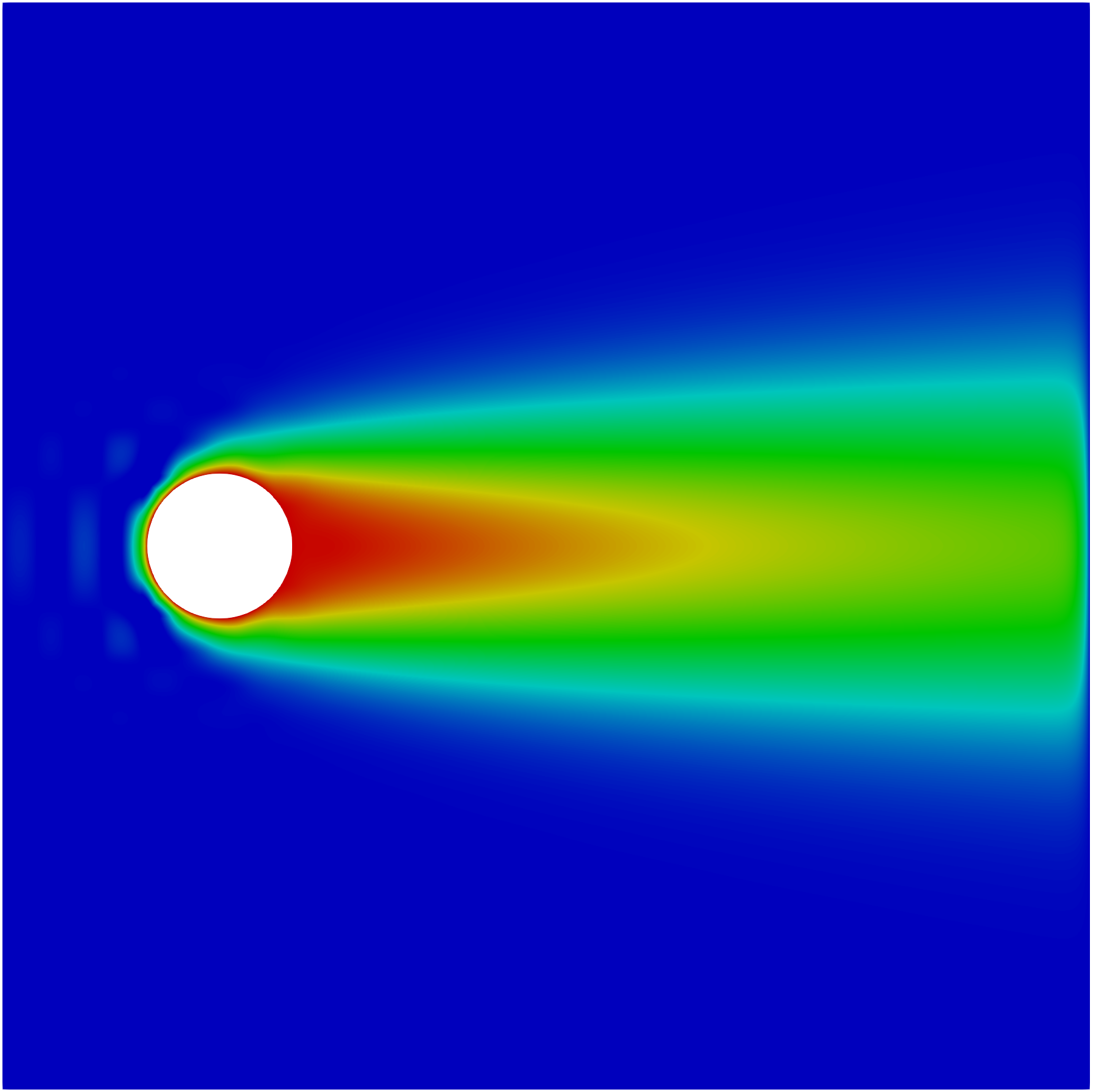}\\
(a) Mesh 0 (DOF: 544) & (b) Mesh 1 (DOF: 1984)\\
\includegraphics[width=0.47\columnwidth]{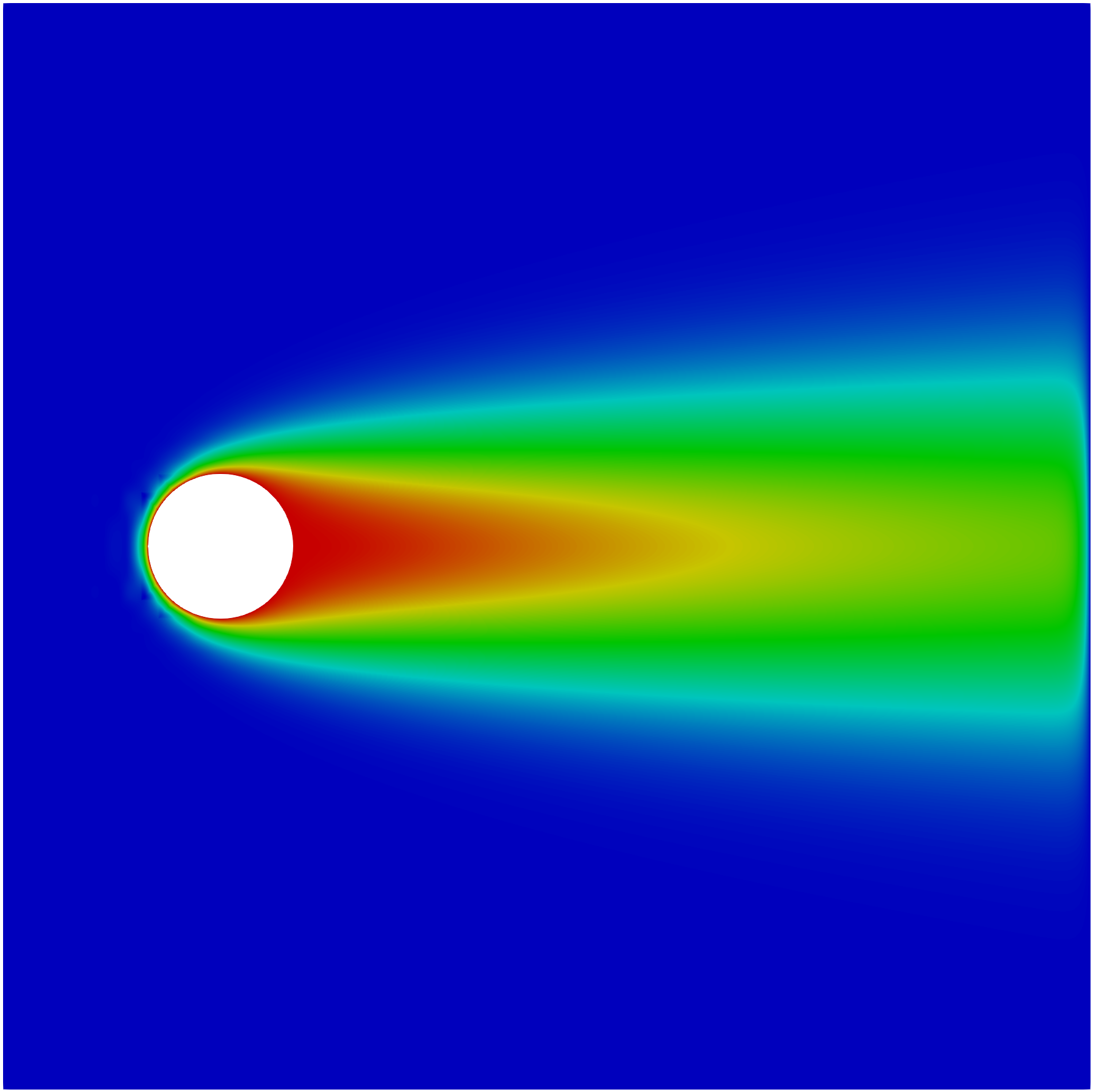}&\hspace{-2mm}
\includegraphics[width=0.47\columnwidth]{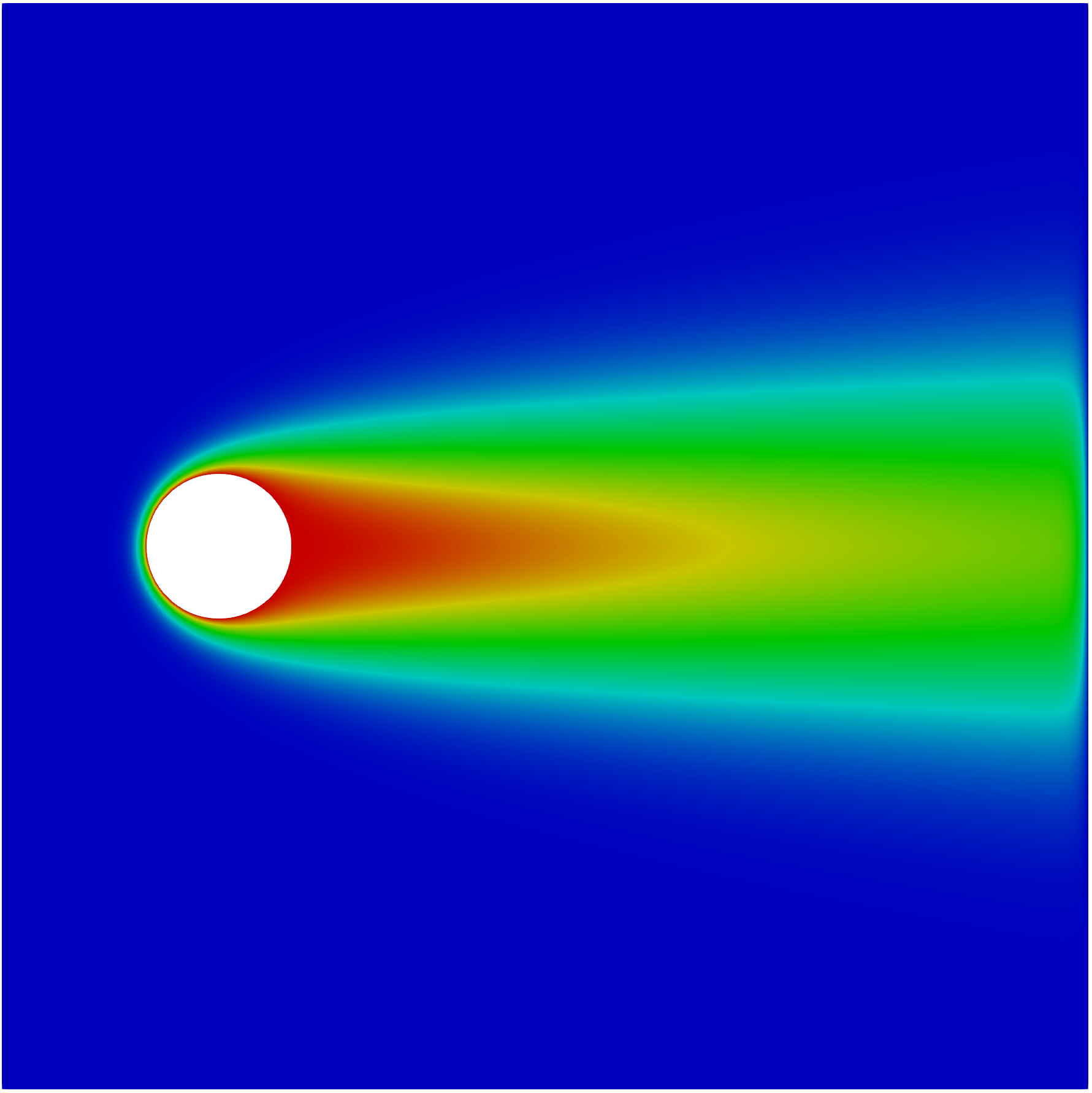}\\
(c) Mesh 2 (DOF: 7548) & (d) Mesh 3 (DOF: 29402)\\
\includegraphics[width=0.47\columnwidth]{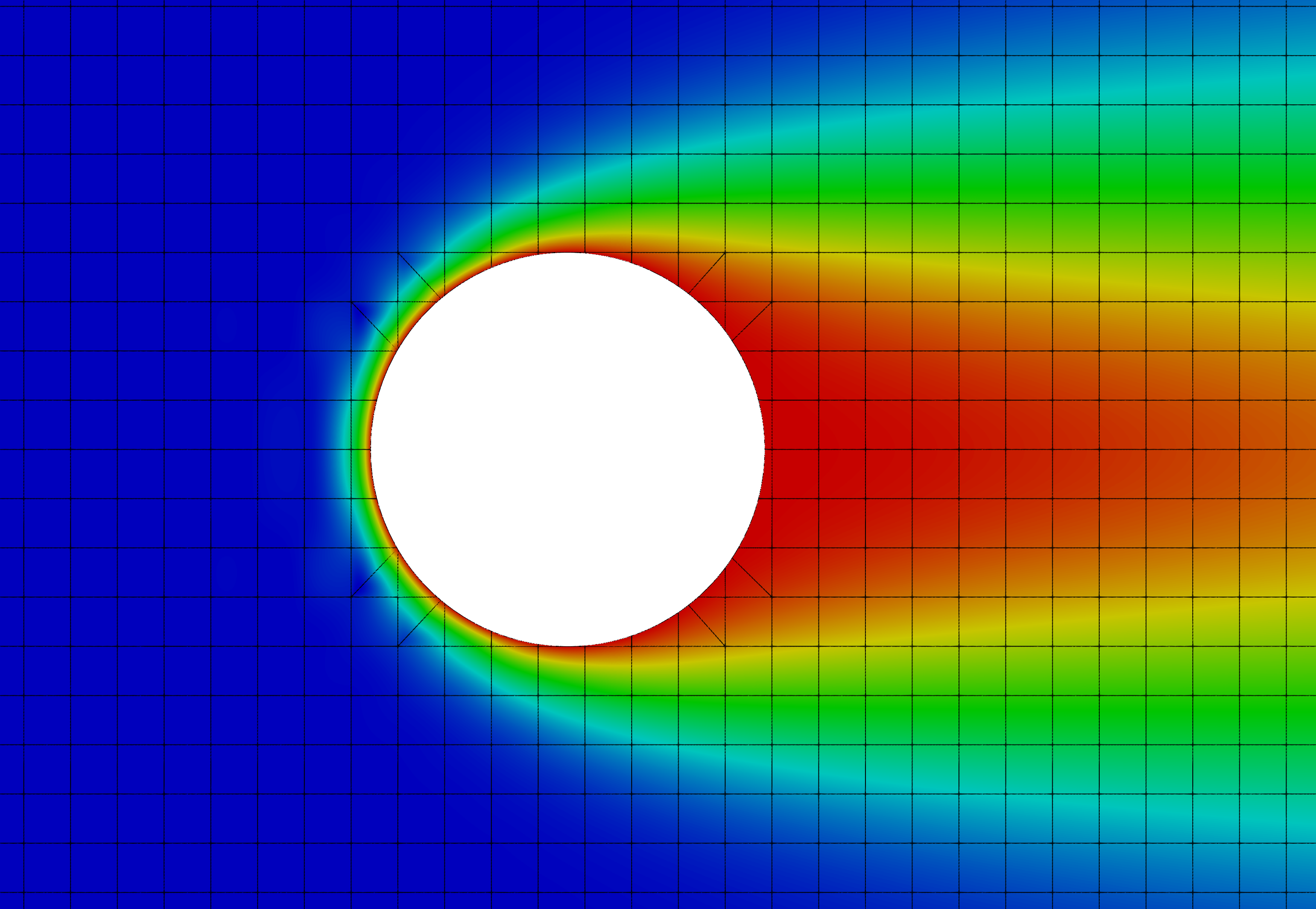}&\hspace{-2mm}
\includegraphics[width=0.47\columnwidth]{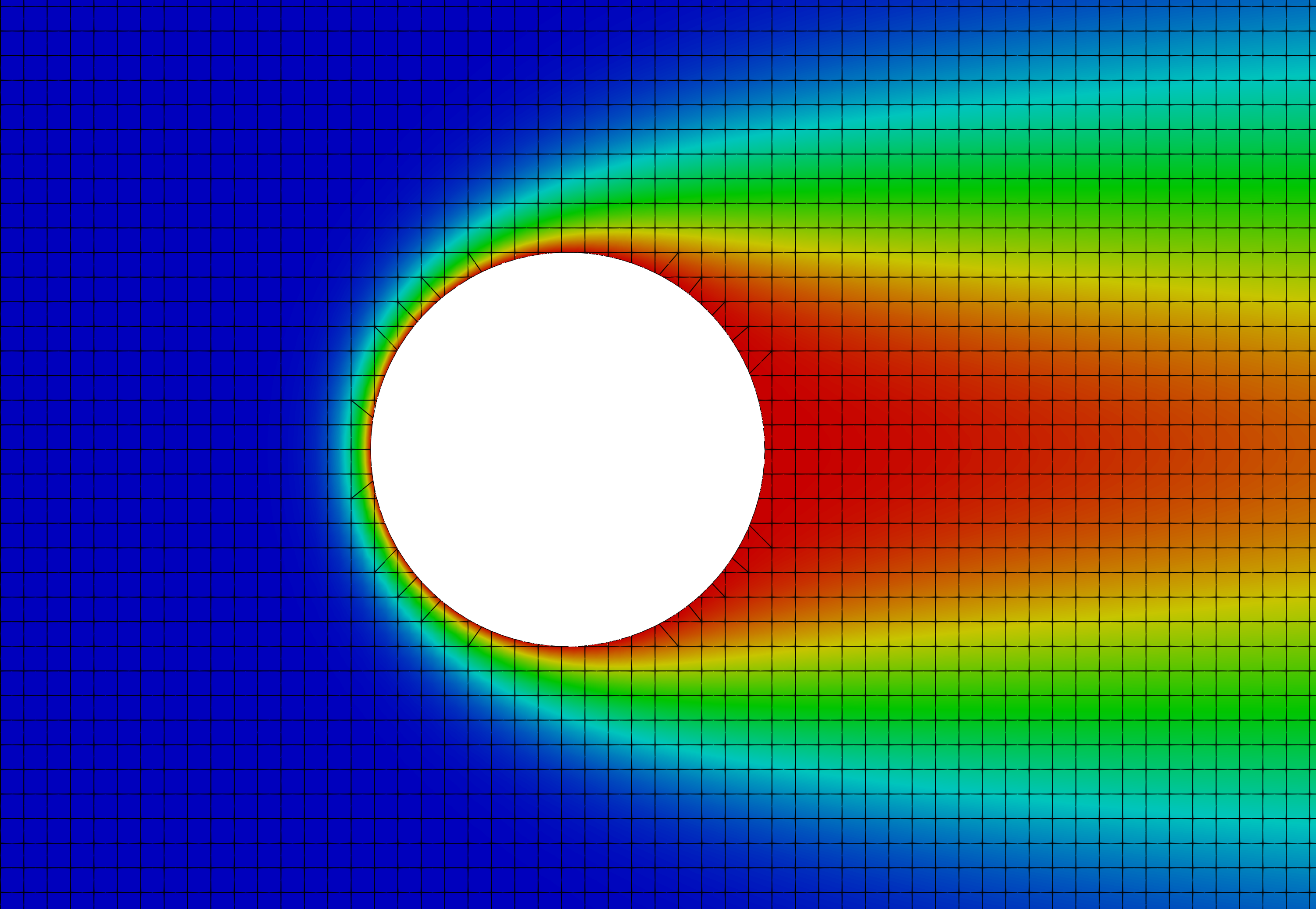}\\
(e) Zoom-in of Mesh 2 & (f) Zoom-in of Mesh 3  \\
\end{tabular}
\includegraphics[width=0.4\columnwidth]{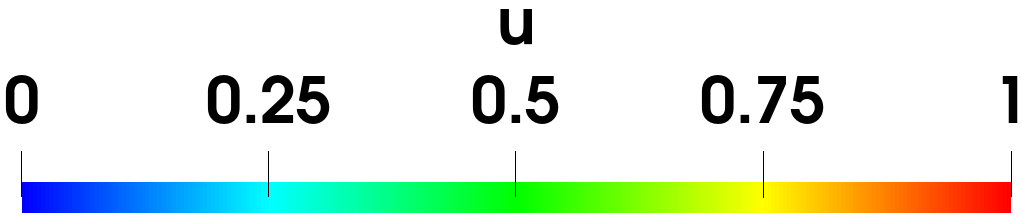}\\
\caption{Simulation results of the advection-diffusion problem using trimming. In the background meshes of (e, f), integration cells of trimmed elements are used for visualization.}
\label{fig:ad_trim_result}
\end{figure}

On other other hand, we can add a conformal layer around the hole to enhance the solution. The initial B\'ezier mesh is shown in Fig.~\ref{fig:ad_mesh}(b), where the thickness of the conformal layer is $1$. As the conformal layer aligns with the hole, we can easily control the mesh resolution around the hole. Even though the background mesh is the same as that in Fig.~\ref{fig:ad_mesh}(a), the solution is significantly improved; see Fig.~\ref{fig:ad_IBCM_result}(a). Moreover, with the background mesh globally refined just once and the conformal layer unchanged, the solution using IBCM is already comparable to the best result using trimming; see Figs.~\ref{fig:ad_IBCM_result}(b) and \ref{fig:ad_trim_result}(d). In other words, with IBCM, we can use much fewer DOF ($<1/10$) to resolve the boundary-layer phenomenon.

\begin{figure}[htb]
\centering
\begin{tabular}{cc}
\includegraphics[width=0.47\columnwidth]{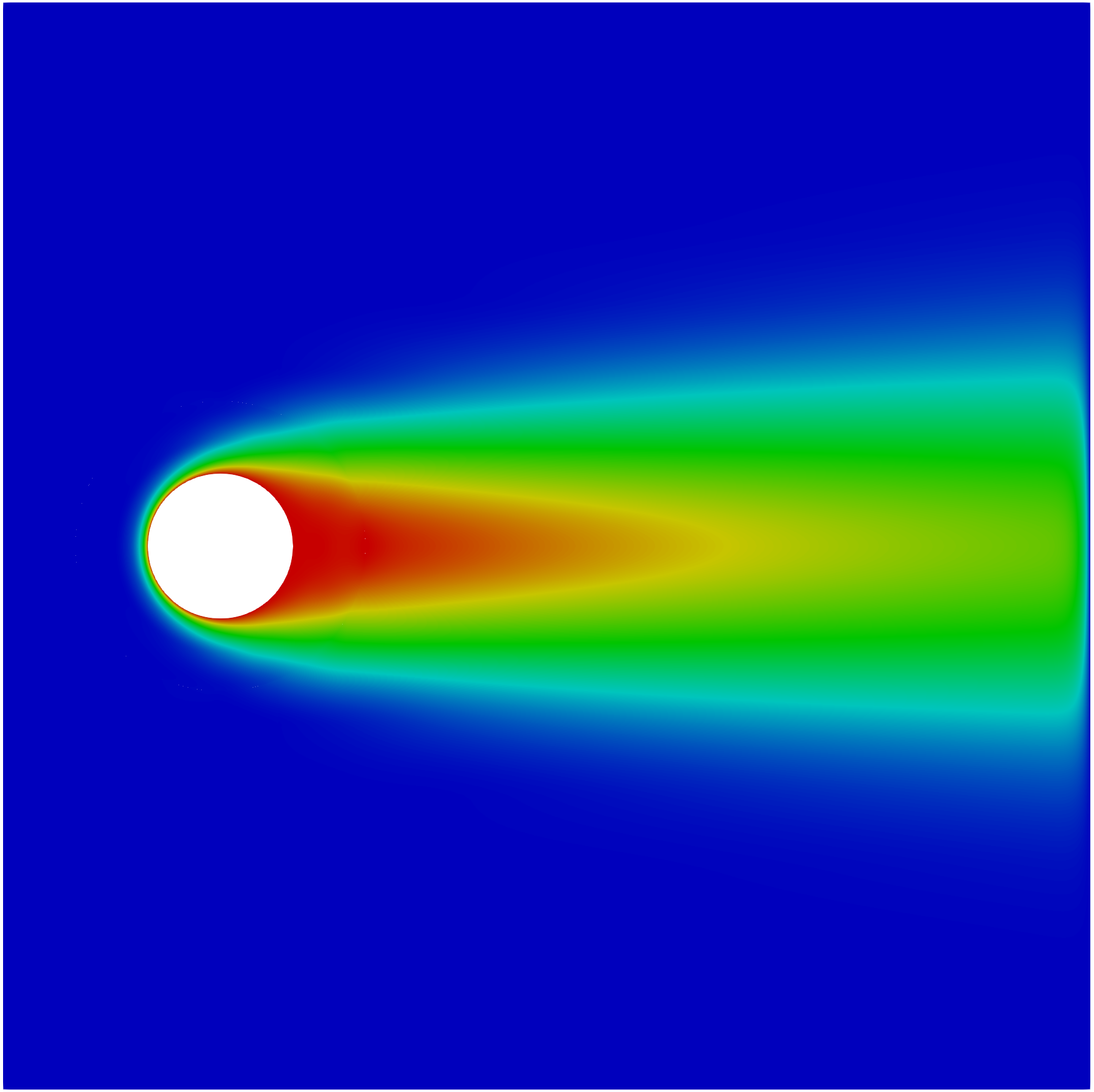}&\hspace{-2mm}
\includegraphics[width=0.47\columnwidth]{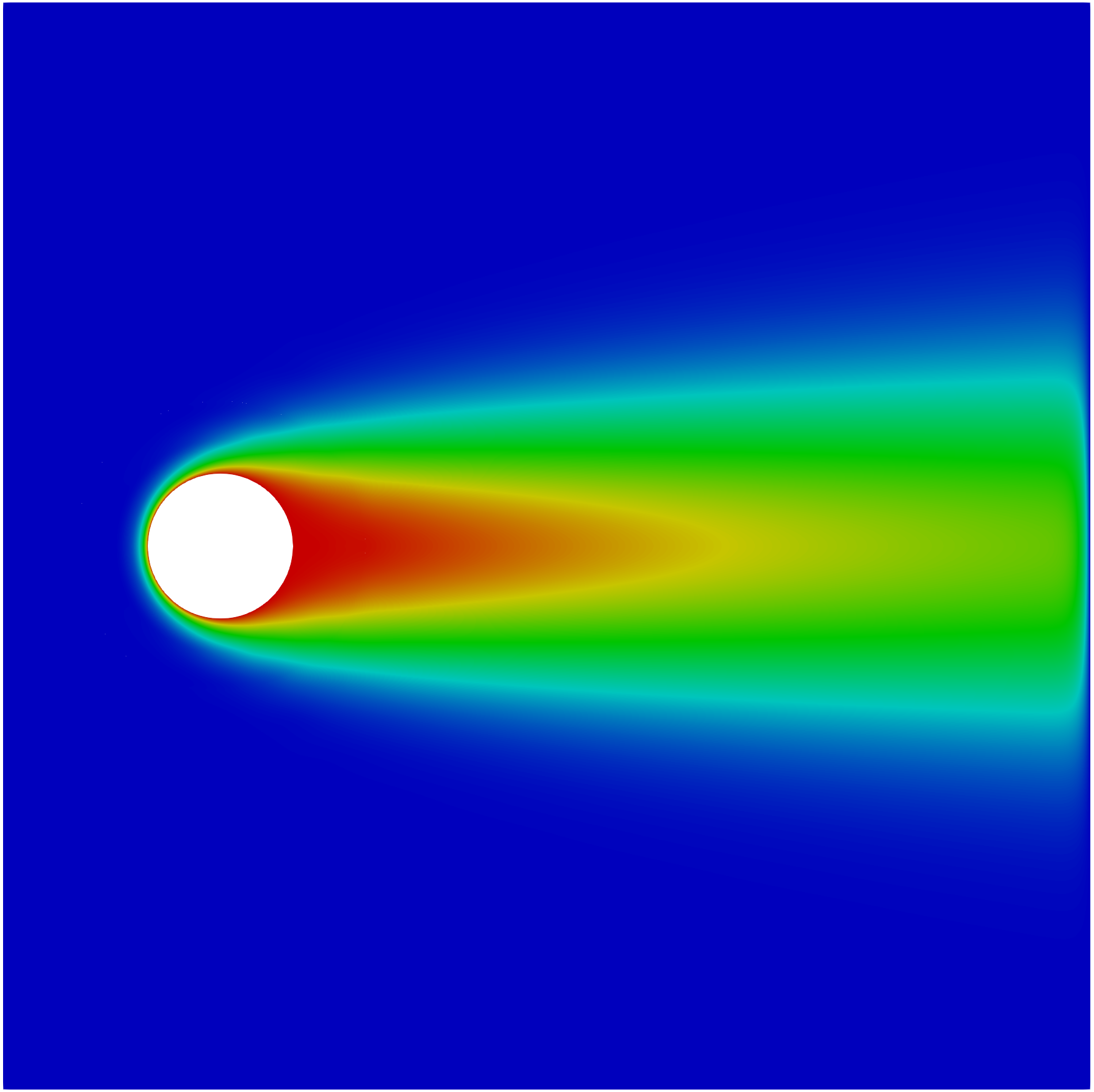}\\
(a) Mesh 0 (DOF: 694) & (b) Mesh 1 (DOF: 2120)\\
\includegraphics[width=0.47\columnwidth]{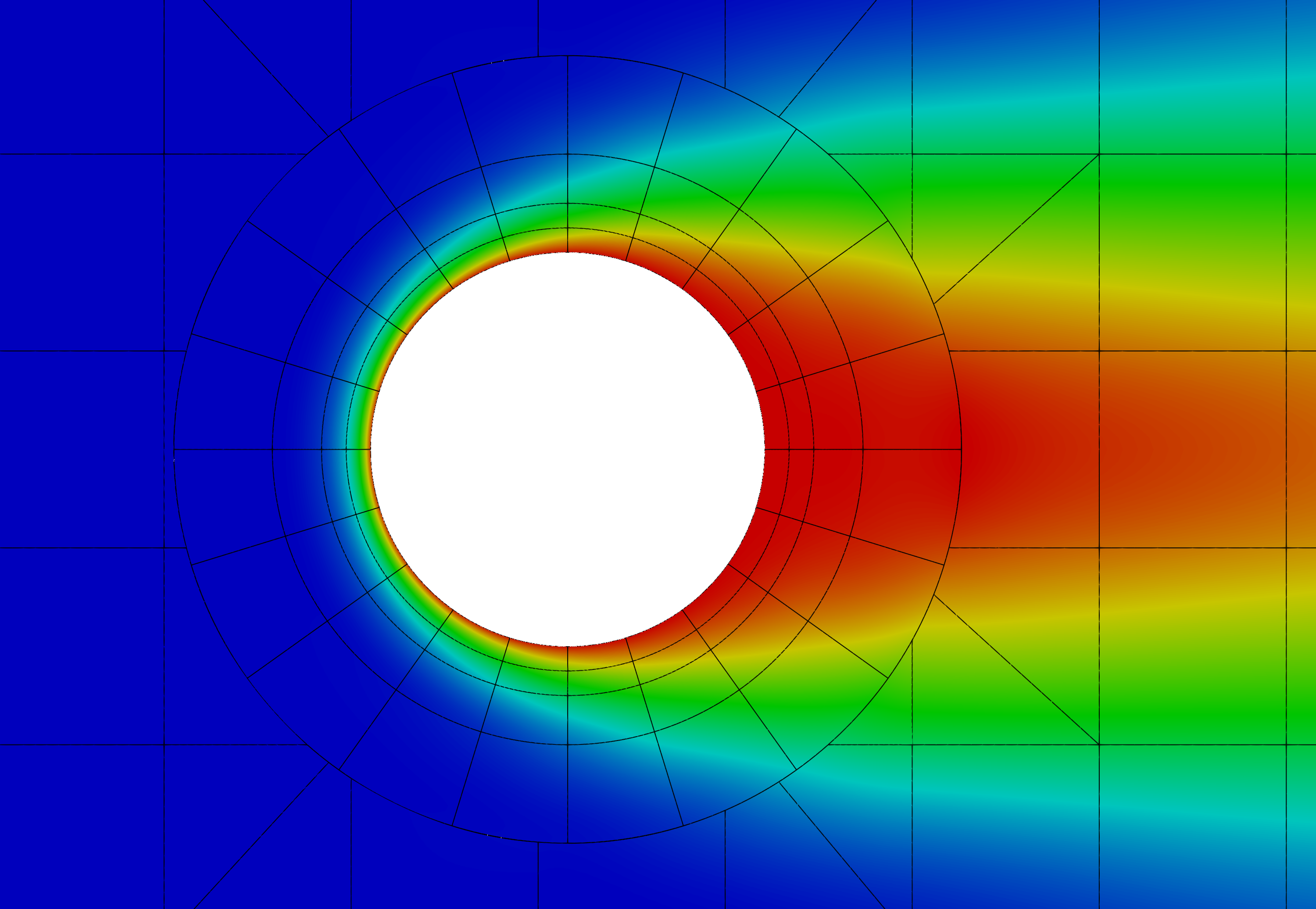}&\hspace{-2mm}
\includegraphics[width=0.47\columnwidth]{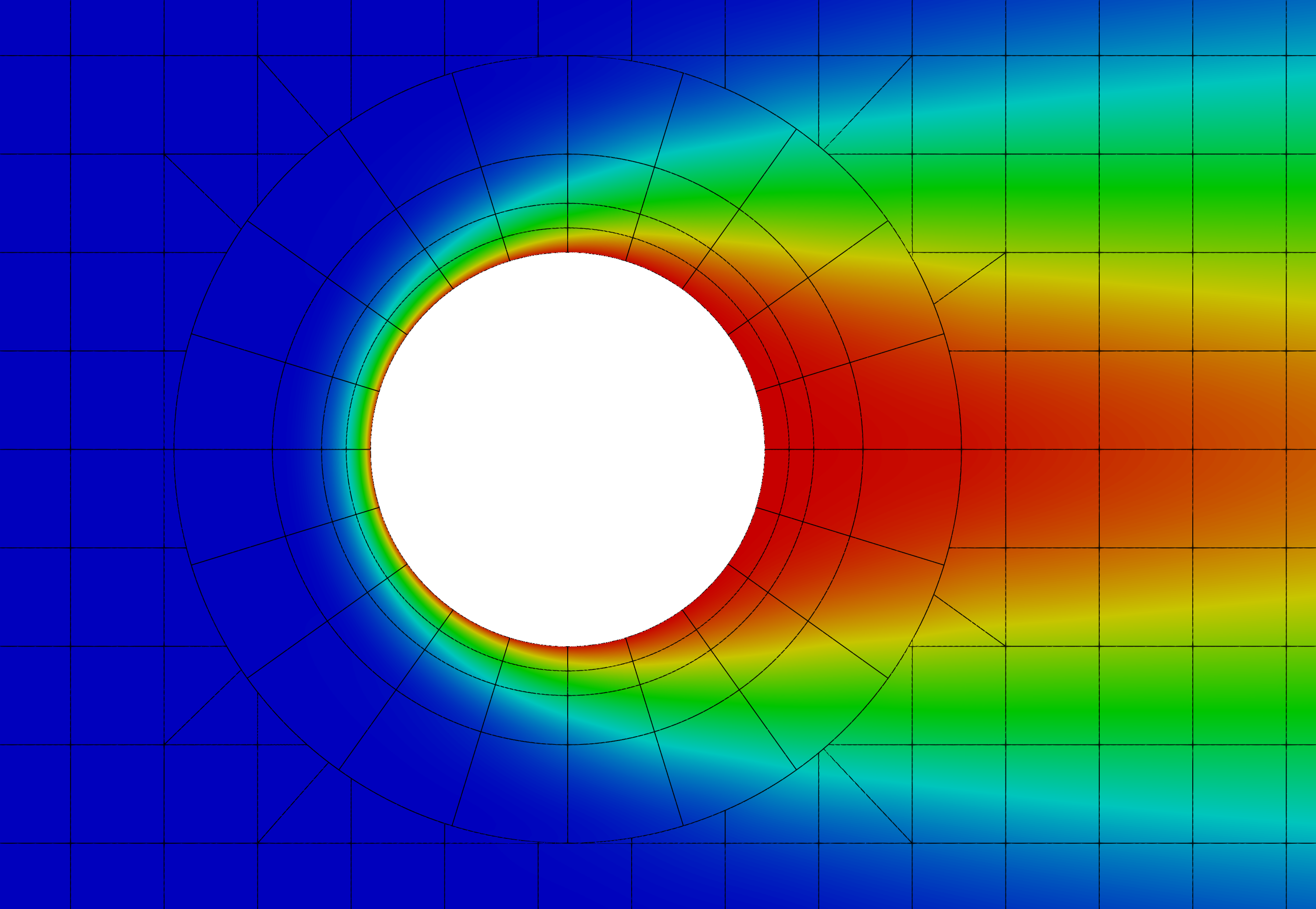}\\
(c) Zoom-in of Mesh 0 & (d) Zoom-in of Mesh 1  \\
\end{tabular}
\includegraphics[width=0.4\columnwidth]{figures/ad_bar.png}\\
\caption{Simulation results of the advection-diffusion problem using IBCM. In the background meshes of (c, d), integration cells of trimmed elements are used for visualization.}
\label{fig:ad_IBCM_result}
\end{figure}

Moreover, we plot the solution field along a horizontal line $y=H/2$; see Fig. \ref{fig:ad_plot_line}. Note that the region $[2,4]$ represents the hole. Fig. \ref{fig:ad_plot_line} shows three results: trimming with Mesh 1 (Fig. \ref{fig:ad_trim_result}(b)), trimming with Mesh 3 (Fig. \ref{fig:ad_trim_result}(d)), and IBCM with Mesh 1 (Fig. \ref{fig:ad_IBCM_result}(b)). We thus quantitatively confirm that IBCM with Mesh 1 achieves a very similar result to trimming with Mesh 3, whereas trimming with Mesh 1 exhibits oscillation near the hole as the mesh resolution is not fine enough.

\begin{figure}[htb]
\centering
\includegraphics[width=\columnwidth]{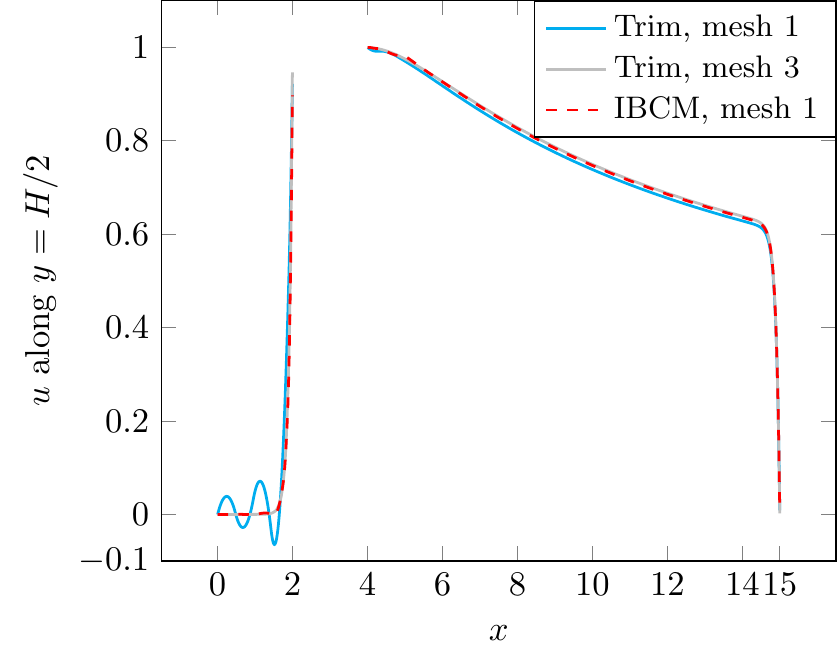}
\caption{The solution $u$ over the horizontal line $y=H/2$. The region $[2,4]$ represents the hole.}
\label{fig:ad_plot_line}
\end{figure}

}

\section{Conclusions and Future Work}
\label{sec:con}

In this paper, we have presented an immersed boundary-conformal method (IBCM) to leverage advantages of both conformal discretization and immersed methods. This is enabled by analysis-aware treatment of trimming and union operations as well as sophisticated construction of conformal layers. To efficiently and robustly deal with cut elements, we present an enhanced decomposition method to reduce the number of quadrature cells needed for a cut element. The union operation weakly couples independent domains through Nitsche's method. With the one-sided flux from the non-trimmed conformal layer, no further treatment is needed regarding stabilization.

On the other hand, the key to constructing a conformal layer lies largely in creating a target curve for an input boundary curve. The target curve can be constructed either as an offset curve or as a curve that has a very different shape from the input. Offset curves seem to be a natural choice, but very often they end up with cusps and self-intersections, and moreover, the parameterization is generally different from the input. Therefore, it requires further repairing and lofting to be able to construct a desired conformal layer. Fortunately, it is flexible for IBCM to choose the shape of the target curve, and there is numerical evidence that the shape does not play a significant role in the solution accuracy. With this flexibility, a target curve that has the same parameterization as the input can be constructed, and thus the conformal layer is readily available.

Two types of geometric constructions with IBCM are presented: a boundary type and an interface type. Both types can be used to represent complex geometric and/or solution features. In a boundary-type IBCM representation, as it features a conformal discretization near the boundary, many benefits are carried over from boundary-fitted methods, such as the intuitive control of mesh resolution and the ability to strongly impose Dirichlet boundary conditions. In an interface-type representation, kinematic constraints on the interface are automatically satisfied. Various benchmark problems are present to show the improved accuracy of IBCM compared to other means of geometric representations, such as solely trimming or union. There also exists evidence that the results are not sensitive to the thickness and shape of a conformal layer, but further study is needed to conclude. Moreover, two examples with complex geometric features are presented to show the flexibility of IBCM in representing complex geometries, strongly imposing Dirichlet boundary conditions, as well as easy adaptation of mesh resolution to solution features. \rt{In the end, we push a step forward to study an advection-diffusion problem with IBCM, where we have shown that IBCM can efficiently resolve the boundary-layer phenomenon near general geometric features in an intuitive manner.}

In the future, on top of all potential directions is the extension to 3D. This is indeed challenging as it requires a robust analysis-aware treatment for trimming and union. Constructing conformal layers for surfaces also becomes much more involving. On the other hand, extending IBCM to other problems is also promising, such as the Stokes problem, nonlinear problems, and shell/plate problems, where coming up with stabilized formulations poses as the most challenging problem.

\begin{acknowledgements}
X. Wei is partially supported by the ERC AdG project CHANGE n. 694515 and the SNSF (Swiss National Science Foundation) project HOGAEMS n.200021\_188589. P. Antolin, and A. Buffa are partially supported by the ERC AdG project CHANGE n. 694515 and the SNSF through the project ``Design-through-Analysis (of PDEs): the litmus test” n.40B2-0 187094 (BRIDGE Discovery 2019).
\end{acknowledgements}

\bibliographystyle{unsrt}
\bibliography{ref}

\begin{thebibliography}{10}

\bibitem{ref:wei19u}
P.~Antolin, A.~Buffa, R.~Puppi, and X.~Wei.
\newblock Overlapping multi-patch isogeometric method with minimal
  stabilization.
\newblock {\em SIAM Journal on Scientific Computing}, 43(1):A330--A354, 2021.

\bibitem{ref:hughes05}
T.~J.~R. Hughes, J.~A. Cottrell, and Y.~Bazilevs.
\newblock Isogeometric analysis: {CAD}, finite elements, {NURBS}, exact
  geometry and mesh refinement.
\newblock {\em Computer Methods in Applied Mechanics and Engineering},
  194(39):4135--4195, 2005.

\bibitem{ref:cottrell09}
J.~A. Cottrell, T.~J.~R. Hughes, and Y.~Bazilevs.
\newblock {\em Isogeometric Analysis: Toward Integration of CAD and FEA}.
\newblock Wiley Publishing, 2009.

\bibitem{ref:martin09}
T.~Martin, E.~Cohen, and R.~M. Kirby.
\newblock Volumetric parameterization and trivariate {B}-spline fitting using
  harmonic functions.
\newblock {\em Computer Aided Geometric Design}, 26(6):648--664, 2009.

\bibitem{ref:zhang13}
Y.~Zhang, W.~Wang, and T.~J.~R. Hughes.
\newblock Conformal solid {T}-spline construction from boundary {T}-spline
  representations.
\newblock {\em Computational Mechanics}, 51(6):1051--1059, 2013.

\bibitem{ref:rank12}
E.~Rank, M.~Ruess, S.~Kollmannsberger, D.~Schillinger, and A.~D\text{\"u}ster.
\newblock Geometric modeling, isogeometric analysis and the finite cell method.
\newblock {\em Computer Methods in Applied Mechanics and Engineering},
  249-252:104--115, 2012.

\bibitem{ref:ruberg12}
T.~Rüberg and F.~Cirak.
\newblock Subdivision-stabilised immersed {B}-spline finite elements for moving
  boundary flows.
\newblock {\em Computer Methods in Applied Mechanics and Engineering},
  209-212:266--283, 2012.

\bibitem{ref:schillinger12}
D.~Schillinger, L.~Ded\`e, M.~A. Scott, J.~A. Evans, M.~J. Borden, E.~Rank, and
  T.J.~R. Hughes.
\newblock An isogeometric design-through-analysis methodology based on adaptive
  hierarchical refinement of {NURBS}, immersed boundary methods, and {T}-spline
  {CAD} surfaces.
\newblock {\em Computer Methods in Applied Mechanics and Engineering},
  249-252:116--150, 2012.

\bibitem{ref:kamensky15}
D.~Kamensky, M.-C. Hsu, D.~Schillinger, J.~A. Evans, A.~Aggarwal, Y.~Bazilevs,
  M.~S. Sacks, and T.~J.~R. Hughes.
\newblock An immersogeometric variational framework for fluid–structure
  interaction: {A}pplication to bioprosthetic heart valves.
\newblock {\em Computer Methods in Applied Mechanics and Engineering},
  284:1005--1053, 2015.

\bibitem{ref:hoang19}
T.~Hoang, C.~V. Verhoosel, C.-Z. Qin, F.~Auricchio, A.~Reali, and E.~H. van
  Brummelen.
\newblock Skeleton-stabilized immersogeometric analysis for incompressible
  viscous flow problems.
\newblock {\em Computer Methods in Applied Mechanics and Engineering},
  344:421--450, 2019.

\bibitem{ref:casquero20im}
H.~Casquero, C.~Bona-Casas, D.~Toshniwal, T.~J.~R. Hughes, H.~Gomez, and Y.~J.
  Zhang.
\newblock The divergence-conforming immersed boundary method: {A}pplication to
  vesicle and capsule dynamics.
\newblock arXiv:2001.08244, 2020.

\bibitem{ref:zli06}
Z.~Li and K.~Ito.
\newblock {\em The Immersed Interface Method: {N}umerical Solutions of {PDE}s
  Involving Interfaces and Irregular Domains}.
\newblock Society for Industrial and Applied Mathematics, 2006.

\bibitem{ref:kim09}
H.~Kim, Y.~Seo, and S.~Youn.
\newblock Isogeometric analysis for trimmed cad surfaces.
\newblock {\em Computer Methods in Applied Mechanics and Engineering},
  198(37):2982--2995, 2009.

\bibitem{ref:nagy15}
A.~P. Nagy and D.~J. Benson.
\newblock On the numerical integration of trimmed isogeometric elements.
\newblock {\em Computer Methods in Applied Mechanics and Engineering},
  284:165--185, 2015.

\bibitem{ref:kudela15}
L.~Kudela, N.~Zander, T.~Bog, S.~Kollmannsberger, and E.~Rank.
\newblock Efficient and accurate numerical quadrature for immersed boundary
  methods.
\newblock {\em Advanced Modeling and Simulation in Engineering Sciences}, 2,
  2015.

\bibitem{ref:antolin19a}
P.~Antolin, A.~Buffa, and M.~Martinelli.
\newblock Isogeometric analysis on {V}-reps: {F}irst results.
\newblock {\em Computer Methods in Applied Mechanics and Engineering},
  355:976--1002, 2019.

\bibitem{ref:massarwi19}
F.~Massarwi, P.~Antolin, and G.~Elber.
\newblock Volumetric untrimming: {P}recise decomposition of trimmed trivariates
  into tensor products.
\newblock {\em Computer Aided Geometric Design}, 71:1--15, 2019.

\bibitem{ref:marussig17}
B.~Marussig, J.~Zechner, G.~Beer, and T.~P. Fries.
\newblock Stable isogeometric analysis of trimmed geometries.
\newblock {\em Computer Methods in Applied Mechanics and Engineering},
  316:497--521, 2017.

\bibitem{ref:marussig18}
B.~Marussig, R.~Hiemstra, and T.~J.~R. Hughes.
\newblock Improved conditioning of isogeometric analysis matrices for trimmed
  geometries.
\newblock {\em Computer Methods in Applied Mechanics and Engineering},
  334:79--110, 2018.

\bibitem{ref:guo18}
Y.~Guo, J.~Heller, T.~J.~R. Hughes, M.~Ruess, and D.~Schillinger.
\newblock Variationally consistent isogeometric analysis of trimmed thin shells
  at finite deformations, based on the {STEP} exchange format.
\newblock {\em Computer Methods in Applied Mechanics and Engineering},
  336:39--79, 2018.

\bibitem{ref:puppi19}
A.~Buffa, R.~Puppi, and R.~V\'{a}zquez.
\newblock A minimal stabilization procedure for isogeometric methods on trimmed
  geometries.
\newblock {\em SIAM Journal on Numerical Analysis}, 58(5):2711–--2735, 2020.

\bibitem{ref:bazilevs07}
Y.~Bazilevs and T.~J.~R. Hughes.
\newblock Weak imposition of {D}irichlet boundary conditions in fluid
  mechanics.
\newblock {\em Computers \& Fluids}, 36(1):12--26, 2007.

\bibitem{ref:embar10}
A.~Embar, J.~Dolbow, and I.~Harari.
\newblock Imposing {D}irichlet boundary conditions with {N}itsche's method and
  spline-based finite elements.
\newblock {\em International Journal for Numerical Methods in Engineering},
  83(7):877--898, 2010.

\bibitem{ref:benek83}
J.~A. Benek, J.~L. Steger, and F.~C. Dougherty.
\newblock A flexible grid embedding technique with application to the {E}uler
  equations.
\newblock {\em AIAA Journal}, pages 373--382, 1983.

\bibitem{ref:liou94}
M.-S. Liou and K.-H. Kao.
\newblock Progress in grid generation: {F}rom {C}himera to {DRAGON} grids.
\newblock NASA Technical Memorandum 106709, 1994.

\bibitem{ref:henshaw94}
W.~D. Henshaw.
\newblock A fourth-order accurate method for the incompressible
  {N}avier-{S}tokes equations on overlapping grids.
\newblock {\em Journal of Computational Physics}, 113(1):13--25, 1994.

\bibitem{ref:tang03}
H.~S. Tang, S.~C. Jones, and F.~Sotiropoulos.
\newblock An overset-grid method for {3D} unsteady incompressible flows.
\newblock {\em Journal of Computational Physics}, 191(2):567–--600, 2003.

\bibitem{ref:kannan07}
R.~Kannan and Z.~J. Wang.
\newblock Overset adaptive {C}artesian/prism grid method for stationary and
  moving-boundary flow problems.
\newblock {\em AIAA Journal}, 45(7):1774--1779, 2007.

\bibitem{ref:duan20}
Z.~Duan and Z.~J. Wang.
\newblock High-order overset flux reconstruction method for dynamic moving
  grids.
\newblock {\em AIAA Journal}, 58(10):4534--4547, 2020.

\bibitem{ref:bouclier16}
R.~Bouclier, J.~Passieux, and M.~Salaün.
\newblock Local enrichment of {NURBS} patches using a non-intrusive coupling
  strategy: {G}eometric details, local refinement, inclusion, fracture.
\newblock {\em Computer Methods in Applied Mechanics and Engineering},
  300:1--26, 2016.

\bibitem{ref:dokken19}
J.~S. Dokken, S.~W. Funke, A.~Johansson, and S.~Schmidt.
\newblock Shape optimization using the finite element method on multiple meshes
  with {N}itsche coupling.
\newblock {\em SIAM Journal on Scientific Computing}, 41(3):A1923--A1948, 2019.

\bibitem{ref:becker11}
R.~Becker, E.~Burman, and P.~Hansbo.
\newblock A hierarchical {NXFEM} for fictitious domain simulations.
\newblock {\em International Journal for Numerical Methods in Engineering},
  86(4-5):549--559, 2011.

\bibitem{ref:johansson19}
A.~Johansson, B.~Kehlet, M.~G. Larson, and A.~Logg.
\newblock Multimesh finite element methods: {S}olving {PDE}s on multiple
  intersecting meshes.
\newblock {\em Computer Methods in Applied Mechanics and Engineering},
  343:672--689, 2019.

\bibitem{ref:nurbsbook}
L.~Piegl and W.~Tiller.
\newblock {\em The NURBS Book (2nd Ed.)}.
\newblock Springer-Verlag, 1997.

\bibitem{ref:marussig18r}
B.~Marussig and T.~J.~R. Hughes.
\newblock A review of trimming in isogeometric analysis: {C}hallenges, data
  exchange and simulation aspects.
\newblock {\em Archives of Computational Methods in Engineering},
  25:1059–--1127, 2018.

\bibitem{ref:mc}
W.~E. Lorensen and H.~E. Cline.
\newblock Marching cubes: {A} high resolution {3D} surface construction
  algorithm.
\newblock {\em ACM SIGGRAPH Computer Graphics}, 21(4):163–--169, 1987.

\bibitem{ref:fries17a}
T.~P. Fries, S.~Omerović, D.~Schöllhammer, and J.~Steidl.
\newblock {Higher-order meshing of implicit geometries—Part I: Integration
  and interpolation in cut elements}.
\newblock {\em Computer Methods in Applied Mechanics and Engineering},
  313:759--784, 2017.

\bibitem{ref:fries17b}
T.~P. Fries and D.~Schöllhammer.
\newblock {Higher-order meshing of implicit geometries, Part II: Approximations
  on manifolds}.
\newblock {\em Computer Methods in Applied Mechanics and Engineering},
  326:270--297, 2017.

\bibitem{ref:nitsche71}
J.~Nitsche.
\newblock Über ein variationsprinzip zur lösung von dirichlet-problemen bei
  verwendung von teilräumen, die keinen randbedingungen unterworfen sind.
\newblock {\em Abhandlungen aus dem Mathematischen Seminar der Universität
  Hamburg}, 36:9--15, 1971.

\bibitem{ref:prenter17}
F.~de~Prenter, C.~V. Verhoosel, G.~J. van Zwieten, and E.~H. van Brummelen.
\newblock Condition number analysis and preconditioning of the finite cell
  method.
\newblock {\em Computer Methods in Applied Mechanics and Engineering},
  316:297--327, 2017.

\bibitem{ref:hansbo05}
P.~Hansbo.
\newblock Nitsche's method for interface problems in computational mechanics.
\newblock {\em Gamm-mitteilungen}, 28:183--206, 2005.

\bibitem{Hoschek1992b}
J.~Hoschek and D.~Lasser.
\newblock {\em Grundlagen der geometrischen Datenverarbeitung}.
\newblock Vieweg+Teubner, 1992.

\bibitem{patrikalakis2009b}
N.~M. Patrikalakis and T.~Maekawa.
\newblock {\em Shape Interrogation for Computer Aided Design and
  Manufacturing}.
\newblock Springer Science \& Business Media, 2009.

\bibitem{Maekawa1999a}
T.~Maekawa.
\newblock An overview of offset curves and surfaces.
\newblock {\em Computer-Aided Design}, 31(3):165--173, 1999.

\bibitem{Pham1992a}
B.~Pham.
\newblock Offset curves and surfaces: {A} brief survey.
\newblock {\em Computer-Aided Design}, 24(4):223--229, 1992.

\bibitem{WALLNER2001a}
J.~Wallner, T.~Sakkalis, T.~Maekawa, H.~Pottmann, and G.~Yu.
\newblock Self-intersections of offset curves and surfaces.
\newblock {\em International Journal of Shape Modeling}, 07(01):1--21, 2001.

\bibitem{ref:rhino}
{Rhino - Rhinoceros 3D}.
\newblock https://www.rhino3d.com, 2021.

\bibitem{ref:hinz18}
J.~Hinz, M.~Möller, and C.~Vuik.
\newblock Elliptic grid generation techniques in the framework of isogeometric
  analysis applications.
\newblock {\em Computer Aided Geometric Design}, 65:48--75, 2018.

\bibitem{Randrianarivony2006phd}
Maharavo Randrianarivony.
\newblock {\em Geometric Processing of {CAD} Data And Meshes as Input of
  Integral Equation Solvers}.
\newblock PhD thesis, Computer Science Faculty Technische Universit{\"a}t
  Chemnitz, 2006.

\bibitem{Bommes2013a}
David Bommes, Marcel Campen, Hans-Christian Ebke, Pierre Alliez, and Leif
  Kobbelt.
\newblock Integer-grid maps for reliable quad meshing.
\newblock {\em ACM Trans. Graph.}, 32(4):98:1--98:12, July 2013.

\bibitem{ref:sederberg03}
T.~W. Sederberg, J.~Zheng, A.~Bakenov, and A.~Nasri.
\newblock T-splines and {T-NURCC}s.
\newblock {\em ACM Trans. Graph.}, 22(3):477--484, 2003.

\bibitem{ref:veiga12}
L.~Beir\text{\~a}o da~Veiga, A.~Buffa, D.~Cho, and G.~Sangalli.
\newblock Analysis-suitable {T}-splines are dual-compatible.
\newblock {\em Computer Methods in Applied Mechanics and Engineering},
  249-252:42--51, 2012.

\bibitem{ref:xli18}
J.~Zhang and X.~Li.
\newblock Local refinement for analysis-suitable++ {T}-splines.
\newblock {\em Computer Methods in Applied Mechanics and Engineering},
  342:32--45, 2018.

\bibitem{ref:casquero20}
H.~Casquero, X.~Wei, D.~Toshniwal, A.~Li, T.~J.~R. Hughes, J.~Kiendl, and Y.~J.
  Zhang.
\newblock Seamless integration of design and {K}irchhoff–{L}ove shell
  analysis using analysis-suitable unstructured {T}-splines.
\newblock {\em Computer Methods in Applied Mechanics and Engineering},
  360:112765, 2020.

\bibitem{ref:vuong11}
A.-V. Vuong, C.~Giannelli, B.~J\text{\"u}ttler, and B.~Simeon.
\newblock A hierarchical approach to adaptive local refinement in isogeometric
  analysis.
\newblock {\em Computer Methods in Applied Mechanics and Engineering},
  200(49):3554--3567, 2011.

\bibitem{ref:giannelli12}
C.~Giannelli, B.~J\text{\"u}ttler, and H.~Speleers.
\newblock {THB}-splines: {T}he truncated basis for hierarchical splines.
\newblock {\em Computer Aided Geometric Design}, 29(7):485--498, 2012.

\bibitem{ref:wei15}
X.~Wei, Y.~Zhang, T.~J.~R. Hughes, and M.~A. Scott.
\newblock Truncated hierarchical {C}atmull--{C}lark subdivision with local
  refinement.
\newblock {\em Computer Methods in Applied Mechanics and Engineering},
  291:1--20, 2015.

\bibitem{ref:buffa16}
A.~Buffa and C.~Giannelli.
\newblock Adaptive isogeometric methods with hierarchical splines: {E}rror
  estimator and convergence.
\newblock {\em Mathematical Models and Methods in Applied Sciences},
  26(1):1--25, 2016.

\bibitem{ref:bracco19}
C.~Bracco, A.~Buffa, C.~Giannelli, and R.~V\'{a}zquez.
\newblock Adaptive isogeometric methods with hierarchical splines: {A}n
  overview.
\newblock {\em Discrete \& Continuous Dynamical Systems - A}, 39(1):241--261,
  2019.

\bibitem{ref:johannessen14}
K.~A. Johannessen, T.~Kvamsdal, and T.~Dokken.
\newblock Isogeometric analysis using lr b-splines.
\newblock {\em Computer Methods in Applied Mechanics and Engineering},
  269:471--514, 2014.

\bibitem{ref:patrizi20}
F.~Patrizi, C.~Manni, F.~Pelosi, and H.~Speleers.
\newblock Adaptive refinement with locally linearly independent lr b-splines:
  Theory and applications.
\newblock {\em Computer Methods in Applied Mechanics and Engineering},
  369:113230, 2020.

\bibitem{ref:sukumar01}
N.~Sukumar, D.L. Chopp, N.~Moës, and T.~Belytschko.
\newblock Modeling holes and inclusions by level sets in the extended
  finite-element method.
\newblock {\em Computer Methods in Applied Mechanics and Engineering},
  190(46):6183--6200, 2001.

\bibitem{ref:stanford19}
J.~W. Stanford and T.~P. Fries.
\newblock A higher-order conformal decomposition finite element method for
  plane {B}-rep geometries.
\newblock {\em Computers \& Structures}, 214:15--27, 2019.

\bibitem{ref:herraez16}
M.~Herráez, C.~González, C.S. Lopes, R.~Guzmán {de Villoria}, J.~LLorca,
  T.~Varela, and J.~Sánchez.
\newblock Computational micromechanics evaluation of the effect of fibre shape
  on the transverse strength of unidirectional composites: An approach to
  virtual materials design.
\newblock {\em Composites Part A: Applied Science and Manufacturing},
  91:484--492, 2016.

\bibitem{ref:hemker96}
P.~W. Hemker.
\newblock A singularly perturbed model problem for numerical computation.
\newblock {\em Journal of Computational and Applied Mathematics},
  76(1):277--285, 1996.

\bibitem{ref:buffa06}
A.~Buffa, T.~J.~R. Hughes, and G.~Sangalli.
\newblock Analysis of a multiscale discontinuous galerkin method for
  convection‐diffusion problems.
\newblock {\em SIAM Journal on Numerical Analysis}, 44(4):1420–--1440, 2006.

\end{thebibliography}

\end{document}